\documentclass[11pt,a4paper]{amsart}

\usepackage{amsmath,amssymb,hyperref}
\usepackage{amsfonts}
\usepackage{enumerate}
\usepackage{amsmath, amsthm}
\usepackage{amscd}
\usepackage{tikz}
\usepackage{amssymb}
\hypersetup{pdfpagemode=FullScreen, colorlinks=true}

\newtheorem{thm}{Theorem}[section]

\newtheorem{prop}[thm]{Proposition}
\newtheorem{lemma}[thm]{Lemma}

\newtheorem{cor}[thm]{Corollary}

\theoremstyle{definition}
\newtheorem{defi}[thm]{Definition}
\newtheorem*{remark}{Remark}
\newtheorem*{example}{Example}

\newtheorem{remrk}[thm]{Remark}
\newtheorem{expl}[thm]{Example}
\newtheorem{construction}[thm]{Construction}
\numberwithin{equation}{section}

\newcommand{\R}{{\mathbb R}}
\newcommand{\N}{{\mathbb N}}
\newcommand{\C}{\mathcal C}

\newcommand{\ZZ}{{\mathbb Z}}
\newcommand{\Q}{{\mathbb Q}}
\newcommand{\J}{{\mathbb J}}
\newcommand{\IH}{{\mathbb H}}
\newcommand{\T}{{\mathcal T}}

\def\neuzl{\hfill\break}
\let\8=\ss
\let\eps=\epsilon
\let\Ups=\Upsilon
\let\Omg=\Omega
\let\omg=\omega
\let\noi=\noindent
\def\blackbox{\vrule width 1ex height 1ex depth 0pt}%
\def\Vol{{\it Vol}}%
\let\ttt=\tt\def\tt{\ttt\rightskip=0pt plus 0.88\hsize}
\def\hl{}

\newenvironment{pf}{\begin{trivlist}\item[]{\bf Proof:\ }}
{\mbox{}\hfill\rule{.08in}{.08in}\end{trivlist}}

\begin{document}
 \author{Thilo Kuessner}
   \address{Mathematisch-Geographische Fakult\"at,
Katholische Universit\"at Eichst\"att-Ingolstadt,
Ostenstra\ss e 28,
85072 Eichst\"att, Germany}
   \email{Thilo.Kuessner@ku.de}
 \author{Janusz Przewocki}
\address{
Institute of Mathematics,
Faculty of Mathematics, Physics and Informatics,
University of Gdansk,
ul. Wita Stwosza 57,
80-308 Gdansk,
Poland }
\email{j.przewocki@gmail.com}
\author{Andreas Zastrow}
\address{
Institute of Mathematics, 
Faculty of Mathematics, Physics and Informatics, 
University of Gdansk, 
ul. Wita Stwosza 57, 
80-308 Gdansk, 
Poland }
   \email{andreas.zastrow@ug.edu.pl}

\title{On measure homology of mildly wild spaces}

\begin{abstract}
We prove injectivity of the canonical map from singular homology to measure homology for certain ``mildly wild" spaces, that is, certain spaces not having the homotopy type of a CW-complex, but having countable fundamental groups.
\end{abstract}
\maketitle

Measure homology $\mathcal{H}_*(X)$, also called Milnor-Thurston homology, of a space $X$ is a variant of the usually studied singular homology groups $H_*(X;\R)$. While the latter are defined as the homology theory of the chain complex of finite linear combinations of singular simplices with the canonical boundary operator, measure homology uses the chain complex of \hl{compactly} determined signed measures of bounded variation on the space of singular simplices $map(\Delta^k,X)$, again with the canonical boundary operator (see \hyperref[mhdefs]{Section \ref*{mhdefs}}). 

Singular chains (with real coefficients) can be considered as finite sums of (real multiples of) Dirac measures, so there is a canonical homomorphism
$$\iota_*\colon H_*(X;\R)\to \mathcal{H}_*(X).$$
It was proven in \cite{zas} and \cite{han} that measure homology satisfies the Eilenberg-Steenrod axioms and thus that $\iota_*$ is an isomorphism whenever $X$ is a CW-complex. This \hl{applies in particular to closed manifolds and interiors of compact manifolds with boundary, because they always have the homotopy type of a CW-complex}.
 
\hl{The main motivation to consider measure homology is its use in the proof of proportionality principles. By ``proportionality principle" one means the result that the equality $$\frac{\Vol(M)}{\Vert M\Vert}=\frac{\Vol(N)}{\Vert N\Vert}$$ holds for Riemannian manifolds with isometric universal covering. Here $\Vol(M)$ denotes the Riemannian volume and $\Vert M\Vert$ is the simplicial volume, a homotopy invariant  defined as the $l^1$-norm (``Gromov norm") of the fundamental class (cf.\ the beginning of
\hyperref[gromovnorm]{Subsec.\ \ref*{gromovnorm}}). The proportionality principle for closed Riemannian manifolds was established by Gromov in \cite[Section 2.3]{gro}}. 

\hl{The question whether the proportionality principle holds in greater generality for finite-volume manifolds, became an active research problem in the following decades. In this setting, the proportionality principle was originally stated in \cite[Section 6.1]{thu} for finite-volume hyperbolic manifolds, where it can be used to give an elegant ``topological" proof of the celebrated Mostow rigidity theorem. Thurston's proof in this case used the isometric isomorphy of measure homology to singular homology, while Gromov gave a proof which avoided using measure homology in this case. It turns out however, that for generalizing the proportionality principle beyond hyperbolic manifolds, one could no longer avoid working with measure homology and the isometric injectivity between singular homology and measure homology became a central ingredient in proving such generalizations, where for the general (non-hyperbolic) case one has to replace the simplicial volume by the Lipschitz simplicial volume. After proofs of the proportionality for rank one locally symmetric spaces in \cite{ls}, pinched negatively curved manifolds in \cite{kk}, and manifolds with an upper curvature bound in \cite{str}, the problem was finally settled by Franceschini, who proved in \cite{fra} that the proportionality principle holds for finite-volume manifolds without any curvature restrictions.}  

\hl{One may then wonder whether the proportionality principle still holds true for spaces which have a meaningful definition of fundamental class (and hence of simplicial volume as Gromov norm as the fundamental class). The typical example here are manifolds glued together along lower-dimensional pieces, where the fundamental class is the formal sum of fundamental classes of the glued manifolds. In the Riemannian setting, the glueing should of course be done isometrically. To generalize the proportionality principle to such spaces one would again have to use measure homology and in particular need isometric injectivity of $\iota_*$.}

An example from \cite{pz} shows that $\iota_*$ is not always injective. However, the example constructed there is in some sense an artificial one: 
it relies on the existence of non-measurable sets and ultimately on the axiom of choice. So one may ask whether for more natural spaces one can still prove injectivity of $\iota_*$ as it holds for CW-complexes.

The picture below shows the convergent arcs space $CA$. It is formed by one arc $l_\infty$ and a sequence of arcs $(l_n)_{n\in \N}$ with the same endpoints as $l_\infty$ and pointwise converging to $l_\infty$. Although the arcs provide a natural cell decomposition, $CA$ is not a CW-complex because its topology is 
not the weak topology from the cell decomposition: the union $\bigcup_{n\in\N}l_n$ is not a closed subset.

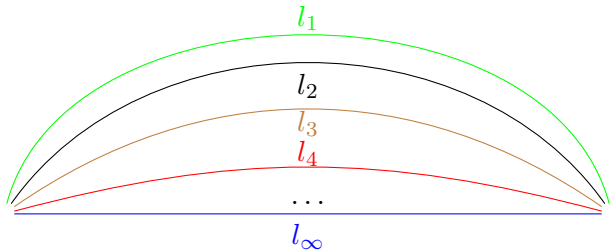
\begin{tikzpicture}
\node (A) at (0,0) {};
\node (B) at (8,0) {};
\draw[-, green] (A) to[out=75, in=105] (B);
\draw[-, ] (A) to[out=55, in=125] (B);
\draw[-, brown] (A) to[out=35, in=145] (B);
\draw[-, red] (A) to[out=15, in=165] (B);
\draw[-, blue] (A) to[out=0, in=180] (B);
\node[below, blue] at (4,0) {$l_\infty$};
\node[green] at (4,2.6) {$l_1$};
\node at (4,1.7) {$l_2$};
\node[brown] at (4,1.2) {$l_3$};
\node[red] at (4,0.8) {$l_4$};
\node at (4,0.15) {$\ldots$};
\end{tikzpicture}

It was shown in \cite[Section 6]{zas} that $H_1(CA;\R)\to {\mathcal H}_1(CA)$ is not surjective. On the other hand, \cite[Theorem 2.8]{pr3} computes the measure homology of $CA$ and the proof implies in particular that $H_*(CA;\R)\to{\mathcal H}_*(CA)$ is injective.

The convergent arcs space is a ``mildly wild" space in the sense that it is semi-locally simply connected and that it has countable fundamental group. When we collapse the arc $l_\infty$ to a point, then we obtain the 
``earring space" (pictured below to the left, until 2021 usually called
``the Hawaiian earring") which is not semi-locally simply connected connected and which has uncountable (and very complicated) 
fundamental group. This is an example of a ``really wild" space. More generally, one could consider shrinking wedges of manifolds: the earring
space is the shrinking wedge of circles, and the Barratt-Milnor sphere (pictured below to the right) is the shrinking wedge of spheres. The Barratt-Milnor sphere is semi-locally simply connected and has countable (actually trivial) fundamental group, however its higher homotopy groups are not countable.

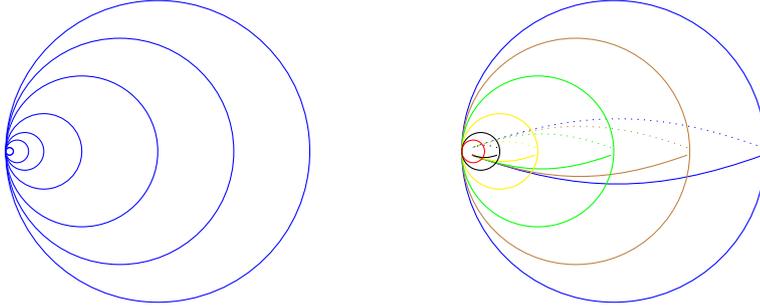
\begin{tikzpicture}

\draw[-,blue](2,0)circle(2);
\draw[-,blue](1.5,0)circle(1.5);
\draw[-,blue](1,0)circle(1);
\draw[-,blue](0.5,0)circle(0.5);
\draw[-,blue](0.25,0)circle(0.25);
\draw[-,blue](0.15,0)circle(0.15);
\draw[-,blue](0.05,0)circle(0.05);

\draw[-,blue](8,0)circle(2);
\node[] (a) at (6,0) {};
\node[] (b) at (10.1,0) {};
\draw[-,blue] (a) to[out=340, in=200] (b);
\draw[-,blue,dotted] (b) to[out=160,in=20] (a);

\draw[-,brown](7.5,0)circle(1.5);
\node[] (c) at (9.1,0) {};
\draw[-,brown](a) to[out=340, in=200] (c);
\draw[-,brown,dotted](c) to[out=160,in=20] (a);
\draw[-,green](7,0)circle(1);
\node[] (d) at (8.1,0) {};
\draw[-,green](a) to[out=340, in=200] (d);
\draw[-,green,dotted](d) to[out=160,in=20] (a);
\draw[-,yellow](6.5,0)circle(0.5);
\node[] (e) at (7.1,0) {};
\draw[-,yellow](a) to[out=340, in=200] (e);
\draw[-,yellow,dotted](e) to[out=160,in=20] (a);
\draw[-,black](6.25,0)circle(0.25);
\node[] (f) at (6.6,0) {};
\draw[-,black](a) to[out=340, in=200] (f);
\draw[-,black,dotted](f) to[out=160,in=20] (a);
\draw[-,red](6.15,0)circle(0.15);
\node[] (g) at (6.3,0) {};
\draw[-,red](a) to[out=340, in=200] (g);
\draw[-,blue](0.05,0)circle(0.05);
\end{tikzpicture}

Although ultimately we would like to say something about injectivity of $\iota_*$ for ``really wild" spaces of uncountable fundamental group like the earring space, in this paper we will pursue
a more modest goal: we
will prove injectivity of the canonical homomorphism for two classes of
``mildly wild" spaces, i.e., spaces  
which have countable fundamental group \hl{and will be semilocally
    simply connected. Observe that these two conditions are not
    completely independent: \hl{If a first countable non-semilocally simply 
    connected space is at least homotopically Hausdorff (i.e.\ there
    does not exist a non-trivial homotopy class of closed paths,
    that could be represented in an arbitrary small neighbourhood), it must by 
    \cite[Prop.4.6]{fz} have an uncountable fundamental group.
    Therefore only spaces that the not homotopically Hausdorff 
    are amongst non-semilocally simply connected spaces candidates for 
    potentially having a countable fundamental group. Such spaces
    exist, but apparently one needs to involve the axiom
    of choice in order to construct corresponding examples} (cf.\
    \hyperref[afortRem]{Remark \ref*{afortRem}},
    \hyperref[afortExpl]{Example \ref*{afortExpl}}). The class of 
        spaces to which our results
    apply do not include the earring space}, 
but they apply to  various generalizations of the convergent arc space.

The proofs of our two cases are independent and will use different methods.

The first result is the following:
\begin{thm}\label{thm1}Let $X$ be a topological space, which is $T_1$, second countable, has countable fundamental group and admits a contractible generalized universal covering space $\widetilde{X}$ in the sense of \cite{fz}. 

Then the kernel of $\iota_*\colon H_*(X;\R)\to {\mathcal{H}}_*(X)$ is contained in the zero-norm subspace with respect to the Gromov norm on $H_*(X;\R)$.

In particular, if for some $k$ the Gromov norm on $k$-th homology is an actual norm, that is $\Vert x\Vert\not=0$ for all $x\in H_k(X;\R)\setminus\left\{0\right\}$, then 
$\iota_k\colon H_k(X;\R)\to {\mathcal{H}}_k(X)$ is injective.\end{thm}

We will recall the definition of the Gromov norm in \hyperref[gromovnorm]{Section \ref*{gromovnorm}}.
The assumption on non-vanishing of the Gromov norm seems to be a more severe restriction than the others.
Therefore our second theorem may be more useful. 
\hl{In this theorem we are considering a space $X$
which has a decomposition $X=\bigcup_{n\in \N}Y_n$,
where each of the $Y_n$ has a CW-topology, but $X$, due to the way how the
$Y_i$ accumulate, has not. However, we will requite that all
intersections of the $Y_i$ are compatible with the CW-structure,
i.e., that all intersections 
\begin{equation}
Y_I:=\bigcap_{i\in I}Y_i, ~~I\subset\N
\label{eq:thm2-1}
\end{equation}
are sub-CW-complexes of each $Y_i$ with $i\in I$. 
\hl{We especially consider the intersection set 
\begin{equation}
Y^\cap = 
\{y \in Y \mid \exists {i\neq j\in \N}\hbox{ with } y \in Y_i \cap Y_j\}
\label{eq:thm2-2}
\end{equation}
and in particular its system of connected components that we denote by 
\begin{equation}
\{Y^\cap_\nu \mid \nu \in \Ups\}. 
\label{eq:thm2-3}
\end{equation}
By our assumption on the 
intersection complexes $Y_I$, each {\it``intersection component"}\/ $Y^\cap_\nu$  
is a union of components (i.e., in particular of subcomplexes) of some
of the $Y_I$ and has therefore a unique cell-subdivision, induced by the
$Y_I$ and $Y_n$. For each $\nu$ the index set 
\begin{equation}
I_\nu :=\{
i \in \N \mid Y_i \cap Y_\nu^\cap\neq\emptyset\} 
\label{eq:thm2-4}
\end{equation}
is defined}. With the help of this terminology we can now state: 
\begin{thm}\label{thm2}Let $X$ be a first-countable $T_4$-space which is a countable union 
$$X=\bigcup_{n\in \N}Y_n$$ of closed subsets $Y_n, n\in\N$ which have the structure of finite CW-complexes and such that 
all intersections $Y_I
$ are sub-CW-complexes of each $Y_i$ with $i\in I$. 
With the above in  \hyperref[eq:thm2-1]{(\ref*{eq:thm2-1})}--\hyperref[eq:thm2-4]{(\ref*{eq:thm2-4})}
defined notation, we assume that
\begin{itemize}
\item[(i)] 
\hl{Each $Y^\cap_\nu$ has a finite cell-subdivision}.
\item[(ii)] 
\hl{Each $Y^\cap_\nu$ is contractible.}
\item[(iii)] \hl{%
every intersection component $Y_\nu^\cap$ has a neighbourhood $U_\nu^\cap\subset 
X$, that deformation retracts onto $Y_\nu^\cap$, and that is relatively open in $\bigcup_{i\in I_\nu}Y_i$},
\item[(iv)] \hl{%
for different intersection components $Y_\nu^\cap, Y_\mu^\cap$,  
we have $U_\nu^\cap\cap U_\mu^\cap=\emptyset$,}%   
\item[(v)]  \hl{%
no $x\in X$ is a limit of a sequence $x_\nu\in Y_{\nu}^\cap$ for pairwise distinct \hl{intersection components} $Y_{\nu}^\cap$.}
\end{itemize}
Then $\iota_*\colon H_k(X;\R)\to {\mathcal{H}}_k(X)$ is injective for $k\ge 2$.
\end{thm}}

The result is also true for $k=0$, where it follows from Theorem 4.1 in \cite{pz}. It should also be true for $k=1$ (work in progress).

\hl{Roughly speaking, the above assumptions of 
\hyperref[thm2]{Theorem \ref*{thm2}} demand that the
intersection sets of the subspaces $Y_i$, apart from having to be consistent
with their CW-structure, subdivide into several
(maybe infinitely many) well-isolated contractible domains, each of
which is the deformation retract from some kind of neighbourhood. 
\hyperref[CAstrich]{Examples \ref*{CAstrich}}--%
\hyperref[Hyperbolas]{ \ref*{Hyperbolas}}
contain additional information on the
origins and finer technical points of 
these assumptions. In this context observe that the
converging arcs space $CA$ is also in this case just a very special
example, because each point $\in CA$ belongs either to all $l_i$,
or just to one of them. The above stated theorem allows also for more
 general scenarios, in particular that different subsets 
of the $Y_i$ intersect in different regions, or that there exists
one point in $\bigcap_{i=1}^m Y_i$, but the intersection
complexes for each of the 
\begin{equation}
2^m - m - 1 \label{eq:AnzTeilm}
\end{equation}
different subsets of $\{1,2,3,\ldots,m\}$ containing at least
two indices gives each time a different different complex.}   

\hl{Our Condition (iii) of \hyperref[thm2]{Theorem \ref*{thm2}}}
rules out examples like the earring space. 

An example of spaces to which we apply these results
are the convergent $Y$-spaces defined in \hyperref[convy]{Definition \ref*{convy}},
which are not CW-complexes and which are sort of a generalization of the convergent arcs space $CA$. They are constructed by gluing countably many copies of a CW-complex along finitely many points such that all the copies $Y_n$ of $Y$ accumulate at one copy $Y_\infty$.

This kind of example would satisfy the assumptions of \hyperref[thm1]{Theorem \ref*{thm1}} \hl{(in particular the one about the non-vanishing of the Gromov norm)} only under additional assumptions (e.g., when $Y$ is a negatively curved manifold, see \hyperref[negcurv]{Section \ref*{negcurv}}).
However we will see in \hyperref[hypex]{Section \ref*{hypex}} that it satisfies all assumptions of \hyperref[thm2]{Theorem \ref*{thm2}} whenever $Y$ is 
a compact, smooth manifold, and we thus get injectivity of $\iota_*$ for any convergent $Y$-space $X$.

In \cite{zas1}, the third author will exhibit an example of a space $X$, which is a countable union of CW-complexes, but where the intersection complexes
are neither contractible, nor contained 
in bigger contractible 
\hl{intersection components},
and where $\iota_*$ is not injective. (This is also the first such example which does not rely on the existence of non-measurable sets.)

In a final section we develop a technical device which is not used in the proof of \hyperref[thm2]{Theorem \ref*{thm2}}, but should be of independent interest for approaches to measure homology. We show that one can reduce the problem of computing measure homology (and the dual notion of measurable cohomology) to computing it for the 
subcomplex of simplices with all vertices in 
a given basepoint (\hyperref[lem2]{Lemma \ref*{lem2}} and \hyperref[boundedx0]{Corollary \ref*{boundedx0}}). To prove this result we are imposing a condition that $X$ can be covered by finitely many Borel sets of compact closure contractible in $X$. This may look like a technical condition, but it may not be avoidable as \hyperref[noninj]{Example \ref*{noninj}}
shows. 

When applied to CW-complexes our argument is similar but simpler than the one in \cite{loe} which did not restrict to simplices with vertices in a
basepoint and therefore needed a larger effort to prove the technical \cite[Lemma A.1]{loe} on existence of a measurable section. Our argument, together with countability of the fundamental group, actually 
also
provides such a measurable section from pointed simplices in $X$ to pointed simplices in its (generalized) universal covering.
It also contains an argument, why under more general conditions
(e.g., for the earring space) such a measurable section cannot exist
(\hyperref[uncount]{Theorem \ref*{uncount}}).

For the proof of \hyperref[thm1]{Theorem \ref*{thm1}} we show that (under the made assumptions) the action of the deck transformation 
group on a generalized universal covering space has a Borel-measurable fundamental domain. This might be of independent interest, here we use it
to show in \hyperref[iotainjective]{Section \ref*{iotainjective}}
that (in the case of countable 
fundamental groups) the homomorphism from measurable bounded cohomology to bounded cohomology is an isometric isomorphism.

We remark that the reader interested in \hyperref[thm2]{Theorem \ref*{thm2}} 
should only need \hyperref[mhdefs]{Sections \ref*{mhdefs}}, 
\hyperref[coning]{ \ref*{coning}} and \hyperref[negcurv]{ \ref*{negcurv}}
to be able to follow the proof of 
\hyperref[thm2]{Theorem \ref*{thm2}} in
\hyperref[proof2]{Section \ref*{proof2}}.

{\sl Conventions}: spaces of simplices will be equipped with the compact-open topology and ``measurable" will always mean Borel-measurable with respect to that topology. ``Measures" will always mean signed measures, i.e., differences of two non-negative measures. A ``$G$-module" will always mean a Banach space $V$ which is a module over the group ring ${\mathbb Z} G$ and such that $\parallel gv\parallel\le \parallel v\parallel$ for all $g\in G,v\in V$.

\section{Preliminaries}
\subsection{Measure homology}\label{mhdefs}
Let us start with recalling the definition of measure homology (or Milnor-Thurston-homology) from \cite[Definition 1.8]{zas}.
\begin{defi}For a topological space $X$ and $k\in \N$ we denote its set of singular $k$-simplices, i.e., 
of continuous maps from the standard simplex $\Delta^k$ to $X$, by 
$map(\Delta^k,X)$. We equip $map(\Delta^k,X)$ with the compact-open-topology and the corresponding $\sigma$-algebra of Borel sets. \end{defi}
\begin{defi}\label{mcyc}For a topological space $X$ and $k\in\N$ let 
$${\mathcal C}_k(X)=\left\{\mu\mid \mu\mbox{ is a compactly determined measure on }map(\Delta^k,X), \parallel\mu\parallel<\infty\right\}.$$
\end{defi}

Here, a compactly determined measure is one that vanishes on any measurable subset of the complement of some (not necessarily measurable) compact set. (We follow the convention that a
compact set need not be Hausdorff but satisfies the Heine-Borel covering property. Such sets are sometimes called quasicompact, therefore the definition in \cite{zas} speaks of quasicompactly determined measures.) The 
variation of a signed measure is
$\Vert \mu\Vert:=\max_A\mu(A)-\min_B\mu(B)$, where the maximum resp.\ minimum are taken over all measurable sets.

It is proved in \cite[Corollary 2.9]{zas} that the canonical boundary 
operator $\partial_k$ for singular simplices, that 
is based on a linear combination 
of face maps $\partial_k=\sum_{i=0}^k(-1)^i\partial^i_k$ extends
to an operator on measures on the set $map(\Delta^k,X)$ and in 
particular to an operator 
$\partial_k\colon {\mathcal C}_k(X)\to {\mathcal C}_{k-1}(X)$. Then one defines measure homology as 
$${\mathcal H}_k(X)=ker(\partial_k)/im(\partial_{k-1}).$$

\subsection{Generalized universal covering spaces}\label{guc}
\begin{defi}\label{Def1.3}(\cite[Section 1.1]{fz}) A generalized universal covering space of a path-connected topological space $X$ is a topological space $\widetilde{X}$ with a continuous surjection $p\colon\widetilde{X} \to X$ such that 
\begin{itemize}
\item[(i)] 
$\widetilde{X}$ is locally path-connected and simply-connected,
\item[(ii)] 
if $Y$ is path-connected and locally path-connected, then every pointed continuous map $f\colon (Y,y)\to(X,x)$ with $f_*(\pi_1(Y,y))=1$ admits unique pointed liftings, that is, for each $\tilde{x}\in p^{-1}\{x\}$ there is a unique pointed continuous map $g\colon(Y,y)\to(\widetilde{X},\tilde{x})$ with $p\circ g=f$.\end{itemize}
\end{defi}
A generalized universal covering space, if it exists, is in one-to-one correspondence with the 
homotopy classes of paths in $X$ which emanate from a fixed $x_0\in X$. (For more details 
see \cite[Section 2]{fz}.) 
 
A generalized universal covering is a Serre fibration, thus one has $\pi_k\widetilde{X}\cong \pi_kX$ for $k\ge2$, see \cite[Section 1.2]{fz}. Moreover the deck transformation group of $p\colon \widetilde{X}\to X$ is isomorphic to the fundamental group $\pi_1X$, and it acts freely and transitively on each fibre, see \cite[Proposition 2.14]{fz}.

For our arguments, the most important property of the generalized universal 
covering space will be that the lifts of a singular simplex $\sigma\colon \Delta^k\to X$ form exactly a $G$-orbit 
of singular simplices in $\widetilde{X}$, where $G\cong\pi_1(X,x_0)$ is the deck 
transformation group. Moreover the lifts of the simplices with all vertices in $x_0\in X$ 
are exactly the simplices with vertices in $G\tilde{x}_0$, for a preimage $\tilde{x}_0\in\widetilde{X}$ of $x_0$.

\subsection{Relatively injective modules and bounded cohomology}

\begin{defi}\label{relinjbc}For a topological space $X$ we let 
$$C_b^k(X):=B(map(\Delta^k,X),\R)=\left\{f\colon map(\Delta^k,X)\to \R\mid f\mbox{\ is\ bounded}\right\}$$ be the vector space of 
bounded cochains. It is a Banach space with the norm $\Vert f\Vert=sup\left\{\vert f(\sigma)\vert\colon \sigma\in map(\Delta^k, X)\right\}$.
The usual coboundary operator $$\delta_kf(\sigma)=\sum_{i=0}^k(-1)^if(\partial_i\sigma)$$
makes $C_b^*(X)$ a cochain complex and its cohomology is denoted by $H_b^*(X)$ and called the bounded 
cohomology of $X$.\end{defi}

If $X$ comes with an action of a group $G$, then $C_b^k(X)$ becomes a $G$-module via the induced action. In particular,
if $\widetilde{X}\to X$ is a generalized universal covering space and $G\cong\pi_1(X,x_0)$ its group of 
deck transformations, then $C_b^k(\widetilde{X})$ is naturally understood as a $G$-module and this will always be meant 
when we refer to $C_b^k(\widetilde{X})$ as a $G$-module. For readers familiar with \cite{mon} we want to mention that,
although $\pi_1(X,x_0)$ can be topologized as a non-discrete topological group acting continuously on $\widetilde{X}$, this is not what we are going
to do and we rather consider $G$ as a discrete group. In particular, for the proof of \hyperref[reso]{Lemma \ref*{reso}} it will 
be sufficient to consider the module $B(G,V)$ of bounded functions rather than the module of continuous, bounded functions and
so we will not need the general results on continuous bounded cohomology from \cite{mon} but only the results on bounded cohomology from \cite{iva}.

It is often useful to compute bounded cohomology via other resolutions. The general setting for this to work are strong resolutions by relatively injective modules.

\begin{defi}\label{relinje}Let $G$ be a topological group. A $G$-module $U$ is called relatively injective if any diagram of the form 

\begin{tikzpicture}
\node at (0,2) {$V_1$};
\draw[->] (0.5,2)--(3.5,2);
\node[above] at (2.5,2) {$i$};
\draw[->] (3.5,1.9)--(0.5,1.9);
\node[below] at (2.5,1.9) {$\sigma$};
\draw[->] (0.5,1.8)--(3.5,0.2);
\node[left, below] at (2,1) {$\alpha$};
\node at (4,2) {$V_2$};
\draw[->,dotted] (4,1.7)--(4,0.3);
\node[right] at (4,1) {$\beta$};
\node at (4,0) {$U$};
\end{tikzpicture} 

\noindent can be completed.
Here $i\colon V_1\to V_2$ is an injective morphism of $G$-modules, $\sigma\colon V_2\to V_1$ is a bounded (not necessarily 
$G$-equivariant) linear operator with $\sigma\circ i=id$ and $\Vert\sigma\Vert\le1$, $\alpha$ is 
a $G$-morphism, and we want $\beta$ to be a $G$-morphism with $\beta\circ i=\alpha$ and $\Vert\beta\Vert\le\Vert\sigma\Vert$.\end{defi}

\begin{defi}\label{strong}A strong resolution of a $G$-module $U$ is an exact sequence of $G$-modules and $G$-morphisms 

\begin{tikzpicture}
\node at (0,2) {$0$};
\draw[->](0.5,2)--(1.5,2);
\node at (2,2) {$U$};
\draw[->](2.5,2)--(3.5,2);
\node[above] at (3,2) {$\delta_{-1}$};
\node at (4,2) {$U_0$};
\draw[->](4.5,2)--(5.5,2);
\node[above] at (5,2) {$\delta_{0}$};
\node at (6,2) {$U_1$};
\draw[->](6.5,2)--(7.5,2);
\node[above] at (7,2) {$\delta_{1}$};
\node at (8,2) {$U_2$};
\draw[->](8.5,2)--(9.5,2);
\node[above] at (9,2) {$\delta_{2}$};
\node at (10,2) {$\ldots$};
\end{tikzpicture}

\noindent for which there
exists a sequence of linear (not necessarily $G$-equivariant) operators $\kappa_n\colon U_n\to U_{n-1}$ such that $\delta_{n-1}\kappa_n+\kappa_{n+1}\delta_n=id$ and $\Vert \kappa_n\Vert\le 1$ for all $n\ge 0$ and $\kappa_0\delta_{-1}=id$.\end{defi}

According to \cite[Lemma 7.2.6]{mon} the trivial $G$-module $\R$ has a strong resolution by 
relatively injective $G$-modules, and any two such resolutions are chain homotopy equivalent. 
In particular the cohomology of the $G$-invariants of the resolution does not depend on the chosen resolution. 
This cohomology is, by definition, the continuous bounded cohomology of $G$, denoted by $H^*_{cb}(G)$. As said, we only
consider the bounded cohomology $H^*_b(G)$ defined by equipping $G$ with the discrete topology.
We will 
need the following two facts, which can be found for example in \cite{iva} or in the more general setting of continuous bounded cohomology 
in \cite{mon}.
\begin{lemma}\label{relin}

i) (\cite[Lemma 3.2.2]{iva}) For any Banach space $V$, the $G$-module $B(G,V)$ of bounded functions with values in $V$ is relatively injective.

ii) (\cite[Lemma 3.3.2]{iva}) Let $$0\to U\to U_1\to U_2\to\ldots$$ be a strong resolution of the $G$-module $U$ and $$0\to V\to V_1\to V_2\to\ldots$$ be a complex of relatively injective $G$-modules, then any $G$-morphism $U\to V$ can be extended to a $G$-morphism of complexes and any two such extensions are $G$-chain homotopic.\end{lemma}

The following lemma is well-known for CW-complexes and more generally for semi-locally simply connected spaces, and 
we are going to show that the same proof also works for spaces that admit a generalized 
universal covering space in the sense of \hyperref[guc]{Section \ref*{guc}}.

\begin{lemma}\label{reso}Let $\widetilde{X}\to X$ be a generalized universal covering space and $G$ its group of deck transformations. Then
$$0\to\R\to C_b^0(\widetilde{X})\to C_b^1(\widetilde{X})\to C_b^2(\widetilde{X})\to\ldots$$
is a strong resolution by relatively injective $G$-modules. In particular one has an isometric isomorphism $H_b^*(X)=H_{cb}^*(G)$.\end{lemma}
\begin{pf} We will prove this by copying the argument in the proof of \cite[Theorem 4.1]{iva}.

By \hyperref[relin]{Lemma \ref*{relin}i)}, $B(G,V)$ is relatively injective for each 
Banach space $V$.

By the axiom of choice there exists a set $F\subset \widetilde{X}$ meeting each $G$-orbit exactly once. Let $map((\Delta^k,v_0),(\widetilde{X},F))$ be the set of those singular simplices which send the first vertex of the standard simplex to $F$. We make $B^k(\widetilde{X},F):=B(map((\Delta^k,v_0),(\widetilde{X},F)),\R)$ a Banach space by equipping it with the sup-norm. Then there is an obvious isomorphism $$C_b^k(\widetilde{X})=B(G,B^k(\widetilde{X},F))$$
and thus $C_b^k(\widetilde{X})$ is a relatively injective $G$-module.

By simple connectivity of $\widetilde{X}$ and \cite[Theorem 2.4]{iva} there is a contracting algebraic homotopy for $C_b^*(\widetilde{X})$. Hence we have a strong resolution.

\end{pf}

\section{Measurable bounded cohomology - proof of Theorem 1}\label{proof1}

\subsection{Definitions}

In the previous section we defined bounded cohomology, now we are going to define measurable bounded cohomology.

\begin{defi}\label{measbc}Let $X$ be a topological space and again $map(\Delta^k,X)$ equipped with the compact-open-topology and the corresponding $\sigma$-algebra of Borel-measurable sets. We let 
$${\mathcal{C}}_b^k(X)=\left\{f\colon map(\Delta^k,X)\to \R\mid f\mbox{\ is\ Borel\ measurable\ and\ bounded}\right\}$$ be the measurable bounded cochains.\end{defi}

The usual coboundary operator makes ${\mathcal{C}}_b^k(X)$ into a cochain complex and its cohomology is denoted by ${\mathcal{H}}_b^*(X)$, see \cite[Section 3.4]{loe}. The inclusion $\iota$ induces a homomorphism
$$\iota^*\colon {\mathcal{H}}_b^*(X)\to H^*_b(X;\R)$$
from the measurable bounded cohomology to the bounded cohomology.

\subsection{Connecting the Gromov norm to measurable bounded cohomology}\label{gromovnorm}
The following arguments are well-known, cf.\ \cite[Section 3]{loe}. We will need them for the proof of \hyperref[thm1]{Theorem \ref*{thm1}}.

For a topological space $X$ there is an $l^1$-norm on its singular chain complex $C_*(X;\R)$ defined by 
$\Vert \sum_{i=1}^r a_i\sigma_i\Vert_1=\sum_{i=1}^r\vert a_i\vert$.
The Gromov norm on homology $H_*(X;\R)$ is defined as 
$\Vert\alpha\Vert=\inf\left\{\Vert z\Vert_1\colon \left[z\right]=\alpha\right\}$, i.e., one takes the infimum of the $l^1$-norm over all cycles $z$ representing the homology class $\alpha$.
We denote $NH_k(X)=\left\{\alpha\in H_k(X;\R)\colon \Vert\alpha\Vert=0\right\}$.
\begin{lemma}\label{gr0}Let $X$ be a topological space and $k\in\N$. If 
$$\iota^*\colon {\mathcal{H}}_b^k(X)\to H_b^k(X;\R)$$
is an epimorphism, then 
$$ker(\iota_*\colon H_k(X;\R)\to {\mathcal{H}}_k(X))\subset NH_k(X).$$
\end{lemma}
\begin{pf}Assume there is some $\alpha\in H_k(X;\R)$ with $\Vert \alpha\Vert \not =0$ and $\iota_*(\alpha)=0$.

By \cite[Section 1.1]{gro} the $l^1$-norm on $H_k(X;\R)$ is dual to the norm on $H_b^k(X)$, which for $\phi\in H_b^k(X)$ is defined as infimum of $\Vert f\Vert$ over all bounded cocycles $f$ representing $\phi$. 
In particular, there is some $\phi\in H^k_b(X)$ with $\langle \phi,\alpha\rangle=1$. By assumption there is some $\psi\in {\mathcal{H}}^k_b(X)$ with $\iota^*\psi=\phi$. Then 
$$1=\langle \phi,\alpha\rangle=\langle\iota^*\psi,\alpha\rangle=\langle \psi,\iota_*\alpha\rangle=0,$$
yielding a contradiction.
\end{pf}

\subsection{Construction of a measurable fundamental domain}\label{funddom}
The following \hyperref[fun]{Lemma \ref*{fun}} will be used in this paper for \hyperref[isoiso]{Proposition \ref*{isoiso}} in \hyperref[iotainjective]{Section \ref*{iotainjective}}, though we think that it might be of independent interest.

Properly discontinuous group actions have a measurable fundamental domain, see \cite[Chapter 7, Exercises for \S2, Ex.\ 12]{bou}. However, the action of the group of deck transformations on a generalized universal covering space is in general not properly discontinuous. We are going to show that (under weak assumptions) one can nevertheless adapt the argument and obtain a measurable fundamental domain.
\begin{lemma}\label{fun}Let $X$ be a second-countable $T_1$-space and \hl{assume that 
the set of points at which $X$ is not semi-locally simply connected is at most countable}. If there exists a generalized universal covering space 
$p\colon\widetilde{X}\to X$, then the action of the deck transformation group
$\Gamma\cong \pi_1(X,x_0)$ on $\widetilde{X}$ has a Borel-measurable fundamental domain.\end{lemma}
\begin{pf}

Let $N$ be the countably many points where $X$ is not semi-locally simply connected. \hl{At} any
$x\in X\setminus N$ and each $\widetilde{x}\in p^{-1}\{x\}$ 
there is an open neighbourhood $\widetilde{U_{\tilde{x}}}\subset \widetilde{X}$ 
such that the restriction of $p$ to that neighbourhood is injective. (Namely one can take a neighbourhood $V_x\subset X$ satisfying $im(\pi_1(V_x,x)\to\pi_1(X,x))=0$ and a connected component $\widetilde{U}_{\tilde{x}}$ 
of its preimage $p^{-1}(V_x)$. Note that this does not necessarily surject onto $V_x$. The intersection of 
$\widetilde{U_{\tilde{x}}}$ with any $\Gamma$-orbit has at most one element.)

For $x\in X$ choose some $\tilde{x}\in p^{-1}\{x\}$ and 
let $U_x=p(\widetilde{U_{\tilde{x}}})\subset X$ be the image 
of $\widetilde{U_{\tilde{x}}}$. Second-countable spaces have 
the Lindel\"of property and hence there is a countable family of $U_x$ that covers $X$.

With these preparations we define a measurable fundamental domain as follows. Let $\left\{U_1,U_2,U_3,\ldots\right\}$ be 
an enumeration of the countable family of $U_x$'s and $\left\{\widetilde{U_1},\widetilde{U_2},\widetilde{U_3},\ldots\right\}$ the corresponding subsets of $\widetilde{X}$. Then $$W_1:=\widetilde{U_1}$$
$$W_2:=\widetilde{U_2}\cap (\widetilde{X}\setminus \Gamma \widetilde{U_1}) $$
$$W_3:=\widetilde{U_3}\cap (\widetilde{X}\setminus(\Gamma \widetilde{U_1}\cup\Gamma \widetilde{U_2}))$$
$$\ldots$$
are all Borel-measurable. (One should pay attention that we are using the possibly uncountable unions 
$\Gamma \widetilde{U_i}$, however these are unions of open sets and so no problem arises.) So
$$\bigcup_{n\in\N}W_n$$
is a measurable set, and one easily checks that it contains exactly one point from each $\Gamma$-orbit 
not meeting $p^{-1}(N)$. Adding one point of each of the countably many $\Gamma$-orbits in $p^{-1}(N)$ we obtain 
a fundamental domain. 

We claim that the $T_1$-property for $X$ implies the $T_1$-property for $\widetilde{X}$. Namely, as pointed out in \cite[Lemma 2.10, Lemma 2.11]{fz} 
a space possessing a generalized universal covering space must be homotopically Hausdorff and then two points lying on the same
fibre of $p$ can be even separated in the $T_2$-sense. For two points not lying on the
same fibre of $p$, the analogue of the arguments contained in these lemmas for the
$T_2$-case (taking complete preimages of neighbourhoods with corresponding
separation properties), gives for our assumption that the points can be at
least separated in the $T_1$-sense, giving the claim.
Finally, the $T_1$-property for $\widetilde{X}$ implies that points are closed and their countable 
union is a Borel set, so that the constructed fundamental domain is Borel-measurable.  \end{pf}

\begin{remark} 
A more explicit construction of the fundamental domain exists for
spaces that satisfy a non-positive curvature condition, 
that for the generalized universal covering
space amounts to a global CAT(0)-condition. 
\hl{However, for the just considered CAT(0)-base spaces} 
we can connect each point via {\hl a  
shortest} geodesic to
a base point \hl{making a choice if there should be different
geodesics of the same length}. Then the domain covered by the lifts of the shortest geodesics starting at one lift of the base point
will form a fundamental domain.
In a purely topological context, path systems that satisfy similar
properties
as CAT(0)-geodesics and could be used for analogous constructions,
have been axiomatically described and introduced in \cite{asm1},\cite{asm2}
under the name ``arc-smooth systems". Actually, in our context,
when adapting \hl{the conditions of arc-smooth systems} 
(that can only be satisfied
for a kind of covering space) to the base space,
we would be happy with a bit less. Instead of having one
uniquely defined path between any two points, it would suffice
to have for each point one uniquely defined path connecting to some
base point, usually continuously depending on the other endpoint,
but for a non-contractible base space there must be border-zones
where this continuity-condition cannot be satisfied;
such path-systems are sometimes called a ``combing".
In our case we would need a combing that is prefix-closed, i.e.,
each path starting on the trace of another combing path $c$
or crossing the trace of another combing path $c$ would have to follow the
same trace as the path $c$ to the base point. With one combing path
starting at each point of the space,
then the  set covered by the lift of at least  one
of the combing paths, starting at one lift of the base-point,
will form a fundamental domain, and for a sensible choice of the
border-zones
there is a chance that the result will be a measurable set.
\end{remark}

\subsection{Measurable coning construction}\label{coning}
The following construction will later be applied to the (generalized) universal covering $\widetilde{X}$ of a topological space $X$.
\begin{defi}\label{cone}Let $(\widetilde{X},x_0)$ be a pointed topological space. It is said to have a measurable (resp.\ continuous) coning construction if there is a sequence of Borel-measurable (resp.\ continuous) maps
$$L_i:map(\Delta^i,\widetilde{X})\to map(\Delta^{i+1},\widetilde{X})$$ such that for each $\sigma\in map(\Delta^i,\widetilde{X})$ the $0$-th vertex of $L_i(\sigma)$ is $x_0$ and 
$$\partial_0L_i(\sigma)=\sigma, $$
\begin{equation}
\partial_kL_i(\sigma)=L_{i-1}(\partial_{k-1}\sigma)\mbox{ for }k=1,\ldots,i+1,
\label{eq:cone3}
\end{equation}
where by $\partial_k\colon map(\Delta^{i+1},X)\to map(\Delta^i,X)$ 
for $k=0,\ldots,i+1$ we mean the face
map omitting the $k$-th vertex.\end{defi}
\begin{lemma}\label{contr}A topological space $\widetilde{X}$ has a continuous coning construction if it is contractible.\end{lemma}
\begin{pf}\hl{A similar construction has been used in the proof of \cite[Proposition 5.3]{fribc}, we repeat it} for convenience of the reader:\neuzl   
Assume $\widetilde{X}$ is contractible. Then there is an $x_0\in\widetilde{X}$ 
and a continuous map $H\colon\widetilde{X}\times\left[0,1\right]\to\widetilde{X}$ with $H(x,0)=x\,,\;\;H(x,1)=x_0$ for all $x\in\widetilde{X}$. For a singular $i$-simplex $$\sigma\colon\Delta^i\to\widetilde{X}$$ the map $$h\colon \Delta^i\times\left[0,1\right]\to\widetilde{X}$$
$$(x,t)\mapsto H(\sigma(x),t)$$
factors over the canonical projection $$\Delta^i\times\left[0,1\right]\to\Delta^{i+1},$$
which collapses $\Delta^i\times\left\{1\right\}$ to the $0$-th vertex of $\Delta^{i+1}$. So the map $h$ defines a singular $(i+1)$-simplex 
$$L_i(\sigma)\colon\Delta^{i+1}\to X$$ and it is easy to check that this assignment has the desired properties.\end{pf}

\begin{lemma}\label{1new}
If a coning construction is available, then for any (measurable or singular)
$k$-cycle $z$ a $(k+1)$-dimensional chain can be defined with
$\partial L_k(z) = z$, provided that $k \ge 1$. 
\end{lemma}
\begin{pf}                                   
When computing the complete boundary $\partial$ by summing over the 
face-maps $\partial_i$, \hyperref[eq:cone3]{Formula (\ref{eq:cone3})}
together with the fact that $
\partial(z)=0$ implies that only the $\partial_0$-boundary remains
and we have $\partial L_k(z)=z$.   

This conclusion requires $k\ge1$. Namely, if we consider
a point $P$ as a zero-cycle, we do not have $L_{-1}$
defined, and despite $\partial P=\emptyset$, the fact that 
$\partial L_0(P)= P - x_0$ implies in particular that 
$\partial L_0(P)= x_0 \neq 0$.
\end{pf}
\subsection{Resolution by measurable bounded cochains}

The following \hyperref[mbc]{Lemma \ref*{mbc}} will be a main ingredient in the proof of \hyperref[thm1]{Theorem \ref*{thm1}}.  Its proof is essentially copied from \cite[Theorem 2.4]{iva}, which proves the analogous result for (non-measurable) bounded cohomology.

\begin{lemma}\label{mbc}Let $\widetilde{X}\to X$ be a generalized universal
covering space and $G$ its group of deck transformations. Assume that 
$\widetilde{X}$ is contractible. Then 
$$0\to\R\to {\mathcal{C}}_b^0(\widetilde{X})\to {\mathcal{C}}_b^1(\widetilde{X})\to {\mathcal{C}}_b^2(\widetilde{X})\to\ldots$$
is a strong resolution by $G$-modules, where the maps in the resolution are $\delta_{-1}\colon\R\to{\mathcal{C}}_b^0(\widetilde{X})$ sending real numbers to constant functions, and for $i\ge 0$ the coboundary operator $\delta_i\colon {\mathcal{C}}_b^i(\widetilde{X})\to {\mathcal{C}}_b^{i+1}(\widetilde{X})$ from \hyperref[relinjbc]{Definition \ref*{relinjbc}}. \end{lemma}
\begin{pf}

By \hyperref[contr]{Lemma \ref*{contr}} we have a measurable (even continuous) coning construction for a fixed base point $x_0\in \widetilde{X}$. Dualizing \hyperref[cone]{Definition \ref*{cone}} via $$(\kappa^i(f))(\sigma):=f(L_{i-1}(\sigma))$$
 for $i\ge 1$ 
yields  
homomorphisms
$$\kappa^i\colon {\mathcal{C}}_b^{i}(\widetilde{X})\to {\mathcal{C}}_b^{i-1}(\widetilde{X})$$
for $i\ge 1$
such that
$$\delta_{i-1}\kappa^i+\kappa^{i+1}\delta_i=id$$
for all $i\ge 0$,
where for $i=0$ we define $\kappa^0\colon {\mathcal{C}}_b^0(\widetilde{X})\to\R$ by sending $f$ to $f(x_0)$ for the chosen base point $x_0$.

Because $L_i$ sends each simplex to another simplex, we clearly have $\Vert \kappa^i\Vert\le 1$ for all $i$.
\end{pf}

\subsection {$\iota^*$ is an isomorphism.}\label{iotainjective}
\begin{prop}\label{isoiso} Under the assumptions of \hyperref[mbc]{Lemma \ref*{mbc}}, if one has a measurable fundamental domain for 
the action of $G$, and if moreover $G$ is countable, then $$\iota^*\colon{\mathcal H}_b^*(X)\to H_b^*(X)$$ is an isometric isomorphism.\end{prop}
\begin{pf}We know that $G$-modules of the form $B(G,V)$ (for a Banach space $V$) are relatively injective, see \hyperref[relin]{Lemma \ref*{relin}i)}. Measurability of the fundamental domain $F$ and countability of $G$ imply that we have an isomorphism
$${\mathcal{C}}_b^k(\widetilde{X})=B(G,{\mathcal{B}}^k(\widetilde{X},F))$$
for $${\mathcal{B}}^k(\widetilde{X},F)=\left\{f\colon map((\Delta^k,v_0),(\widetilde{X},F))\to \R\mid f\mbox{\ is\ Borel\ measurable\ and\ bounded}\right\}$$
and thus relative injectivity of ${\mathcal{C}}_b^k(\widetilde{X})$. So
$$0\to\R\to {\mathcal{C}}_b^0(\widetilde{X})\to {\mathcal{C}}_b^1(\widetilde{X})\to {\mathcal{C}}_b^2(\widetilde{X})\to\ldots$$
is a strong resolution by relatively injective modules and the claim follows in view of \hyperref[relin]{Lemma \ref*{relin}}.\end{pf}

\begin{remark} If $G$ is not countable, then ${\mathcal{C}}_b^k(\widetilde{X})$ is a proper subset of $B(G,{\mathcal{B}}^k(\widetilde{X},F))$ and we do not know whether it is relatively injective.\end{remark} 
\subsection{Proof of Theorem 0.1.}
The proof of \hyperref[thm1]{Theorem \ref*{thm1}} now follows from \hyperref[gr0]{Lemma \ref*{gr0}}, \hyperref[fun]{Lemma \ref*{fun}} and
\hyperref[isoiso]{Proposition \ref*{isoiso}}.

\subsection{Example}\label{negcurv}
\begin{defi}\label{convy}Let $Y$ be a metric space. We call a metric space $X$ a convergent $Y$-space if it is a union 
$$X= Y_\infty\cup\bigcup_{n\in\N}Y_n$$ with $Y_\infty=Y$,
and if for each $n\in N$ there is a homeomorphism $$f_n\colon Y_\infty\to Y_n$$ such that there are finitely many points $y\in Y$ with $$f_n(y)=y\ \forall\ n\in \N$$ (``identification points")
and for all other points $y,y'$ 
one has $f_n(y)\not=y'$ and $f_n(y)\not= f_m(y')$ for all $n\neq m$,
but $$\lim_{n\to\infty} {\rm dist}(f_n(y),y)=0$$ 
\hl{uniformly}.
\end{defi}

\begin{lemma}\label{counta}If $Y$ has the homotopy type of a countable CW-complex, then the homotopy groups of a convergent $Y$-space are countable.\end{lemma}

\begin{pf}
We use
the well-known fact that the homotopy groups of a countable CW-complex are countable, see \cite[Theorem IV.6.1]{lw}. Although the convergent $Y$-space $X$ is
not locally path-connected, hence not a CW-complex, one can find a locally path-connected space $X^{lpc}$ with the same homotopy groups, as indicated in \cite[Section 2.1]{ks}. Let $\mathcal{O}=\left\{V\subset X\ \mbox{open}\right\}$
be the topology of $X$. For an open set $V\in\mathcal{O}$ and $x\in V$ let $U(V,x)$
be the path-component of $V$ containing $x$. The sets $U(V,x)$ for varying $x$ and $V$ form the basis of a topology
$\mathcal{O}^{lpc}$ on the set $X$. We denote the so-defined topological space by $X^{lpc}$. The identity map
$$id\colon X^{lpc}\to X$$
is continuous but in general not open. 

Since CW-complexes are locally path-connected, the above-defined
topology is induced by a metric which can just be defined
as follows:       
$$ {\rm dist}_X(x,y) :=
\left\{\vcenter{\halign{$#$&#\hfil\cr
{\rm dist}_Y(x,y),&{\vtop{\noindent\hsize=0.6\hsize if $x$ and $y$ lie in the same 
subspace $Y_i\;\;   (i \in \N \cup \{\infty\})$}\strut}\cr
{\rm min}_j({\rm dist}_Y(&
{\vtop{\noindent\hsize=0.6\hsize$x,x_j) +{\rm dist}_Y(x_j,y)),$ if $x\in Y_i\,,\;\;
                                y \in Y_{i'}\,,\;\;i\neq i'$   
and $j$ enumerates all identification points $x_j$.
}}\cr}}\right.$$             
By its construction also $X^{lpc}$ has a natural subdivision into cells.
However, since any CW-structure would have automatically to be locally
infinite, according to \cite[Theorem 1.5.17]{fp} the above defined
metric topology of $X^{lpc}$ cannot coincide with the topology of $X^{lpc}$,
that is induced as CW-topology by its natural cell-structure.
However from 
\cite[\S16, Theorem 1]{Dwkr} it follows
that the two just mentioned topologies are homotopy equivalent.
Therefore for the purpose of computing homotopy groups we can equivalently
use both of these topologies for $X^{lpc}$, and then finally 
conclude that  
according to \cite[Corollary 2.5]{ks} the map $id\colon X^{lpc}\to X$
induces isomorphisms
$$\pi_k(X^{lpc})\cong \pi_k(X)$$
for all $k$. Under the assumptions
of \hyperref[counta]{Lemma \ref*{counta}}, $X^{lpc}$ is 
homotopy equivalent to a countable
CW-complex, thus its homotopy groups are countable, and so are those of $X$.\end{pf}

\begin{lemma}\label{asphe}If $Y$ is aspherical and has the homotopy type of a countable CW-complex, then a convergent $Y$-space is aspherical.\end{lemma}
\begin{pf}
By the proof of \hyperref[counta]{Lemma \ref*{counta}} we know that $\pi_k(X^{lpc})\cong\pi_k(X)$. 
Thus it suffices to prove asphericity for CW-complexes that are obtained by identifying finite subsets of countably many aspherical CW-complexes. Since the image of a sphere can only intersect finitely many cells of $X^{lpc}$ 
(when equipped with its CW-topology) it actually suffices to prove this for 
a union of finitely many aspherical CW-complexes along finite subsets.

First consider the one-point union $Y_1\vee Y_2$ of two path-connected, aspherical CW-complexes. There is a well-known construction
(see \cite[Prop. 4.64]{hat}), which to every map $f\colon A\to B$ associates a fibration $p \colon E_f\to B$ 
and a homotopy equivalence $A\to E_f$. Namely, $$E_f:=\left\{(a,\gamma)\in A\times B^{\left[0,1\right]}\mid \gamma(0)=f(a)\right\}$$
and $p(a,\gamma)=\gamma(1)$. The fibre of $p$ is called the homotopy fibre of $f$.  
In our setting, we see that the homotopy fibre of the inclusion $$Y_1\vee Y_2\to Y_1\times Y_2$$ 
is the union of $PY_1\times \Omega Y_2$ and $\Omega Y_1\times PY_2$ along their 
intersection $\Omega Y_1\times\Omega Y_2$. (Here $PY$ means 
the path space and $\Omega Y$ the loop space.) For CW-complexes $Y_1,Y_2$ it is known that there is a weak homotopy equivalence $w$ from
the join $\Omega Y_1*\Omega Y_2$ to the homotopy fibre of the inclusion $$Y_1\vee Y_2\to Y_1\times Y_2,$$
cf.\ the final paragraph of the 
proof of \cite[Theorem 2.2]{gan}. 
If $Y_1,Y_2$ are aspherical, i.e., $\pi_k(Y_1)=\pi_k(Y_2)=0$ for $k\ge 2$, then $$\pi_k(\Omega Y_1)=\pi_{k+1}(Y_1)=0\mbox{ and }\pi_k(\Omega Y_2)=\pi_{k+1}(Y_2)=0\mbox{ for }k\ge 1,$$
i.e., $\Omega Y_1$ and $\Omega Y_2$ are weakly homotopy equivalent to discrete spaces. The loop space of a countable CW-complex has the homotopy type of a CW-complex by a theorem of Milnor, see \cite[Corollary 5.3.7]{fp}. Thus a
weak homotopy equivalence is actually
a homotopy equivalence by Whitehead's theorem.
So $\Omega Y_1$ and $\Omega Y_2$
have the homotopy type of discrete spaces, hence the join $\Omega Y_1*\Omega Y_2$ has the homotopy type of a wedge of circles. In particular, $\Omega Y_1*\Omega Y_2$ is aspherical and the weak homotopy equivalence $w$ yields that also the homotopy fibre of $$Y_1\vee Y_2\to Y_1\times Y_2$$ 
is aspherical. Moreover, asphericity of
$Y_1$ and $Y_2$ implies that $Y_1\times Y_2$ is aspherical.
This implies by the long exact sequence of homotopy groups 
$$\ldots\to \pi_k(\textit{homotopy\ fibre})\to\pi_k(Y_1\vee Y_2)\to \pi_k(Y_1\times Y_2)\to\ldots$$
that also $Y_1\vee Y_2$ is aspherical.

Next, if we identify two vertices in the same path component of a CW-complex $Y$, then the resulting CW-complex is homotopy-equivalent to the one-point union $Y\vee S^1$. Since $S^1$ is aspherical, we obtain asphericity of $Y\vee S^1$ from asphericity of $Y$.
Finally, by induction we can extend asphericity to the CW-complex obtained by identifying finite subsets.
\end{pf}

\begin{lemma}\label{ExGenCov}If $Y$ is semi-locally simply connected and first-countable, then any convergent $Y$-space has a generalized universal covering.\end{lemma}
\begin{pf}The convergent $Y$-space $X$ is semi-locally simply connected, but not locally path-connected. $X^{lpc}$ is semi-locally simply connected and locally path-connected, thus it has a (classical) universal covering $\widetilde{
X^{lpc}}$. We claim that $\widetilde{X^{lpc}}$ is a generalized universal covering of $X$.

According to \cite[Proposition 5.1]{fz} (and the characterization of generalized universal coverings from \cite[Section 1]{fz}) for a first-countable space it suffices to check the path lifting property for $\widetilde{X^{lpc}}\to X$.
But any path in $X$ lifts to a unique path in $X^{lpc}$ (see \cite[Lemma 2.3(v)]{ks}), and thus (for a given lift of the 
initial point) to a unique path in $\widetilde{X^{lpc}}$.
\end{pf}
        
Let us show how \hyperref[thm1]{Theorem \ref*{thm1}} can be applied to  at least \hl{ two special classes of spaces, the first of which is a 
a class of convergent
$Y$-spaces, but the second is not and could not be treated by
by our \hyperref[thm2]{Theorem \ref*{thm2}}}.
\begin{cor}\label{negcurvinj}      
\begin{itemize}
\item[(a)]  
If $Y$ is a \hl{closed} Riemannian manifold of negative sectional curvature, and $X$ is a convergent $Y$-space having only one identification point $y$, 
then $$\iota_k\colon H_k(X;\R)\to {\mathcal{H}}_k(X)$$
is injective in degrees $k\ge 2$.
\item[(b)] Let $(F_i)_{i\in\N}$ be a countable collection of
closed surfaces with a Riemannian metric of constant curvature $-1$.
Assume, each of the surfaces $F_i$ has a simple closed geodesic $l_i$, and each of them
has the same length $\ell$. On 
$\bigsqcup_{  i \in \N}  F_i$
identify these geodesics $l_i$ by isometries to one closed
geodesic  $l$, and then
put on the resulting point-set $X$ a first countable $T_4$-topology, where each 
subspace $F_i$ obtains as subspace-topology the one that it had before, 
and is closed as a subset of $X$. In addition we require that 
\begin{itemize}
\item[(i)]
$l \subset X$ has an open neighbourhood $U$, so that the intersection with
      each $F_i$ is of topological type $l_i \times (-1,+1)$, i.e., is an open band.
\item[(ii)] 
If $F_i$ gets equipped with its hyperbolic metric, 
     then for each open neighbourhood $V$ of $l$, there exists 
     $\eta > 0$,  so that the minimal width of all strips $l_i \times (-1,+1)$
     that are contained in $V$ with respect to this hyperbolic distance 
     is at least $\eta$, independently of $i$.
\end{itemize}
Then $\iota_k : H_k(X;\R) \to {\mathcal H}_k(X)$ is injective in degrees 
$k \ge 2$.
\end{itemize}    
\end{cor}
%\begin{pf}
{\sl The main ingredient for the proofs of both parts}\/ of this corollary
will be \hyperref[thm1]{Theorem \ref*{thm1}}. One of the assumptions of \hyperref[thm1]{Theorem \ref*{thm1}} is that the Gromov
norm is an actual norm, and showing this will require, that we precisely
understand the algebraic structure of these homology groups and the
geometric meaning their cycles. We will achieve this by proving the following
\hyperref[newEndProp14]{Proposition \ref*{newEndProp14}}. 

The proof of \hyperref[negcurvinj]{Corollary \ref*{negcurvinj}(a)} we will then be able to complete still
within this section. The proof of \hyperref[negcurvinj]{Corollary \ref*{negcurvinj}(b)} will also require one
lemma from the proof of \hyperref[thm2]{Theorem \ref*{thm2}}, and we will therefore finish the 
proof of \hyperref[negcurvinj]{Corollary \ref*{negcurvinj}(b)} in 
\hyperref[newSectSam]{Section \ref*{newSectSam}} 
after having proven
\hyperref[thm2]{Theorem \ref*{thm2}}. 
%\end{pf}
\begin{remrk}\label{newEndRem13a}
In preparing for \hyperref[newEndProp14]{Proposition \ref*{newEndProp14}},
we introduce some notation that we will use henceforth in discussing
the class of spaces $X = \bigcup^\infty_{i=1}F_i$ to be  considered in 
\hyperref[negcurvinj]{Corollary \ref*{negcurvinj}(b)}.
        If $\N$ enumerates our subsurfaces $F_i$ of $X$, we define 
$$       
    S:= \{ i \in \N \mid  l_i \hbox{ is non-separating in }F_i \}~~~  
    \hbox{ and } $$$$
           T := \{ i \in \N \mid  l_i \hbox{ is separating in } F_i \}.
$$           
      Naturally we have $S \cup T = \N$, but one of the sets $S$ or $T$ might
      be finite or even empty. If $T$ is non-empty, we pick some
      $t_0 \in T$. Each surface with $i \in T$ will get cut along $l_i$
      into two halves, that will be denoted by $F_i^+$ and $F_i^-$. 
      Both surfaces inherit their orientation from $F_i$, and concerning
      their relative fundamental cycles, we place the annotation $()^+$ and
      $()^-$ so, that $\partial [F_i^+] = [l]$, but $\partial [F_i^-] = -[l]$. \par
Apart from $T$
we      also introduce the notation $T^+$ and $T^-$,   where
      as sets of indices we have $T=T^+=T^-$, but with adding the upper indices
$()^\pm$,          we 
will be indicating, to which of the half-surfaces $F_i^+$ or $F_i^-$
we are referring. By using the squared version of the union-symbol
(``$\sqcup$") to denote disjoint union, while we naturally
have $\N = S \cup T = S \sqcup T$, we will also use constructions
such as $S \sqcup T^+ \sqcup (T^- \setminus \{t_0\})$.
\end{remrk}
\begin{prop}\label{newEndProp14}
\begin{itemize}
\item[(a)]
Let $X$ be a space satisfying the conditions 
of \hyperref[negcurvinj]{Corollary \ref*{negcurvinj}(a)}. Then for $k \ge 2$ we have 
$$H_k(X;\R) \cong  H_k(Y_\infty;\R) \oplus \bigoplus_{i=1}^\infty H_k(Y_i;\R)$$
      where the corresponding isomorphism is naturally induced by the
      system of projections $X \to Y_i$ (being the identity on $Y_i$
      and projecting the other $Y_j$-spaces to the attaching point). 
\item[(b)]        
      Let $X$ be defined as in \hyperref[negcurvinj]{Corollary \ref*{negcurvinj}(b)}. Then $H_k(X;\R) = 0$, if
       $k \ge 3$ and      
\begin{equation}
H_2(X,\R) = \bigoplus_{j\in S \sqcup T^+ \sqcup (T^- \setminus \{t_0\})} \R,
\label{eq:einFis}   
\end{equation}
       where the element ``1" in
      each $\R$-factor of \hyperref[eq:einFis]{ (\ref*{eq:einFis})} either 
      represents the fundamental cycle of a surface $F_j$ (for $j \in S$), 
      or the relative fundamental cycle of a half-surface 
      $F_j^\pm$ (for $j \in T^+$ or $T^-$).
\end{itemize}
\end{prop}
\noi{\it Remark: }  
Usually the Mayer-Vietoris sequence is considered to be the most
suitable tool when it comes to reducing homology groups of a given space 
to the homology groups of subspaces, whose union is the given space.
However, one should also be aware, that every subspace of the
earring space that just consists of finitely many rings has
by being isomorphic to a finite bouquet as its first $\ZZ$-homology
group a finite direct sum $\bigoplus_i \ZZ$, although the first homology group
of the earring space has by  \cite[Thm.1.3]{Eda-JLond}/\cite[Rem.3.1]{Eda-Fund2016}
a similar structure as our formula \hyperref[eq:afort1]{(\ref*{eq:afort1})}.
We therefore treat the Mayer-Vietoris sequence as a non-suitable tool when it comes
to understanding the passage from finite to infinite unions
of subspaces, and avoid to use it at all in the following argument:   
  
\begin{pf}
{\bf (a)} By \cite[Chapter 2, Theorem 10.6]{munkr} a compact manifold is 
triangulable, and with similar arguments as discussed 
in the proof of \hyperref[counta]{Lemma \ref*{counta}}, for the set 
of singular simplices in a convergent $Y$-space the possible accumulation 
of the various layers is irrelevant, and the result will be the same 
as if the topology is replaced by the $()^{lpc}$-topology. In addition,
since all our spaces are manifolds and have a disk-neighbourhood near
the attaching-point, the precise topology near the attaching point
(topological glueing vs.\ metric glueing) is irrelevant, because
all these topologies are homotopy equivalent, so 
that, again with the references given in the proof of \hyperref[counta]{Lemma \ref*{counta}},
we can actually use for computing the homology groups the same topology 
as when uniting the triangulations
of all $Y_i$-spaces ($i \in \N \cup \{\infty\}$) to one global 
triangulation of $X$. So we can use simplicial homology and apply it
to the resulting simplicial complex that must 
be homotopy equivalent to $X^{lpc}$ for computing the singular 
homology groups of $X$. In the context of simplicial complexes it is well-known
that the homology groups $H_k$ (apart from $H_0$) of a wedge-sum
of arbitrary cardinality can be computed as the direct sum
of the homology groups of the complexes that got united. 
While the argument for $k=1$ needs a little extra consideration,
for $k\ge2$ it is quite obvious, because the relevant chain groups
and boundary morphisms of a wedge sum are not in any way different from 
that of a disjoint union, and since for a disjoint union
all boundary relations remain inside the complexes that got united,
the direct-sum structure of the chain groups factors through
to the homology groups.
\par\noi{\bf  
(b)} We use analogous arguments as for (a) for reducing the computation
of this homology group to the computation of the homology group of a
simplicial complex that arises, when we triangulate each of our 
surfaces $F_i$ in such a way that the geodesics $l_i$ to be identified
are some edge-path in each of triangulations, and the subdivision
into edges agrees with the desired identification. Since that way
computing the homology groups reduces to simplicial techniques, it is
already clear that by lack of presence of higher-dimensional simplices 
all the homology groups for $k \ge 3$ must vanish. \par
   As for the possibility to obtain non-trivial 2-dimensional cycles:
Constructing cycles means already constructing non-trivial elements for the homology groups,
since no boundary coming from  
3-dimensional chains needs to be modded out.
Every cycle can only consist of finitely many simplices, in particular
of 2-simplices taken from at most finitely many of our $S$-surfaces,
$T^+$-half-surfaces and $T^-$-half-surfaces. Since in the interior 
of a triangulated surface each edge is only in the boundary 
of two 2-simplices, if a cycle should take one simplex from one
of our (half-)surfaces with a non-zero coefficient, it must take
all simplices from that (half-)surface with the same coefficient 
for maintaining the chance of having zero boundary. If we follow this 
rule for all (half-)surfaces appearing with non-zero coefficients,
and if we name these coefficients by $\alpha_i$ for the $S$-surfaces,
$\beta^+_i$ for the $T^+$-surfaces and $\beta^-_i$ for the
$T^-$-surfaces (only finitely many from these coefficient can
be non-zero), then 
the boundary of our prospective cycle $c$ will be
\begin{equation}
\partial c = \sum_{i \in T} (\beta^+_i - \beta^-_i)[l].
\label{eq:einStar}                
\end{equation}             
Hence we can take the point of view, as given when stating this 
\hyperref[newEndProp14]{Proposition \ref*{newEndProp14}(b)}, that we can, apart from the finiteness restriction,
choose the coefficients in our cycle in all (half-)surfaces
apart from $F^-_{t_0}$ freely, and the coefficient 
$\beta^-_{t_0}$ will then         
have to be chosen so that $\partial c=0$ and formula \hyperref[eq:einStar]{ (\ref*{eq:einStar})} determines the value
of $\beta^-_{t_0}$. Since we had to choose the simplices in each surface
so that all neighbouring simplices had to be chosen with the same coefficient, 
it is also clear, as claimed in the 
statement of this proposition, that the generators
of the $\R$-summands can be interpreted as fundamental cycles of the    
surfaces or as relative fundamental cycles for the half-surfaces.     
\end{pf}
\begin{prop}\label{newEndProp15}
Let $Y$ be a closed Riemannian manifold of negative sectional curvature.
Then the Gromov norm of every non-nullhomologous 
cycle $z \in H_k(Y;\R)$ is bounded away from zero.
\end{prop}

\begin{pf}
By a well-known argument from \cite[Section 1.1]{gro}, 
nontriviality of the Gromov norm is implied if we have
surjectivity of $H_b^k(Y)\to H^k(Y;\R)$ in degrees $k\ge 2$. Namely, for a homology class $\alpha$ let $\beta$ be a 
cohomology class with $\langle \beta,\alpha\rangle=1$. Then surjectivity of $H_b^k(Y)\to H^k(Y;\R)$ implies $\Vert\beta\Vert<\infty$, and from $1\le \Vert\beta\Vert\Vert\alpha\Vert$ we obtain $\Vert\alpha\Vert>0$.

So we have to prove surjectivity of $H_b^k(Y)\to H^k(Y;\R)$ in degrees $k\ge 2$. 

$Y$ as a closed manifold with negative curvature 
has an upper curvature bound $K < 0$. Its universal cover
$\widetilde{Y}$ is a simply connected complete geodesic manifold with upper 
curvature bound $K$. Therefore 
(cf.\ e.g.\ \cite[Thm.II.1A.6/Def.II.1.2]{BrHfl}) it satisfies
a local CAT($K$)-condition, and as a simply connected 
local CAT($K$)-manifold it also satisfies a global
CAT($K$)-condition by  the Cartan-Hadamard Theorem
(cf.\ e.g.\ \cite[Thm.II.4.1(2)]{BrHfl}). Any two points can be joined
by a geodesic, and the geodesics are unique, since the CAT($K$)-condition
rules out the existence of bigons.

This allows us to define a straightening process for simplices
via the lift to the universal covering of a simplex:
We call a $1$-simplex in $\widetilde{Y}$ straight if it coincides with the
geodesic connection of its endpoints.
Higher-dimensional straight simplices are then defined by successively taking straight cones over straight subsimplices as in 
\cite[Section 1.2]{gro}. 
This implies that every simplex $\sigma$ in $Y$ is homotopic (rel.\ vertices) to a unique straight simplex $str(\sigma)$ in $Y$. In particular we can straighten any cycle $c$ recursively by straightening its $k$-skeleton for $k=1,\ldots, dim(c)$. Dually this yields that any cocycle $c$ is cohomologous to the ``straightened" cocycle $c\circ str$.

The volume of straight simplices (of dimension $\ge 2$) in negatively curved $n$-manifolds is uniformly bounded (see 
\cite[Section 1.2]{gro} 
or \cite[Proposition 1]{iy}) 
by a constant $V(n,K)$ depending on the negative upper curvature bound $K$. 
This implies 
by \cite[Section 1.2, Theorem (C)]{gro} 
that $$\Vert c\Vert_\infty\le V(n,K)comass(\omega)$$ 
for a differential form $\omega$ representing $c$. Thus $c\circ str$ is a bounded cocycle for any cocycle $c$ in degree $\ge 2$. In particular, $H_b^k(Y)\to H^k(Y;\R)$ is surjective in degrees $k\ge 2$, a sufficient criterion
as mentioned at the beginning of this proof. 
\end{pf}

\begin{prop}\label{newEndProp16}
Let $X$ be a space satisfying the assumptions of \hyperref[negcurvinj]{Corollary \ref*{negcurvinj}(a)}
and  $k \ge 2$. Then the Gromov norm of every non-nullhomologous 
cycle $z \in H_k(X;\R)$ is bounded away from zero. 
\end{prop}
\begin{pf}
We realize our abstract cycle $z$ by a concrete chain of singular
simplices $c$ and have to show that the $l^1$-norm of $c$ is bounded away 
from 0 by parameters that only depend on the homology class of $z$.

First we want to show that, without loss of generality, we can reduce
to straight chains $c$ representing $z$. 
All our assertions about straightening processes from the
proof of \hyperref[newEndProp15]{Proposition \ref*{newEndProp15}}
remain true in this case. We just have to extend the definition of
straight simplices to the spaces considered here, where 
from \hyperref[newEndProp15]{Proposition \ref*{newEndProp15}}
we already have the definition of straightening their building blocks: 
\neuzl
We call a $1$-simplex in $X$ straight if its intersections with 
$x_0$ decompose it into 1-simplices, 
each of which is straight in the respective $Y_n$ 
(disregarding those intersections, where it remains in the subspace
$Y_n$).
Therefore the process from 
\hyperref[newEndProp15]{Proposition \ref*{newEndProp15}
can be applied to the segment of each simplex}. 
The straightening process can accidentally transform two singular
simplices into the same one, but not vice versa. Therefore
it cannot increase the $l_1$-norm. 

We hence may, without loss of generality, assume that $c$
is a straight cycle. From the proof of 
\cite[Lemma 5]{kue} 
we know that for every simplex in $X$ its straightening has at most one ``central simplex" (in the terminology of \cite{kue}) 
and that all other parts of the straightened simplex
are degenerate. I.e., with respect to what we have shown in \hyperref[newEndProp14]{Proposition \ref*{newEndProp14}(a)},
our straight (singular) cycle $c$ has only that kind of 
simplices, where each of them is clearly associated 
to one of the spaces $Y_i$, only. If this central 
part is lifted to the universal covering
$\widetilde{Y_i}$ of $Y_i$, it becomes a geometric
straight object whose $k$-dimensional volume can be measured, and will by 
\cite[Section 1.2]{gro} 
or \cite[Proposition 1]{iy} 
be uniformly bounded  by a constant $V(n,K)$ depending on the negative upper
curvature bound $K$ (which exists because the manifold is compact).
Let us now recall that our cycle $z$ by \hyperref[newEndProp14]{Proposition \ref*{newEndProp14}} decomposes
into components, that are associated to finitely many of the $Y_i$.
We will denote them, without loss of generality by 
$Y_1,Y_2,\ldots, Y_m$, and are with this notation not ruling out 
that one of them might be actually $Y_\infty$. Accordingly we
have 
$$z = \sum^m_{i=1} z_i \hbox{ with } z_i \in H_k(Y_i,\R) \subset H_k(X;\R).$$                                                           
On the level of chains our straight chain      
$c = \sum_{j\in J}\gamma_j\cdot\sigma_j $ 
analogously decomposes its finite index set 
$$    J = J_1 \sqcup J_2 \sqcup\ldots\sqcup J_m ,$$  
where $J_i$ comprises precisely those indices $j$ where
the straight simplex $\sigma_j$ takes its central part in $Y_i$. 
So we have
$$z \sim c = \sum^m_{i=1} \sum_{j\in J_i} \gamma_j\cdot \sigma_j~~~\hbox{and} 
$$$$z_i \sim c_i =  \sum_{j\in J_i} \gamma_j\cdot \sigma_j.$$
Although $\sigma_j$ with $j \in J_i$ may have degenerated parts outside $Y_i$,
we will use the same notation for the simplex $\in S_k(Y_i)$, where these
parts are cut off, so that $\sum_{j\in J_i} \gamma_j\cdot \sigma_j$
also defines a cycle $z_i \in H_k(Y_i;\R)$.

When we apply the principle that the Gromov norm of the homology class
is the infimum of the Gromov norm of the representing cycles,
we obtain the inequality  
$$0 \le \Vert z_i\Vert \le \sum_{j\in J_i}   |\gamma_j| .$$                      
When we now sum these inequalities over all $i \in \{1,...,m\}$, we get
on the right hand side $\Vert c\Vert_1$, i.e., the $l^1$-norm of our chain $c$,
and in the middle position    
$\sum^m_{i=1}\Vert z_i\Vert.$ 
So all it needs to finish        
this proof is to rule out that this sum might be zero. However, all the
$z_i$ are homology classes in $H_k(Y;\R)$, where in \hyperref[newEndProp15]{Proposition \ref*{newEndProp15}} we have shown,
that the Gromov norm is an actual norm. So our sum could only be zero, 
if all these homology classes would be zero, but since $z=\sum^m_{i=1} z_i$
                and $z$ was assumed to be non-nullhomologous, this is      
ruled out. Consequently this lemma is proven.      
\end{pf}                    
                     
\neuzl\hbox{~}\neuzl{\it            
Completing the proof of \hyperref[negcurvinj]{Corollary \ref*{negcurvinj}(a)}}  
\neuzl\hbox{~}\par
Let $n:= dim(Y)$.

We need to verify the assumptions of \hyperref[thm1]{Theorem \ref*{thm1}} in order to finish the
proof. The fact that the space $X$ that we constructed is a second
countable $T_1$-space we treat as obvious, the fact that 
the Gromov norm on all homology groups with $k \ge 2$ is an actual norm 
we have shown in \hyperref[newEndProp16]{Proposition \ref*{newEndProp16}}. From \hyperref[ExGenCov]{Lemma \ref*{ExGenCov}} we know, that
$X$ has a generalized universal covering space. Therefore, all it still
needs, is to show that this space is contractible. Therefore,
we need to some extent make up our mind, how this
covering space looks like:

With the arguments that are already used in the 
proof of 
\hyperref[newEndProp15]{Proposition \ref*{newEndProp15}},
the universal cover $\widetilde{Y}$ is a complete
Riemannian manifold satisfying a global CAT($K$)-condition
($K<0$).\par
Since any two points can be joined
by a geodesic, all geodesics emanating from the centre of some
$\R^n$-neighbourhood will reach all points of $\widetilde{Y}$. Since $\widetilde{Y}$ is 
completely geodesic (i.e., every geodesic runs till infinity),
and thanks to the global CAT($K$)-condition the geodesics emanating
from one point can neither close up nor intersect later, we conclude
that $\widetilde{Y}$ must actually be homeomorphic to $\R^n$. 

For proving that not just $Y$, but also $X$ has a contractible
generalized universal covering, we construct this space explicitly: \neuzl
Having shown that $\widetilde{Y}$ is $\R^n$, it is clear, that $Y$
is aspherical. Hence we can apply 
\hyperref[asphe]{Lemma \ref*{asphe}}
and conclude, 
that also our convergent $Y$-space $X$ will be aspherical. Similar as 
in the proof of 
\hyperref[ExGenCov]{Lemma \ref*{ExGenCov}},
a convergent $Y$-space has a generalized universal covering, which
is up to homotopy equivalence the classical universal covering. So, up
to homotopy equivalence our universal covering $\widetilde{X}$ is an
aspherical CW-complex, that in addition as universal covering has also
trivial $\pi_1$. Therefore this space is weakly contractible, and by the
possible application of the CW-topology and Whitehead's theorem therefore
also contractible. Having that way shown that $X$ has a contractible 
universal covering, actually we completed to check the assumptions of 
\hyperref[thm1]{Theorem \ref*{thm1}}
and can now use it, to conclude that for the spaces
as described in 
\hyperref[negcurvinj]{Corollary \ref*{negcurvinj}(a)}
the canonical homomorphism   $\iota_k\colon H_k(X;\R)\to {\mathcal{H}}_k(X)$   
will be injective for $k>0$.                
\hfill\blackbox
\par
The proof of \hyperref[negcurvinj]{Corollary \ref*{negcurvinj}(b)}
will be completed in \hyperref[newSectSam]{Section \ref*{newSectSam}}
that is devoted to that purpose. Apart from the construction
performed there, it also takes 
\hyperref[simplnu1]{Lemma \ref*{simplnu1}}/\hyperref[simplnu1po]{Remark 
\ref*{simplnu1po}} to complete the proof of this corollary,
but this corollary cannot directly be concluded from the 
proof of \hyperref[thm2]{Theorem \ref*{thm2}}~--- other than 
\hyperref[negcurvinj]{Corollary \ref*{negcurvinj}(a)},
which does not only also follow from \hyperref[thm2]{Theorem \ref*{thm2}},
but will be generalized in
\hyperref[hypex]{Section \ref*{hypex}}.
\begin{remrk}\label{newSectEndRem2.17}
Please note that that the strategy to prove 
\hyperref[negcurvinj]{Corollary \ref*{negcurvinj}(b)}
will in many aspects be analogous. However, we are lacking there the
analogue of 
\hyperref[counta]{Lemmas \ref*{counta}}--\hyperref[asphe]{\ref*{asphe}},
that we had been using in the last
paragraph of the just finished proof of 
\hyperref[negcurvinj]{Corollary \ref*{negcurvinj}(a)},
but that required that the identifications are taking place
along points. However, in the course of having
to prove 
\hyperref[negcurvinj]{Corollary \ref*{negcurvinj}(b)}
we have our identifications along
lines. Therefore, in order
to prove 
\hyperref[negcurvinj]{Corollary \ref*{negcurvinj}(b)},
we will be constructing the generalized universal
covering explicitly in 
\hyperref[newSectEndConstr4.3]{Construction \ref*{newSectEndConstr4.3}}
and conclude its contractibility
from this explicit construction. However, we wish to point out that 
in case of the proof of 
\hyperref[negcurvinj]{Corollary \ref*{negcurvinj}(a)}
an analogous explicit construction of 
the universal covering $\widetilde{X}$, based on glueing copies
of $\widetilde{Y}\approx\R^n$ along points with a similar tree-like
glueing structure would also have been possible, in order to end the
proof of 
\hyperref[negcurvinj]{Corollary \ref*{negcurvinj}(a)}.
\end{remrk}

\section{Proof of Theorem 2}\label{proof2}

\subsection{Examples and Neighbourhoods}\label{techlem}

According to a remark from the paragraph after the second
figure in the introduction, we start with an example that sheds light
on the relation between the properties ``semilocally simply connected"
and ``having a countable fundamental group".

\begin{expl}\label{afortExpl}
The purpose of this example is to show that a non-semilocally 
simply connected space can have a countable fundamental group. 
We start the construction of this space with the
harmonic archipelago (``$HA$"), which is a space that is embedded in 
three-space, and that at first sight looks like a disk 
with infinitely many humps built in, so that it will not    any longer
be homeomorphic to an ordinary disk. For a figure see   
\cite[Fig.4]{Omg-org}, \cite[Fig.1]{Omg-publ}, \cite[Fig.1(left)]{Arch}
or \cite{Br-Blogpost}.
For our purposes we do not compactify the space by adding
a vertical line as accumulation line 
for the humps; it is irrelevant for the fundamental group, either.

This space does not only satisfy to be non-semilocally simply 
connected, but also the stronger property of being 
not homotopically Hausdorff (i.e., there exists a non-trivial
homotopy class that can be represented in an arbitrarily small
neighbourhood, cf.\ \cite[Def.5.2]{CC3} and \cite[Def.2.9]{fz}): 
Any closed path can be homotoped into the
earring subspace of any of our pictures, and since only finitely many 
rings are sticking out of an arbitrary small neighbourhood of the accumulation
point, the passages of our paths via these rings can be homotoped via 
finitely many humps into our neighbourhood.

The first singular homology group of this space has been computed
in \cite[Thm.1.2]{KaRe12}.
This theorem conversely relies on \cite[Thm.1.3]{Eda-JLond}, where
a gap of the proof has been reported, but this gap has apparently
been mended in \cite[Rem.3.1]{Eda-Fund2016}.                       
The following formula results by 
substituting the last formula from \cite[Proof of Prop.2.4]{KaRe12}
into the formula from \cite[Thm.1.2]{KaRe12}:
\begin{equation} 
H_1(HA;\ZZ) \cong \bigl( \prod_{p~{\rm prime}} A_p \bigr) \oplus
\bigl(\sum_{i\in{\mathfrak{c}}} \Q\bigr)  \label{eq:afort1}
\end{equation}
where $A_p$ is the $p$-adic completion of the direct sum of $p$-adic integers
$\bigoplus_{\mathfrak{c}}\J_p$, and $\mathfrak{c}$ is the continuum cardinal.

With this observation we construct the following homomorphism

\begin{equation}
\pi_1(HA) \to (\pi_1(HA))^{\rm ab}= H_1(HA;\ZZ) \cong \bigl( \prod_{p~{\rm prime}}A_p \bigr) \oplus
\bigl(\sum_{i\in{\mathfrak{c}}} \Q\bigr) \to \Q  
\label{eq:afort2}
\end{equation}
where the first homomorphism is abelianization, and the 
last homomorphism is based on choosing some basis and
then projecting all but one basis-vector from the $\Q$-summand to zero.

As well obtaining the representation of the homology group     
in the form of formula \hyperref[eq:afort1]{(\ref{eq:afort1})}
as finding a basis for the $\Q$-part rely on the axiom of choice. 

However, having some surjective homomorphism $\pi_1(HA) \to \Q$, we can
attach disks to our space $HA$, so that the entire kernel of that
homomorphism becomes trivial and the fundamental group becomes
accordingly $\Q$, in particular countable. The definition of the topology
of this infinite attaching process can be analogously handled
as the construction of the space $RX$ in \cite[Def.7]{Spnc}, where 
the set of loops can be chosen as some group-theoretical generating
set for the kernel of homotopy classes of paths, and for each 
algebraic generator a path-representative
remaining in the earring subspace is to be chosen;
the distinction between  A-segments and B-segments as in 
\cite[Def.6]{Spnc} is not 
necessary.

Due to this attaching process the property that 
any homotopy class can be represented in the earring subspace is not 
overturned, and also not, that it can be actually represented
in each subspace of the earring subspace that only comprises the rings
from a certain index onwards. So the space is still
non-homotopically Hausdorff, but the fundamental group
is only $\Q$, i.e. countable.

\end{expl}

The following remark estimates the prospect of obtaining an
example as \hyperref[afortExpl]{Example \ref*{afortExpl}}
without having to involve the axiom of choice:

\begin{remrk}\label{afortRem} 

Assume that we have a space that is metric (or at least satisfies
the first axiom of countability). If it is not
semilocally simply connected, it allows to construct an essential map 
of the earring space into it according to the following method:
Take a countable neighbourhood-basis of a point, where it is not semilocally
simply connected. Turn this basis into a descending sequence of neighbourhoods 
by performing finite intersections. Being not semilocally simply connected
means that there must exist a non-trivial homotopy class 
in every arbitrary small neighbourhood. We pick a closed path
representing such a homotopy class in each neighbourhood of our descending
sequence of open neighbourhoods, and define our map from the earring space
to our space by putting a parametrisation of the path that is associated
to the $k$-th neighbourhood of our descending sequence on the $k$-th 
ring of the earring space. Although an infinite glueing, it will be
a continuous map on the entire earring space. 

The above construction mechanism unfortunately gives us little information,
how the homotopy classes that are taken on the different rings,
will be related to each other. If they should be really independent
(with respect to the group structure of the fundamental group
of our space), we can expect that our fundamental group will be
uncountable: Namely with the according independence, since  
there exist $2^\N$ different 0-1-sequences, we are able by composing
paths in the earring subspace, that, as prescribed by the 0-1-sequence,
run though or omit the according rings, to construct candidates
of uncountably many different homotopy classes. But only with independence
in a very strong sense one can conclude that they are really all different
homotopy classes in our space. Therefore it seems worth to look 
next at the other extreme case, that at first sight might appear
as a good track to obtaining a small fundamental group:
 
Assume that all homotopy classes, to which all rings of our
earring space map, are actually the same.
Then we can map an entire harmonic archipelago essentially into our space,
taking the homotopies that must exist between the images of the rings
to extend the map that we constructed already for the earring subspace
to the humps between two adjacent rings. However, by the references
that we gave in \hyperref[afortExpl]{Example \ref*{afortExpl}}, 
we know that also the harmonic archipelago 
has uncountable fundamental group. 

Both cases (as well the special case, where we could map an entire
harmonic archipelago, as the general case, where we could 
only map an earring space ``$HE$"), amount to the same question namely:
What kind of countable quotients does the earring group $\pi_1(HE)$
permit, with or without involving the axiom of choice?
Namely, if our non-semilocally simply connected space $X$ does
only have a countable fundamental group, then the essential
map of earring space into that space would imply that the
image of $\pi_1(HE)$ must be also countable, i.e.\ a countable quotient
of $\pi_1(HE)$. 

There are various combinatorial ways offered in literature
to understand the fundamental group of the earring space
(free $\sigma$-products in \cite{n-slender}, Big free groups in 
\cite{CC1}, topologist's 
product in \cite{Gr}, \cite{Hoj-solo} and \cite{Arch}, 
by order-type words  in \cite[Sect.1.2]{Omg-org} and 
\cite[Lemma 1.1]{Omg-publ}, as a proper subgroup of the canonical inverse
limit of finitely generated free groups in \cite[Corl.2.7]{MM}, word sequences 
in \cite{Z(Constr)} and \cite{Z(Speck)}, tame words in \cite{2Olg}). 
Intuitively the group can be understood as consisting of
all infinite products of ring-passages of all order-types, 
with the single restriction that every generator (that
in the language of word sequences will be called ``letter") can 
at most occur a finite number of times.

However an infinite arrangement of letters inside a word sequence
does not correspond to a group-theoretical product, and if a 
homomorphism respects the infinite product structure
coming from letter-arrangements (i.e., ``is letter-induced"), this
is a special property that not all homomorphisms fulfil. This is best
seen at the example of a space called ``the harmonic onion"
(``$H\Omg$", cf.\ \cite[Fig.3b)]{Omg-org}, \cite[Fig.6]{Omg-publ}): 
As for the harmonic archipelago, it is a space that is embedded
in three-space, and as for the harmonic archipelago it has 
a subspace homeomorphic to the earring space (in this case given by the
boundary curves of the disks), that apparently
generates the entire fundamental group. As for the harmonic archipelago
only finitely many homotopies through the set of disks
can be realized by a continuous homotopy.
The fundamental group of the harmonic onion 
is known to be uncountable, 
in particular since the harmonic archipelago was in  \cite[Prop.13]{Arch} been formally shown
to be homotopy equivalent to the mapping-cone from 
\cite[Fig.1(right)]{Arch}, which in turn is evidently homotopy equivalent
to the harmonic onion. Denoting the rings of the earring subspace 
of the harmonic archipelago by $c_i$ and the boundary curves of the 
harmonic onion by $b_i$, this homotopy equivalence corresponds
to the letter-induced isomorphism
\begin{equation}
   c_i c^{-1}_{i+1} \longleftrightarrow b_i 
\label{eq:afort1fis}
\end{equation}
on the level of fundamental groups between these two spaces. 
Consequently, 
the inclusion-induced homomorphism
\begin{equation}
\pi_1(HE) \to \pi_1(H\Omg) ~\hbox{   is a surjective homomorphism}
\label{eq:afort2fis}
\end{equation}
to an uncountable group that maps each letter and each
finite product of letters just to the neutral element.
          
This homomorphism also shows, how little, if we want to solve the           
task of constructing a non-trivial countable quotient of the
fundamental group of the earring space when mapping each ring
in a non-trivial way, we can rely on the values that we associate
on the rings itself. 
                                                    
By the way, the constructions in \cite{Arch} are only via the above-quoted
homomorphism  
\hyperref[eq:afort1fis]{(\ref{eq:afort1fis})}
associated to the original Archipelago Space, and could have
much more directly been associated to the harmonic onion-construction.

Our homomorphism as constructed above in that 
\hyperref[eq:afort2fis]{(\ref{eq:afort2fis})}
shows that it is possible
to map the group of word sequences so that all finite products
of letters map trivial, and only a subset of the 
infinite products maps to nontrivial elements.

However, the converse is not possible:  When treating 
only the place with a ``$\cdot$" as a group-theoretical
product, the following equation connects three elements of
$\pi_1(HE)$, two of which are given by infinite word sequences 
\begin{equation}    
  a_i\,a_{i+1}\,a_{i+2}\,a_{i+3} \ldots = a_i \mathrel{\cdot}
  a_{i+1}\,a_{i+2}\,a_{i+3} \ldots
\label{eq:afort3}
\end{equation}
Due to this and a lot of equations of that kind, 
it is clear that because we decided to map some finite
products and any letter to non-trivial elements, not all infinite 
products can map to zero, and all infinitely many nontrivial
values that we associated to our ring-passages, must have some
influence on how to map the infinite word sequences. In addition, a very
natural way to obtain a countable quotient of the earring group
(to take the letter-induced homomorphism that maps all
letters with index bigger than some $n$ to zero and the remaining letters
to the generators of some finitely generated group) contradicts
to our construction principle, that we have to map all letters
in a non-trivial way.

If we manage to construct a homomorphism from the group of word-sequences
to a countable group, it must have an uncountable kernel, that must have
a non-trivial intersection with the subset of infinite word 
sequences, and if there should be one non-zero-element in the image, 
at least one uncountable coset must map to a non-zero
element. 

The question is, whether this border between trivially mapping
and non-trivially mapping word-sequences can be arranged compatibly
with the group-arithmetic without involving the axiom of choice,
in particular, since the associated non-trivial values on 
all infinitely many letters influence this border.

Apparently the earring group is difficult to be taken apart,
the two results on its atomic properties \cite{Eda-atomic} 
and \cite{Nakamura-atomic} can be treated as an indication therefore.

As a special case  of mapping to countable 
groups, in \cite{Tlas} the homomorphisms from the earring
group to finite groups have been discussed: Any assignment on the letters
extends (although most of the assignments will not look as if 
the values on the initial strings of word sequences would converge),
but in order to extend the arrangement to infinite word sequences
\cite{Tlas} needed ultrafilter-methods motivated by non-standard analysis,
i.e., involved the axiom of choice.

Therefore it seems unlikely that such a map from the earring group 
to a countable group with the additional demand that each letter
has to be mapped in a non-trivial way can be arranged without
involving the axiom of choice, but conversely each  non-semilocally
simply connected space with countable fundamental group indirectly
defines such a map.
\end{remrk} 

We continue with an example that motivates the need of assumption
(v) for \hyperref[thm2]{Theorem \ref*{thm2}}:

\begin{expl}\label{CAstrich}

We consider a space $CA'$    obtained      by perturbing the definition 
of the space $CA$, that was described in the introduction:  
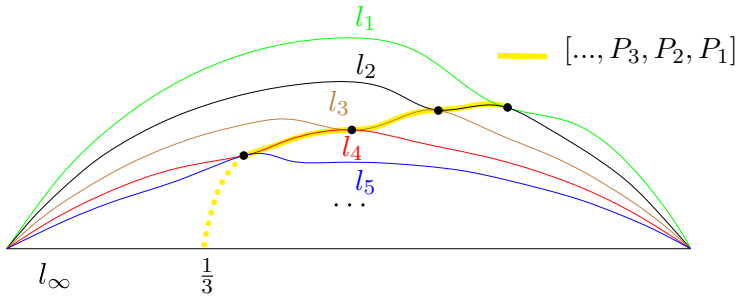
\neuzl

\def \globalscale {1.500000}
\begin{tikzpicture}[y=1cm, x=1cm, yscale=\globalscale,xscale=\globalscale, every node/.append style={scale=1}, inner sep=0pt, outer sep=0pt]
  \begin{scope}[line width=0.0106cm,miter limit=4.0,shift={(0.0048, 0.8452)}]
    \path[draw=black,line cap=butt,line join=miter,line width=0.0106cm,miter limit=4.0,cm={ 1.0,-0.0,-0.0,-1.0,(0.0053, 5.3391)}] (0.0, 4.5445).. controls (2.3398, 4.5445) and (3.6602, 4.5445) .. (6.0001, 4.5445);

  \end{scope}
  \begin{scope}[fill=black,line width=0.0106cm,miter limit=4.0,shift={(-1.0685, 0.8126)}]
    \begin{scope}[fill=black,line width=0.0106cm,miter limit=4.0,shift={(1.3543, -4.0262)}]
      \path[fill=black,line width=0.0106cm,miter limit=4.0] (0.0682, 4.7248).. controls (0.0682, 4.7248) and (0.0682, 4.7277) .. (0.0645, 4.7277).. controls (0.0587, 4.7277) and (0.0393, 4.7256) .. (0.0322, 4.7248).. controls (0.0302, 4.7248) and (0.0273, 4.7244) .. (0.0273, 4.7194).. controls (0.0273, 4.7165) and (0.0302, 4.7165) .. (0.0339, 4.7165).. controls (0.0467, 4.7165) and (0.0471, 4.7141) .. (0.0471, 4.712) -- (0.0463, 4.7066) -- (0.0128, 4.5747).. controls (0.012, 4.5718) and (0.0116, 4.5702) .. (0.0116, 4.5661).. controls (0.0116, 4.5508) and (0.0232, 4.5417) .. (0.0356, 4.5417).. controls (0.0442, 4.5417) and (0.0508, 4.547) .. (0.0554, 4.5566).. controls (0.0599, 4.5665) and (0.0633, 4.5818) .. (0.0633, 4.5822).. controls (0.0633, 4.5851) and (0.0608, 4.5851) .. (0.0599, 4.5851).. controls (0.0575, 4.5851) and (0.0571, 4.5838) .. (0.0566, 4.5801).. controls (0.0521, 4.5632) and (0.0471, 4.5475) .. (0.0364, 4.5475).. controls (0.0285, 4.5475) and (0.0285, 4.5557) .. (0.0285, 4.5594).. controls (0.0285, 4.5661) and (0.0289, 4.5673) .. (0.0302, 4.5723) -- cycle(0.0682, 4.7248);

    \end{scope}
  \end{scope}
  \begin{scope}[fill=black,line width=0.0106cm,miter limit=4.0,shift={(-1.0685, 0.8126)}]
    \begin{scope}[fill=black,line width=0.0106cm,miter limit=4.0,shift={(1.4329, -4.0657)}]
      \path[fill=black,line width=0.0106cm,miter limit=4.0] (0.1067, 4.595).. controls (0.0967, 4.6066) and (0.0938, 4.6095) .. (0.0868, 4.6144).. controls (0.0757, 4.6227) and (0.0637, 4.626) .. (0.0537, 4.626).. controls (0.0306, 4.626) and (0.0141, 4.6062) .. (0.0141, 4.5843).. controls (0.0141, 4.5628) and (0.0302, 4.5429) .. (0.0529, 4.5429).. controls (0.0785, 4.5429) and (0.0967, 4.5636) .. (0.1034, 4.5735).. controls (0.1129, 4.5619) and (0.1162, 4.559) .. (0.1228, 4.5541).. controls (0.1344, 4.5458) and (0.1459, 4.5429) .. (0.1563, 4.5429).. controls (0.179, 4.5429) and (0.1955, 4.5623) .. (0.1955, 4.5843).. controls (0.1955, 4.6062) and (0.1798, 4.626) .. (0.1567, 4.626).. controls (0.1311, 4.626) and (0.1133, 4.6049) .. (0.1067, 4.595) -- cycle(0.1124, 4.5884).. controls (0.1199, 4.6004) and (0.136, 4.619) .. (0.1583, 4.619).. controls (0.1774, 4.619) and (0.1906, 4.602) .. (0.1906, 4.5843).. controls (0.1906, 4.5669) and (0.1761, 4.5528) .. (0.1592, 4.5528).. controls (0.1418, 4.5528) and (0.1302, 4.5669) .. (0.1124, 4.5884) -- cycle(0.0972, 4.5805).. controls (0.0901, 4.5685) and (0.074, 4.5495) .. (0.0513, 4.5495).. controls (0.0322, 4.5495) and (0.0194, 4.5665) .. (0.0194, 4.5843).. controls (0.0194, 4.602) and (0.0339, 4.6157) .. (0.0508, 4.6157).. controls (0.0678, 4.6157) and (0.0798, 4.6016) .. (0.0972, 4.5805) -- cycle(0.0972, 4.5805);

    \end{scope}
  \end{scope}

  % FIRST ARC LABELS
  \begin{scope}[fill=black,line width=0.0106cm,miter limit=4.0,shift={(0.1128, 0.7997)}]
    \begin{scope}[fill=black,line width=0.0106cm,miter limit=4.0,shift={(2.9068, -1.7423)}]

    \end{scope}
  \end{scope}
  \begin{scope}[fill=black,line width=0.0106cm,miter limit=4.0,shift={(0.1128, 0.7997)}]
    \begin{scope}[fill=black,line width=0.0106cm,miter limit=4.0,shift={(2.9855, -1.7819)}]

      % L_1 as a text
      \draw (0.0616,4.662) node[green] {$l_1$};

    \end{scope}
  \end{scope}

  % LABELS OF THE SECOND PATH
  \begin{scope}[fill=black,line width=0.0106cm,miter limit=4.0,shift={(0.1184, 0.8078)}]
    \begin{scope}[fill=black,line width=0.0106cm,miter limit=4.0,shift={(2.9068, -2.1718)}]

    \end{scope}
  \end{scope}
  \begin{scope}[fill=black,line width=0.0106cm,miter limit=4.0,shift={(0.1184, 0.8078)}]
    \begin{scope}[fill=black,line width=0.0106cm,miter limit=4.0,shift={(2.9855, -2.2113)}]

     % L_2 as a text
     \draw (0.0616,4.66) node[black] {$l_2$};

    \end{scope}
  \end{scope}
  \begin{scope}[fill=black,line width=0.0106cm,miter limit=4.0,shift={(0.013, 0.9251)}]
    \begin{scope}[fill=black,line width=0.0106cm,miter limit=4.0,shift={(2.8515, -3.4564)}]
      \path[fill=black,line width=0.0106cm,miter limit=4.0] (0.0504, 4.5586).. controls (0.0504, 4.5661) and (0.0442, 4.5727) .. (0.0368, 4.5727).. controls (0.0289, 4.5727) and (0.0227, 4.5661) .. (0.0227, 4.5586).. controls (0.0227, 4.5508) and (0.0289, 4.5446) .. (0.0368, 4.5446).. controls (0.0442, 4.5446) and (0.0504, 4.5508) .. (0.0504, 4.5586) -- cycle(0.0504, 4.5586);

     % L_3 as a text
     \draw (0.0616,5.46) node[brown] {$l_3$};

    \end{scope}
  \end{scope}
  \begin{scope}[fill=black,line width=0.0106cm,miter limit=4.0,shift={(0.013, 0.9251)}]
    \begin{scope}[fill=black,line width=0.0106cm,miter limit=4.0,shift={(2.9687, -3.4564)}]
      \path[fill=black,line width=0.0106cm,miter limit=4.0] (0.0504, 4.5586).. controls (0.0504, 4.5661) and (0.0442, 4.5727) .. (0.0368, 4.5727).. controls (0.0289, 4.5727) and (0.0227, 4.5661) .. (0.0227, 4.5586).. controls (0.0227, 4.5508) and (0.0289, 4.5446) .. (0.0368, 4.5446).. controls (0.0442, 4.5446) and (0.0504, 4.5508) .. (0.0504, 4.5586) -- cycle(0.0504, 4.5586);

     % L_4 as a text
     \draw (0.0616,5.06) node[red] {$l_4$};

    \end{scope}
  \end{scope}
  \begin{scope}[fill=black,line width=0.0106cm,miter limit=4.0,shift={(0.013, 0.9251)}]
    \begin{scope}[fill=black,line width=0.0106cm,miter limit=4.0,shift={(3.0857, -3.4564)}]
      \path[fill=black,line width=0.0106cm,miter limit=4.0] (0.0504, 4.5586).. controls (0.0504, 4.5661) and (0.0442, 4.5727) .. (0.0368, 4.5727).. controls (0.0289, 4.5727) and (0.0227, 4.5661) .. (0.0227, 4.5586).. controls (0.0227, 4.5508) and (0.0289, 4.5446) .. (0.0368, 4.5446).. controls (0.0442, 4.5446) and (0.0504, 4.5508) .. (0.0504, 4.5586) -- cycle(0.0504, 4.5586);

     % L_5 as a text
     \draw (0.0616,4.76) node[blue] {$l_5$};

    \end{scope}
  \end{scope}

 % Yellow line

    \path[draw=yellow,line cap=butt,line join=miter,line width=0.0756cm,miter limit=4.0, line cap=round, dash pattern=on 0pt off 4pt] (2.0821, 2.4515).. controls (1.9988, 2.3815) and (1.9307, 2.2933) .. (1.8841, 2.195).. controls (1.8412, 2.1046) and (1.8167, 2.0068) .. (1.7916, 1.91).. controls (1.776, 1.8499) and (1.7602, 1.7897) .. (1.7503, 1.7285).. controls (1.7449, 1.695) and (1.7413, 1.6613) .. (1.7394, 1.6275);

  \path[draw=yellow,line cap=butt,line join=miter,line width=0.0756cm,miter limit=4.0] (4.4219, 2.8853).. controls (4.3566, 2.9108) and (4.286, 2.9223) .. (4.216, 2.9189).. controls (4.1428, 2.9154) and (4.0715, 2.896) .. (3.9992, 2.8843).. controls (3.9559, 2.8773) and (3.9122, 2.873) .. (3.8685, 2.8701).. controls (3.7975, 2.8652) and (3.7261, 2.8637) .. (3.6558, 2.8528).. controls (3.5744, 2.8401) and (3.4954, 2.8151) .. (3.4192, 2.7839).. controls (3.3533, 2.757) and (3.2887, 2.7252) .. (3.2194, 2.7091).. controls (3.1528, 2.6936) and (3.0838, 2.6929) .. (3.0155, 2.6902).. controls (2.9495, 2.6877) and (2.8836, 2.6832) .. (2.818, 2.6754).. controls (2.7353, 2.6656) and (2.6531, 2.6505) .. (2.5734, 2.6262).. controls (2.5049, 2.6052) and (2.4386, 2.5775) .. (2.3708, 2.5542).. controls (2.3089, 2.5331) and (2.2458, 2.5156) .. (2.1843, 2.4936).. controls (2.1502, 2.4814) and (2.1165, 2.4678) .. (2.0835, 2.4529);

   \draw (1.77,1.4) node[black] {$\frac{1}{3}$};

  \path[draw=green,line cap=butt,line join=miter,line width=0.0106cm,miter limit=4.0] (0.01, 1.6398).. controls (0.0884, 1.7519) and (0.1663, 1.8633) .. (0.2451, 1.9732).. controls (0.3239, 2.0831) and (0.4035, 2.1911) .. (0.4952, 2.2986).. controls (0.5869, 2.4062) and (0.6908, 2.5133) .. (0.8248, 2.6304).. controls (0.9587, 2.7474) and (1.1226, 2.8743) .. (1.2883, 2.9816).. controls (1.454, 3.0889) and (1.6214, 3.1765) .. (1.7945, 3.2504).. controls (1.9675, 3.3243) and (2.1484, 3.3851) .. (2.3483, 3.4295).. controls (2.5483, 3.474) and (2.7685, 3.5025) .. (2.9907, 3.5041).. controls (3.2129, 3.5058) and (3.437, 3.4807) .. (3.6491, 3.3613).. controls (3.8611, 3.2418) and (4.061, 3.0281) .. (4.2871, 2.9267).. controls (4.5133, 2.8253) and (4.7657, 2.8363) .. (4.9907, 2.7436).. controls (5.2158, 2.6509) and (5.4135, 2.4546) .. (5.568, 2.2754).. controls (5.7225, 2.0961) and (5.8339, 1.934) .. (5.9004, 1.8309).. controls (5.9669, 1.7278) and (5.9885, 1.6838) .. (6.0101, 1.6398);

  \path[draw=black,line cap=butt,line join=miter,line width=0.0106cm,miter limit=4.0] (0.01, 1.6398).. controls (0.085, 1.732) and (0.1601, 1.8243) .. (0.2369, 1.9088).. controls (0.3137, 1.9934) and (0.3924, 2.0702) .. (0.4856, 2.157).. controls (0.5789, 2.2439) and (0.6866, 2.3406) .. (0.7981, 2.4236).. controls (0.9096, 2.5066) and (1.0266, 2.5769) .. (1.1461, 2.6425).. controls (1.2656, 2.7081) and (1.3886, 2.7697) .. (1.5445, 2.8314).. controls (1.7005, 2.8932) and (1.8893, 2.9552) .. (2.0876, 3.0038).. controls (2.2859, 3.0524) and (2.4961, 3.088) .. (2.7071, 3.1064).. controls (2.9182, 3.1247) and (3.1312, 3.126) .. (3.3037, 3.0681).. controls (3.4761, 3.0102) and (3.6103, 2.891) .. (3.7673, 2.8679).. controls (3.9243, 2.8448) and (4.1051, 2.9166) .. (4.2648, 2.9105).. controls (4.4245, 2.9045) and (4.5628, 2.8206) .. (4.6934, 2.7425).. controls (4.8241, 2.6644) and (4.9469, 2.5921) .. (5.0622, 2.5186).. controls (5.1775, 2.445) and (5.2874, 2.3687) .. (5.3794, 2.2959).. controls (5.4714, 2.2231) and (5.5466, 2.153) .. (5.6489, 2.0441).. controls (5.7513, 1.9352) and (5.8807, 1.7875) .. (6.0101, 1.6398);

  \path[draw=brown,line cap=butt,line join=miter,line width=0.0106cm,miter limit=4.0] (0.01, 1.6398).. controls (0.0901, 1.715) and (0.1702, 1.7902) .. (0.2424, 1.8518).. controls (0.3146, 1.9135) and (0.3789, 1.9617) .. (0.4467, 2.0092).. controls (0.5145, 2.0566) and (0.5859, 2.1034) .. (0.6886, 2.1663).. controls (0.7913, 2.2292) and (0.9252, 2.3082) .. (1.0489, 2.3719).. controls (1.1725, 2.4357) and (1.288, 2.485) .. (1.4015, 2.5299).. controls (1.515, 2.5748) and (1.6277, 2.6158) .. (1.7607, 2.6543).. controls (1.8936, 2.6929) and (2.0468, 2.729) .. (2.1787, 2.7558).. controls (2.3107, 2.7826) and (2.4215, 2.8) .. (2.5646, 2.7755).. controls (2.7077, 2.751) and (2.8832, 2.6846) .. (3.0442, 2.6853).. controls (3.2052, 2.6859) and (3.3518, 2.7535) .. (3.4737, 2.803).. controls (3.5955, 2.8526) and (3.6944, 2.8846) .. (3.8292, 2.8552).. controls (3.964, 2.8258) and (4.1356, 2.7352) .. (4.2955, 2.6635).. controls (4.4554, 2.5918) and (4.6037, 2.5389) .. (4.7413, 2.4803).. controls (4.8789, 2.4216) and (5.0073, 2.3565) .. (5.1402, 2.28).. controls (5.2731, 2.2034) and (5.4112, 2.1152) .. (5.5219, 2.0386).. controls (5.6327, 1.9619) and (5.7161, 1.8968) .. (5.7929, 1.8323).. controls (5.8697, 1.7678) and (5.9399, 1.7038) .. (6.0101, 1.6398);

  \path[draw=red,line cap=butt,line join=miter,line width=0.0106cm,miter limit=4.0] (0.01, 1.6398).. controls (0.1445, 1.7292) and (0.279, 1.8186) .. (0.4249, 1.9028).. controls (0.5708, 1.9869) and (0.728, 2.0659) .. (0.9083, 2.1443).. controls (1.0887, 2.2228) and (1.2921, 2.3007) .. (1.5135, 2.3534).. controls (1.7349, 2.406) and (1.9743, 2.4334) .. (2.1324, 2.473).. controls (2.2905, 2.5126) and (2.3674, 2.5644) .. (2.5083, 2.6108).. controls (2.6491, 2.6571) and (2.8539, 2.6979) .. (3.0417, 2.6914).. controls (3.2294, 2.6849) and (3.4002, 2.6311) .. (3.603, 2.5839).. controls (3.8058, 2.5367) and (4.0385, 2.4965) .. (4.2714, 2.4382).. controls (4.5043, 2.38) and (4.7364, 2.3039) .. (4.9261, 2.2306).. controls (5.1158, 2.1573) and (5.2619, 2.0874) .. (5.3849, 2.023).. controls (5.5079, 1.9586) and (5.606, 1.9006) .. (5.7062, 1.8376).. controls (5.8065, 1.7747) and (5.9083, 1.7073) .. (6.0101, 1.6398);

  \path[draw=blue,line cap=butt,line join=miter,line width=0.0106cm,miter limit=4.0] (0.01, 1.6398).. controls (0.1787, 1.7235) and (0.3473, 1.8071) .. (0.5129, 1.8818).. controls (0.6785, 1.9564) and (0.841, 2.0221) .. (1.0154, 2.0805).. controls (1.1898, 2.1389) and (1.3761, 2.1901) .. (1.5801, 2.2688).. controls (1.7841, 2.3475) and (2.0059, 2.4536) .. (2.1567, 2.4768).. controls (2.3076, 2.5) and (2.3876, 2.4404) .. (2.4994, 2.4153).. controls (2.6113, 2.3901) and (2.755, 2.3994) .. (2.8916, 2.403).. controls (3.0281, 2.4066) and (3.1575, 2.4045) .. (3.2818, 2.3988).. controls (3.4061, 2.3931) and (3.5255, 2.3838) .. (3.649, 2.3725).. controls (3.7725, 2.3613) and (3.9003, 2.348) .. (4.0282, 2.3262).. controls (4.1561, 2.3044) and (4.2794, 2.2752) .. (4.391, 2.2482).. controls (4.5026, 2.2211) and (4.6, 2.1969) .. (4.7092, 2.1681).. controls (4.8185, 2.1393) and (4.9396, 2.1061) .. (5.0584, 2.0661).. controls (5.1772, 2.0261) and (5.2938, 1.9794) .. (5.3953, 1.9381).. controls (5.4968, 1.8967) and (5.5832, 1.8607) .. (5.6832, 1.8119).. controls (5.7831, 1.7631) and (5.8966, 1.7014) .. (6.0101, 1.6398);

  % Points P_1, P_2, etc

  \path[draw=black,fill=black] (3.7968, 2.8615) circle (0.033cm);

  \path[draw=black,fill=black] (3.038, 2.6899) circle (0.033cm);

  \path[draw=black,fill=black] (2.0897, 2.4634) circle (0.033cm);

  \path[draw=black,fill=black] (4.404, 2.8887) circle (0.033cm);

  \path[draw=yellow,fill=black,line width=0.0756cm] (4.3157, 3.3465) -- (4.7519, 3.3384);

  \node[text=black,anchor=south west,line width=0.0132cm] (text3) at (4.8973, 3.2603){$[..., P_3, P_2, P_1]$};

\end{tikzpicture}

\neuzl
$l_\infty$ is 
identically used for $CA$ and $CA'$, and is assumed to be embedded 
into $\R^2$ as the straight line segment $[(0,0),(1,0)]$. However
in order to obtain $CA'$ we change the definition of $l_1$
in the neighbourhood of the $x$-coordinate $2/3$, and for each other 
$l_n\;(n\neq 1,\infty)$ in the neighbourhood of the $x$-coordinates 
${1\over3}+{1\over3n}$ and ${1\over3}+{1\over3(n+1)}$, in such a way that
$l_n$ and $l_{n+1}$ obtain a joint tangent point 
$P_n$ at the $x$-coordinate ${1\over3}+{1\over3n}\,,\;\;\forall n\in\N$.
These tangent points \hl{are treated as 
parts of the intersections $l_n \cap l_{n+1}$ and form
intersection components},
and in $CA'$ the CW-structure of each $l_n$ consists of two 
or three 1-cells. $[P_{n+1},P_n]$ is a curve-segment in $l_{n+1}$,
and the concatenation of curve-segments 
$[\ldots,P_3,P_2,P_1]$, \hl{that is drawn with highlighing the
background like a thickened curve}
%marked in yellow in the above figure},
can be parametrized as a simple path $\sigma : (0,1] \to CA'$
such that $\sigma(1/n) = P_n$. The definition of this path 
can be continuously extended by letting $\sigma(0) := (1/3,0) \in 
l_\infty \subset  CA'$. 

Condition (v) of \hyperref[thm2]{Theorem \ref*{thm2}} is needed to rule out
that a path can in that way jump into the interior of a cell
of one of our CW-complexes, see
\hyperref[simplnu1]{Lemma \ref*{simplnu1}}. $CA'$ fulfils all other
assumptions of \hyperref[thm2]{Theorem \ref*{thm2}}.
\end{expl}  

Recall that any subcomplex of a CW-complex 
has deformation retractible neighbourhoods (\cite[Theorem II.6.1]{lw}). 
However, the following example
shows that the fact that we require our spaces to be coverable 
by CW-complexes does not imply that an intersection
complex must necessarily have such a deformation retractible
neighbourhood, and therefore this needs to be explicitly demanded
by condition (iii) of \hyperref[thm2]{Theorem \ref*{thm2}}. 

\begin{expl}\label{sin-touch} 
Let $$Y_0 :=  \{(x,y) \in \R^2 \mid y = x\cdot\sin(1/x),\;\;
x \in (0,1]\} \cup \{(0,0)\},$$
and for $i \in \N $ we first choose a sequence $a_i \in (1,2]$
so that $(a_i)_{i\in\N}$ is a strictly monotonically  
decreasing sequence converging to 1, and then let
$$Y_i :=  \{(x,y) \in \R^2 \mid y = a_i\cdot x,\;\;
x \in [0,1]\}\,,\;\; X :=  \bigcup^\infty_{\hl{i=0}} Y_i. $$
                 Each of the spaces $Y_i$ $(i \in \N_0)$
is homeomorphic to a line segment, and      ${(0,0)}$ is 
                            the only intersection point.
That way clearly      $X$ satisfies all conditions
of \hyperref[thm2]{Theorem \ref*{thm2}}, except possibly
the existence of a 
deformation retractible neighbourhood for $\{(0,0)\}$. 

However, let us first assume 
that there exists such a neighbourhood $ U_0$
together with the corresponding deformation retraction
$r_t : U_0 \to U_0\,,\;\; r_1 = {\rm id}|_{U_0}\,,\;\; r_0(U_0)=\{(0,0)\}$.
Since $Y_0$ has infinitely many tangent points with the diagonal
of the first quadrant which are accumulating at $(0,0)$, we are not losing
generality if we assume that $U_0 \owns ({2\over\pi},{2\over\pi}) \in Y_0$
and $U_0 \owns ({2\over\pi},{2a_i\over\pi})\in Y_i$ for $i \ge
i_0$ with some $i_0 \in \N$. Now
$ t \buildrel w\over\mapsto
r_t(({2\over\pi},{2\over\pi}))$ must be a path, which connects
$(0,0)$ to $({2\over\pi},{2\over\pi})$. Such a path must pass through
an  entire           segment of $Y_0$, in particular it must pass
through $({1\over\pi},0)$. Let $ t_1 :=
\min\{ t \in(0,1) \mid w([t,1]) \subset (Y_0 \cap \{ (x,y) \in \R^2 \mid x \ge
1/ \pi \})\}$. We now cover each point of the trace of $w|_{[t_1,1]}
$ by an open set
from $X$, small enough so as to avoid  $(0,0)$ and the segments
of $Y_i$  with $x$-values close to $1\over\pi$. 
We conclude, using the continuity of $r_t$, that the
complete preimages of these covering sets must give a covering
of $\{({2\over\pi},{2\over\pi})\}\times[t_1,1] \subset
U_0\times[0,1]$. Inside  each of these covering sets we find a smaller
covering set of product type $U_s \times (s-\eps_s,s+\eps_s)$. 
Compactness of $\{({2\over\pi},{2\over\pi})\}\times[t_1,1]$ allows us to
reduce this covering to a finite one, and, denoting the first components
of this finite subcovering by $ U_{s_1}\,,\ldots,\,U_{s_n}$, their
intersection  
$\bigcap^n_{i=1}U_{s_i}$ 
is still an open neighbourhood of $({2\over\pi},{2\over\pi}) \in X$.
Fix $i$ with 
$({2\over\pi},{2a_i\over\pi}) \in \bigcap^n_{i=1}U_{s_i}$
and consider the
path $w': [t_1,1] \to X,\;\;
t \mapsto r_t(({2\over\pi},{2a_i\over\pi}))$.
This path connects a point close to $({1\over\pi},0)
$ with $({2\over\pi},{2a_i\over\pi})$ inside $X$. Now
all points $\in X$  
close to              $({1\over\pi},0)$        belong to $Y_0$.
By the topology of $X$ each path            between such a point and 
$({2\over\pi},
{2a_i\over\pi})$ 
must pass through 
$(0,0)$. 
But $w'$, being subordinated to the above constructed covering, 
cannot pass through $(0,0)$.
                 This contradiction             rules out
that $X$ allows to construct a deformation retractible neighbourhood
$U_0$  of $(0,0)$.

\end{expl}

In view of the two preceding examples, 
conditions (iii) and (v) of \hyperref[thm2]{Theorem \ref*{thm2}}
should be understood to mean that morally we require 
the topology of our spaces to be to a certain extent well-behaved 
in the direction transversal to the CW-complexes, 
because otherwise we seem not to obtain the wanted results.

\medskip   
\hl{Recall that the intersection components that we introduced  
in \hyperref[thm2]{Theorem \ref*{thm2}} may be drawn from 
up to countably many different CW-complexes $Y_i$. 
Intersection scenarios where the induced cell-subdivision
of an intersection component is infinite, can easily be arranged, 
but are ruled out by \hyperref[thm2]{Condition \ref*{thm2}(i)}.
Our main motive to exclude them is the first reduction step 
that we perform in \hyperref[reduction]{Subsect.\ref*{reduction}},
and where we wanted to keep our spaces $Y_i$ to have finite triangulation.
The following example has been included to shed a bit of light on,
what problems would have to be dealt with in case we would have
also permitted infinite intersection complexes}.

\begin{expl}\label{stattviExpl}
We construct two spaces $CA''_1$ and  $CA''_2$, and similarly as for 
our \hyperref[CAstrich]{Example \ref*{CAstrich}}, their
construction is again a variation of the
convergent arcs space. Both examples have in common, that we 
keep the line-segment $l_\infty$, but perturb the other line-segments 
for $x$-values $\in [0,{1\over3}+\eps)$ in that way, that 
we lower the $x$-coordinates of the  $l_i$ with odd $i$ and raise them
for the $l_i$ with even $i$, so that precisely for the $x$-coordinates
$\in [0,{1\over3}]$ each of the subspaces $l_i$ with even $i$ 
coincides with $l_{i-1}$. These changes, when applied to 
$CA$, already describe $CA''_1$.  Concerning $CA''_2$ there is
an additional addendum: we add between the points 
$(1/3,0)$ and $(1,0)$ an arc to the space that runs through
the area of points with negative $x$-values, and treat this arc
together with the line-segment $[(0,0),(1/3,0)]$ as an additional
subspace $l'_\infty$. 

%Beim n\"achsten Bild:  cyan --> 0.3   black, opacity=0.3  
%                       green--> 0.4   black, opacity=0.4
%                       blue ---> 0.7  black, opacity=0.7
%                       red   ---> 0.5 black, opacity=0.5
%                     ex-brown ----> 0.8 black, opacity=0.8

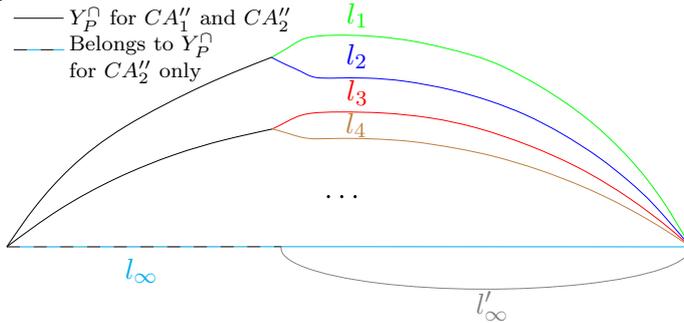
\def \globalscale {1.500000}
\begin{tikzpicture}[y=1cm, x=1cm, yscale=\globalscale,xscale=\globalscale, every node/.append style={scale=1}, inner sep=0pt, outer sep=0pt]
  \begin{scope}[line width=0.0106cm,miter limit=4.0,shift={(0.0059, 0.8195)}]

    \path[draw=cyan,fill,line width=0.0132cm,miter limit=4.0,cm={ 1.0,-0.0,-0.0,-1.0,(0.0053, 5.3391)}]
    (2.4062, 4.5432) -- (2.5937, 4.5432) -- (2.7812, 4.5432) -- (2.9687, 4.5432) -- (3.1562, 4.5432) -- (6.0001, 4.5445);

    \path[draw=black,fill=black,line width=0.0132cm,miter limit=4.0,cm={ 1.0,-0.0,-0.0,-1.0,(0.0053, 5.3391)},  dash pattern=on 5pt off 5pt]
    (0.0, 4.5445) -- (1.5, 4.5445) -- (2.4062, 4.5432);

    \path[draw=cyan,fill=black,line width=0.0132cm,miter limit=4.0,cm={ 1.0,-0.0,-0.0,-1.0,(0.0053, 5.3391)},  dash pattern=on 5pt off 5pt, dash phase=5pt]
    (0.0, 4.5445) -- (1.5, 4.5445) -- (2.4062, 4.5432);

  \end{scope}
  \begin{scope}[fill=black,line width=0.0106cm,miter limit=4.0,shift={(-0.0551, 0.956)}]
    \begin{scope}[fill=black,line width=0.0106cm,miter limit=4.0,shift={(2.8515, -3.4564)}]
      \path[fill=black,line width=0.0106cm,miter limit=4.0] (0.0504, 4.5586).. controls (0.0504, 4.5661) and (0.0442, 4.5727) .. (0.0368, 4.5727).. controls (0.0289, 4.5727) and (0.0227, 4.5661) .. (0.0227, 4.5586).. controls (0.0227, 4.5508) and (0.0289, 4.5446) .. (0.0368, 4.5446).. controls (0.0442, 4.5446) and (0.0504, 4.5508) .. (0.0504, 4.5586) -- cycle(0.0504, 4.5586);

    \end{scope}
  \end{scope}
  \begin{scope}[fill=black,line width=0.0106cm,miter limit=4.0,shift={(-0.0551, 0.956)}]
    \begin{scope}[fill=black,line width=0.0106cm,miter limit=4.0,shift={(2.9687, -3.4564)}]
      \path[fill=black,line width=0.0106cm,miter limit=4.0] (0.0504, 4.5586).. controls (0.0504, 4.5661) and (0.0442, 4.5727) .. (0.0368, 4.5727).. controls (0.0289, 4.5727) and (0.0227, 4.5661) .. (0.0227, 4.5586).. controls (0.0227, 4.5508) and (0.0289, 4.5446) .. (0.0368, 4.5446).. controls (0.0442, 4.5446) and (0.0504, 4.5508) .. (0.0504, 4.5586) -- cycle(0.0504, 4.5586);

    \end{scope}
  \end{scope}
  \begin{scope}[fill=black,line width=0.0106cm,miter limit=4.0,shift={(-0.0551, 0.956)}]
    \begin{scope}[fill=black,line width=0.0106cm,miter limit=4.0,shift={(3.0857, -3.4564)}]
      \path[fill=black,line width=0.0106cm,miter limit=4.0] (0.0504, 4.5586).. controls (0.0504, 4.5661) and (0.0442, 4.5727) .. (0.0368, 4.5727).. controls (0.0289, 4.5727) and (0.0227, 4.5661) .. (0.0227, 4.5586).. controls (0.0227, 4.5508) and (0.0289, 4.5446) .. (0.0368, 4.5446).. controls (0.0442, 4.5446) and (0.0504, 4.5508) .. (0.0504, 4.5586) -- cycle(0.0504, 4.5586);

    \end{scope}
  \end{scope}

  \path[draw=green,line cap=butt,line join=miter,line width=0.0106cm,miter limit=4.0] (2.3318, 3.2884).. controls (2.348, 3.2972) and (2.3643, 3.306) .. (2.3983, 3.326).. controls (2.4323, 3.346) and (2.484, 3.3771) .. (2.526, 3.4042).. controls (2.5681, 3.4313) and (2.6018, 3.4552) .. (2.6844, 3.4694).. controls (2.767, 3.4835) and (2.8992, 3.4883) .. (3.0534, 3.4836).. controls (3.2077, 3.479) and (3.3839, 3.4649) .. (3.5695, 3.4325).. controls (3.7552, 3.4001) and (3.952, 3.349) .. (4.0827, 3.3109).. controls (4.2133, 3.2729) and (4.2822, 3.2463) .. (4.3745, 3.2025).. controls (4.4667, 3.1587) and (4.5841, 3.097) .. (4.7131, 3.0181).. controls (4.842, 2.9392) and (4.9825, 2.8431) .. (5.1124, 2.7376).. controls (5.2422, 2.6321) and (5.3624, 2.5161) .. (5.4696, 2.3982).. controls (5.5768, 2.2803) and (5.6715, 2.16) .. (5.7597, 2.0289).. controls (5.8479, 1.8978) and (5.9296, 1.756) .. (6.0112, 1.6141);

  \path[draw=black,line cap=butt,line join=miter,line width=0.0106cm,miter limit=4.0] (0.0112, 1.6141).. controls (0.0895, 1.7309) and (0.1679, 1.8477) .. (0.2594, 1.9656).. controls (0.3509, 2.0835) and (0.4556, 2.2026) .. (0.5798, 2.3198).. controls (0.7041, 2.4371) and (0.8478, 2.5527) .. (1.0527, 2.6786).. controls (1.2576, 2.8046) and (1.5235, 2.9409) .. (1.7471, 3.0444).. controls (1.9706, 3.1478) and (2.1517, 3.2183) .. (2.3328, 3.2889);

  \path[draw=blue,line cap=butt,line join=miter,line width=0.0106cm,miter limit=4.0] (2.3309, 3.2904).. controls (2.3724, 3.2684) and (2.4139, 3.2464) .. (2.4577, 3.226).. controls (2.5016, 3.2056) and (2.5478, 3.1867) .. (2.5872, 3.1667).. controls (2.6265, 3.1466) and (2.6599, 3.1248) .. (2.7069, 3.1143).. controls (2.7539, 3.1039) and (2.815, 3.1046) .. (2.8941, 3.1058).. controls (2.9732, 3.107) and (3.0702, 3.1089) .. (3.1725, 3.1057).. controls (3.2747, 3.1026) and (3.3822, 3.0945) .. (3.5072, 3.0773).. controls (3.6323, 3.06) and (3.7748, 3.0336) .. (3.9025, 3.0034).. controls (4.0301, 2.9732) and (4.1429, 2.9392) .. (4.2495, 2.9031).. controls (4.356, 2.8669) and (4.4575, 2.8282) .. (4.5598, 2.7819).. controls (4.6621, 2.7355) and (4.7673, 2.6806) .. (4.8775, 2.6161).. controls (4.9876, 2.5517) and (5.1049, 2.4765) .. (5.2066, 2.403).. controls (5.3082, 2.3296) and (5.3948, 2.2574) .. (5.4693, 2.1934).. controls (5.5438, 2.1295) and (5.6061, 2.0739) .. (5.6944, 1.9788).. controls (5.7827, 1.8836) and (5.897, 1.7488) .. (6.0112, 1.6141);

  \path[draw=red,line cap=butt,line join=miter,line width=0.0106cm,miter limit=4.0] (2.3425, 2.6564).. controls (2.3822, 2.6768) and (2.422, 2.6971) .. (2.4618, 2.7175)(2.4615, 2.7172).. controls (2.5123, 2.7453) and (2.5632, 2.7735) .. (2.6413, 2.7886).. controls (2.7195, 2.8038) and (2.825, 2.806) .. (2.9412, 2.8062).. controls (3.0575, 2.8065) and (3.1844, 2.8047) .. (3.3177, 2.7969).. controls (3.4509, 2.789) and (3.5947, 2.7745) .. (3.7339, 2.7525).. controls (3.8731, 2.7304) and (4.013, 2.6999) .. (4.1561, 2.6592).. controls (4.2992, 2.6185) and (4.4491, 2.5665) .. (4.5989, 2.5072).. controls (4.7488, 2.448) and (4.9017, 2.3802) .. (5.0516, 2.3003).. controls (5.2015, 2.2204) and (5.3494, 2.1278) .. (5.4777, 2.0395).. controls (5.6061, 1.9512) and (5.7148, 1.8673) .. (5.8005, 1.7971).. controls (5.8861, 1.7269) and (5.9487, 1.6705) .. (6.0112, 1.6141);

  \path[draw=black,line cap=butt,line join=miter,line width=0.0106cm,miter limit=4.0] (0.0112, 1.6141).. controls (0.1613, 1.7271) and (0.3115, 1.8401) .. (0.4726, 1.9435).. controls (0.6337, 2.047) and (0.8059, 2.1408) .. (0.9782, 2.2224).. controls (1.1506, 2.304) and (1.3233, 2.3732) .. (1.4806, 2.4282).. controls (1.6379, 2.4833) and (1.7791, 2.5238) .. (1.9201, 2.5594).. controls (2.0611, 2.595) and (2.2014, 2.6255) .. (2.3418, 2.656);

  \path[draw=brown,line cap=butt,line join=miter,line width=0.0106cm,miter limit=4.0] (2.3406, 2.6569).. controls (2.3906, 2.643) and (2.4405, 2.629) .. (2.4914, 2.613).. controls (2.5422, 2.597) and (2.594, 2.579) .. (2.6601, 2.5724).. controls (2.7262, 2.5657) and (2.8066, 2.5705) .. (2.9, 2.5722).. controls (2.9935, 2.5739) and (3.1, 2.5726) .. (3.2133, 2.5686).. controls (3.3266, 2.5646) and (3.4468, 2.5579) .. (3.5723, 2.5437).. controls (3.6978, 2.5295) and (3.8286, 2.5077) .. (3.9456, 2.4864).. controls (4.0625, 2.4652) and (4.1686, 2.444) .. (4.292, 2.412).. controls (4.4154, 2.3801) and (4.5576, 2.3372) .. (4.6894, 2.2916).. controls (4.8212, 2.246) and (4.9443, 2.1971) .. (5.0432, 2.1552).. controls (5.1421, 2.1134) and (5.2176, 2.0783) .. (5.2985, 2.038).. controls (5.3794, 1.9977) and (5.4657, 1.9521) .. (5.5481, 1.9057).. controls (5.6306, 1.8593) and (5.7092, 1.8121) .. (5.7857, 1.7633).. controls (5.8623, 1.7146) and (5.9368, 1.6644) .. (6.0112, 1.6141);

  \path[draw=black,opacity=0.5,line width=0.0132cm,miter limit=4.0,dash pattern=on 0.0527cm off 0.00132cm] (6.0038, 1.6016)arc(0.0:60.8333:1.7931 and -0.3616)arc(60.8333:121.6667:1.7931 and -0.3616)arc(121.6667:182.5:1.7931 and -0.3616);

  \node[text=green,anchor=south west,line width=0.0556cm] (text1) at (2.9938, 3.5538){$l_1$};

  \node[text=blue,anchor=south west,line width=0.0556cm] (text2) at (2.9884, 3.193){$l_2$};

  \node[text=red,anchor=south west,line width=0.0556cm] (text3) at (2.9938, 2.87){$l_3$};

  \node[text=brown,anchor=south west,line width=0.0556cm] (text4) at (2.983, 2.58){$l_4$};

  \node[text=cyan,anchor=south west,line width=0.0556cm] (text4-3) at (1.0493, 1.295){$l_\infty$};

  \node[text=black,opacity=0.5,anchor=south west,line width=0.0556cm] (text4-7) at (4.1256, 0.9606){$l'_\infty$};

  \path[draw=black,fill=black,line width=0.0132cm] (0.0718, 3.6351) -- (0.5102, 3.6351);

  % Dashed with colours
  \path[draw=black,fill=black,line width=0.0132cm, dash pattern=on 5pt off 5pt] (0.0718, 3.3566) -- (0.5178, 3.3566);

  \path[draw=cyan,fill=black,line width=0.0132cm, dash pattern=on 5pt off 5pt, dash phase=5pt] (0.0718, 3.3566) -- (0.5178, 3.3566);

  \node[text=black,anchor=south west,line width=0.0132cm] (text7) at (0.5592, 3.5166){\scriptsize $Y^\cap_P$ for $CA_1''$ and  $CA_2''$ };

   \node[text=black,anchor=south west,line width=0.0132cm, align=left] (text8) at (0.5592, 3.1833){\scriptsize \parbox[c][8pt][c]{3cm}{Belongs to $Y^\cap_P$ \\ for $CA_2''$ only}};

\end{tikzpicture}

Of course, all points on all $l_i$ with $x$-coordinate $1/3$ (with maybe
$l_\infty$ in case of $CA''_1$ excepted) have to be treated 
as zero-cells with respect to the CW-structure of the $l_i$. 

Both spaces $CA''_i$ have two intersection components, 
$Q = \{(0,1)\} =: Y^\cap_Q$ is one of them, and the other intersection
component $Y^\cap_P$, \hl{that is drawn in black in the above 
picture}, consists of all segments on all $l_i$ for the
$x$-coordinates $[0,1/3]$ (including $l_\infty$ for
$CA''_2$, but excluding it for $CA''_1$). In this sense
both spaces do not satisfy \hyperref[thm2]{Condition \ref*{thm2}(i)}.
Conversely, both of them satisfy \hyperref[thm2]{Condition \ref*{thm2}(v)},
that cannot be violated for spaces having only two intersection
components. 

Actually each of our spaces does, apart from 
\hyperref[thm2]{Condition \ref*{thm2}(i)}, not satisfy another condition as
well: For $CA''_1$ it is impossible to construct a deformation retraction
of a relatively open neighbourhood of $Y^\cap_P$, 
because such a deformation would have to move
the points from $\{(x,0)\mid0< x \leq 1/3)\}$, which are on the other hand
limit points of points from the intersection component
that have to be fixed. 

This is a problem which disappears when looking at the space
$CA''_2$, but with $CA''_2$ there is another problem, namely,
that the cell-subdivision structure that $Y^\cap_P$
inherits, does not induce its subspace-topology, since 
the weak topology of the cells in $Y^\cap_P$ would not induce the
segment $\{(x,0)\mid x \leq 1/3)\}$ to be an accumulation line of the
corresponding segments of the other lines $l_i$.

However the following lemma 
shows, that this cannot happen as long as we demand 
\hyperref[thm2]{Condition \ref*{thm2}(i)}, i.e., that our intersection
components inherit only a finite cell-structure from
the CW-complexes $Y_i$:

\end{expl}

\begin{lemma}\label{stattviLemma}
The subspace-topology of each intersection component coincides
with the CW-topology that is induced by the finite cell-structure
that the intersection components inherit according to 
\hyperref[thm2]{Condition \ref*{thm2}(i)}.
\end{lemma}
\begin{pf}
We cannot rule out that infinitely many $Y_i$ intersect in our 
intersection component $Y^\cap_\nu$, but since the cell-subdivision
structure is finite, and one needs at most to consider one of the
$Y_i$-spaces to contribute one of the cells, it is clear that a finite
number of these intersecting $Y_i$-spaces suffices to give a complete 
closed covering of our intersection component. Now a set is open in the
CW-topology iff the complete preimage of the characteristic
map of all cells is open, and if we assume a set to satisfy this condition,
it is clear that it is open in each $Y_i$-space, because they were
assumed to be CW-complexes, i.e., they have the same topology 
as induced by their CW-structure. On the other hand, it is standard
general topology that, if a set is relatively open with respect
to all covering sets of a finite closed covering, then it is open
in the space as well. The conclusion in the reverse direction is even
more direct.  
\end{pf}
\begin{expl}\label{Mf-triang} 
Let $X$ be a finitely triangulated $n$-dimensional manifold.
It is the purpose of this example to study to what extent such a space,
if we try to interpret each $n$-dimensional simplex of the triangulation
together with its boundary faces                         as a CW-complex,
has the chance to fulfil the assumptions of our 
\hyperref[thm2]{Theorem \ref*{thm2}}:

Note that the properties of a triangulation ensure that all intersection
complexes              are contractible, and by standard techniques
of constructing regular neighbourhoods in triangulated spaces, they also
have deformation retractible neighbourhoods. Also non-accumulation
in the sense of           \hyperref[thm2]{condition \ref*{thm2}(v)}
is fulfilled by the existence of a finite triangulation. 
\hyperref[thm2]{Condition \ref*{thm2}(ii)} \hl{does not demand 
the contractibility of the intersection complexes $Y_I$, but
only for the intersection components $Y^\cap_\nu$, and in this
case we have just one intersection component which is the entire
$(n-1)$-skeleton $X^{(n-1)}$, that even in the case of a contractible
$X$ cannot be expected to be contractible}. 
Indeed 
\hyperref[Mf-triang]{Example \ref*{Mf-triang}}
fulfils all the conditions of
\hyperref[thm2]{Theorem \ref*{thm2}}, apart from  
\hl{having contractible intersection components}.

However with respect to the example given in \cite{zas1}, 
having a topology that is based on the accumulation of 
triangulable layers is definitely a too weak assumption 
in order to prove the desired injectivity of the canonical
homomorphism from singular homology as in our 
\hyperref[thm1]{Theorems \ref*{thm1}} and \hyperref[thm2]{ \ref*{thm2}}.
\end{expl}

\begin{expl}\label{Hyperbolas} 
Let $X$ be embedded into $\R^2$ such that \hl{%
$$ X := \{(x,y\} \mid(\Vert (x,y)\Vert \le 10)\wedge( x=0 \vee y=0 
\vee \exists n\in\N \hbox{ with }x \cdot y = \pm1/n  \} 
.$$}
In words: the space consists of segments of         
both coordinate axes and four sequences             
of segments of hyperbolas converging in the four quadrants against these
coordinate axes. Both coordinate axes and each of the hyperbolas is to 
be understood as a space $Y_i$ which has therefore the homeomorphism
type of a line segment, i.e., of a finite CW-complex. We do not need to 
decide on the precise CW-structure (i.e., into how many 1-cells
these line-segments are subdivided), apart from that 
$\{(0,0)\}$ must be a zero-cell as
\hl{the only intersection component}.

This           example       fulfils all conditions of
\hyperref[thm2]{Theorem \ref*{thm2}}. 
The deformation
retractible neighbourhood of  $\{(0,0)\}$
as required by \hyperref[thm2]{condition \ref*{thm2}(iii)}
consists of       segments on the 
$x$-axis and the $y$-axis around $(0,0)$. Observe, that this
neighbourhood is not open in the topology of $X$, but
\hyperref[thm2]{Theorem \ref*{thm2}} can be applied
to this example, because the neighbourhood is, at least,
relatively open.

We have chosen to go for such weak
assumptions, because we wanted examples like this 
still to be covered by our theorems. An open deformation-retractible
neighbourhood does in this case not exist: It would have to
contain some parts of some hyperbolas, which are not path-connected
to the origin.
\end{expl}

\begin{lemma}\label{dlaCorl3-3} 
Let $X$ be a convergent $Y$-space as defined in 
\hyperref[convy]{Definition \ref*{convy}}, and $P$ one of its identification
points, then $P$ is
\hl{an intersection component}, and
has an open deformation retractible neighbourhood as required
by condition (iii) of
\hyperref[thm2]{Theorem \ref*{thm2}}.

\end{lemma} 
\begin{pf} Let $P$ be an identification point of the convergent
$Y$-space $X$. In order to construct a deformation retractible
neighbourhood $U$, let us first use the CW-structure of $Y_\infty$
(see, e.g., \cite[Theorem II.6.1]{lw})
to construct a strongly deformation retractible neighbourhood $U'\subset
Y_\infty$. Without loss of generality, another identification point than
$P$ is not contained in $U'$.
Let $r'_t : U' \to U'\,,\;\;r'_1 = {\rm id}|_{U'}\,,\;\;
r'_0(U')=\{P\}$.
We want to show that
$$U:= U' \cup \bigcup^\infty_{i=1}f_i(U')$$ 
satisfies all assertions of this lemma.

First let us use the condition on sequences to show that $U$ is open:\neuzl
Assume that a sequence $(P_n)_{n\in\N}$ from the complement of $U$
converges to $P \in U$. 
Since by definition of convergent $Y$-spaces 
all $Y_i$ apart from $Y_\infty$ are isolated layers, 
$P \in U'\subset Y_{\infty}$. Since $U'$ is open in $Y_\infty$, 
a subsequence of $(P_n)_{n\in\N}$ must satisfy that 
$P_n \in Y_{i(n)}\setminus f_{i(n)}(U')$ with $i(n)$ strictly
monotonically increasing. Hence $f^{-1}_{i(n)}(P_n) \in Y_{\infty}\setminus U'
$, and by compactness  of a finite CW-complex this sequence 
will have an accumulation point in $Y_{\infty}\setminus U'$.
But by uniformity of $Y_i \to Y_\infty$, this would also have to be a limit
of a subsequence of the $P_n$-sequence, contradicting the assumption
that $P_n \to P$. \hl{Hence $U$ is open}. 

Now on 
$U$ we define our deformation retraction 
\begin{equation} 
r_t(x) := \left\{\vcenter{\halign{$#$\hfil&~#~&$#$\hfil\cr
       r'_t(x)& if &x \in U'\subset Y_\infty,\cr 
       f_i(r'_t(f_i^{-1}(x)))& if &                  x    \in 
f_i(U') \subset Y_i.\cr}}
\right.  
\label{eq:r-zweizl1}
\end{equation}
If $x $ is  the identification point $P$ (i.e., belongs to all $Y_i$),
the definition is not contradictory, because we constructed
$r$ as a strong deformation retraction. In order to show that
this definition is continuous with respect to the metric of $X$, we consider
two convergent sequences $(t_n)_{n\in \N} \in [0,1]$ and
$(x_n)_{n\in \N} \in X$                 with $t_n\to t$ and
$x_n\to x$. 
Since continuity of the restrictions of $r_t$ to each $Y_i$-layer
is clear, the interesting case quickly reduces, as in the preceding
paragraph, to the case when $x_n \in Y_{i(n)}$. 
By the uniformity with which $Y_i$ converge to $Y_\infty$, we can
conclude that dist$(x_n, f^{-1}_{i(n)}(x_n))\to0$.
Therefore also $f^{-1}_{i(n)}(x_n)\to x$.     Similarly
the investigation of convergence of the sequence $r_{t_n}(x_n)$
reduces to having to investigate $f^{-1}_{i(n)}(r_{t_n}(x_n))$, 
which by \hyperref[eq:r-zweizl1]{(\ref{eq:r-zweizl1})}
can be transformed to $
  r'_{t_n}(f^{-1}_{i(n)}(x_n))
$. 
However, having seen before that the arguments of $r'_{t_n}$
converge to $x$ and knowing that $r'$ is continuous in both
arguments, we get that the sequence under consideration converges
to $r_t(x)$, as desired to show continuity also in that case.
\end{pf}

\begin{expl}\label{dla-glmKonv}
For each $n\in\N$ inside $\R^2\subset \R^3$
we connect the points (0,0), (${1\over{2n}},n$),
(${1\over n},{1\over n}$), ($1,{1\over n}$) by a polygon. 
We smooth out these lines at the corners,       interpret them  
as   lines in $\R^3$ and perturb them a bit so that they are disjoint 
apart from their common start-point (0,0). With the latter modifications
we can try to interpret this union of lines as a convergent $Y$-space, where
$Y$ is the unit interval, the straight-line segment
$[(0,0),(1,0)]$ is the space $Y_\infty$, the other $Y_n$ are the just constructed
lines, $(0,0)$ is the only identification point, and the vertical
projection onto the $Y_\infty$-segment lying on the $x$-axis are
the inverses of the required homeomorphisms $f_n$.
All conditions for the construction of convergent $Y$-spaces
are satisfied, apart from that we have only pointwise
but not uniform convergence of the spaces $Y_n \to Y_{\infty}$.
However, we also note that this space is so embedded
into $\R^3$ that almost all $Y_n$ will be with two disjoint
neighbourhoods represented inside any ball-neighbourhood 
of the origin, so, similarly as for
\hyperref[sin-touch]{Example \ref*{sin-touch}}, although each
$Y_n$ has a deformation-retractible neighbourhood of (0,0),
there does not exist a deformation-retractible
neighbourhood for the entire space. Therefore it was essential
to require uniformity in the definition of a convergent $Y$-space,
since it was essentially needed in the proof of the preceding
\hyperref[dlaCorl3-3]{Lemma \ref*{dlaCorl3-3}} showing that
\hyperref[thm2]{Theorem \ref*{thm2}} can be applied
to convergent $Y$-spaces.
\end{expl}

\begin{expl}\label{CA-bis} 
First we embed $CA \setminus \buildrel\circ\over{l_\infty}$, 
i.e., our space $CA$ from the introduction with the interior of the accumulation
line removed, into
the $(x,y)$-plane of $\R^3$, so that the two intersection points
are $(0,0,0)$ and $(1,0,0)$. Then we let $CA''' :=$  
$$\{ (x,y,0) \mid (x,y) \in CA \setminus \buildrel\circ\over{l_\infty}
\}
\cup
\{ (x,0,z) \mid z \in [-1,+1],\;\;x \in \{0,1,
{\textstyle{1\over2},{1\over3},{1\over4},\ldots\}} 
\}$$
The CW-covering of this space is given by the spaces
$l_i$ from $CA$, and on the other hand by the spaces 
$\{(0,0)\} \times  [-1,+1]$ and  $\{(1/n,0)\}\times  [-1,+1]$ for all
$n\in\N$. 
The points $\{(0,0,0)\}$ and $\{(1,0,0)\} $ are the only intersection
complexes, and $CA'''$ satisfies all assumptions of  
\hyperref[thm2]{Theorem \ref*{thm2}}.

The remarkable property of $CA'''$ is its deformation
retractible neighbourhood of $(0,0,0)$:
$$ U = \{ (x,y,0) \mid (x,y) \in CA \setminus \buildrel\circ\over{l_\infty}\,,\;\;
x < \delta 
\} \cup (\{(0,0)\}\times (-\delta,+\delta))
$$  
Namely, according to the topological definition of the
boundary within $CA'''$, as $\partial U$
we do not only obtain the 
expected boundary points $(0,0,-\delta)$, $(0,0,+\delta)$,
and the points on the $l_i$ having $x$-coordinate $\delta$, 
but $\partial U$ contains also all points of type $(1/n,0,0)$ 
with $1/n<\delta$, and therefore also $(0,0,0) \in \partial U$. 
We have to take into account that our decision (see 
\hyperref[thm2]{condition \ref*{thm2}(v)}/\hyperref[Hyperbolas]{Example \ref*{Hyperbolas}}) 
to only require the deformation retractible neighbourhoods to 
be relatively open, does not rule out that this can happen.

Therefore some constructions (e.g.,  
\hyperref[eq:r-zweizl3]{(\ref{eq:r-zweizl3})}--\hyperref[eq:r-zweizl4]{(\ref{eq:r-zweizl4})}) will make
it necessary to look, apart from the deformation
retractible neighbourhood $U$, at another neighbourhood $U'$,
that is open in $X$ and satisfies $U'\cap X = U$.
In the case $X=CA'''$ the intersection with the metric neighbourhood,
$U((0,0,0),\delta) \cap CA'''$, that also \hl{will contain some
segments of some of the subspaces} $(1/n,0)\times [-1,+1]$,
could be used as $U'$. 
\end{expl}

\begin{lemma}\label{simplnu1}Let the assumptions of \hyperref[thm2]{Theorem \ref*{thm2}} be satisfied. Then any
(continuous) singular simplex $\sigma\colon\Delta^k\to X$ satisfies the following condition:

\begin{itemize}
\item[(a)] 
for each $p\in\Delta^k$, if we denote $I(p)\subset\N$ the maximal index set with $\sigma(p)\in Y_{I(p)}$, then there is an open neighbourhood $U(p)\subset\Delta^k$, such that
restriction of $\sigma$ gives a well-defined and continuous map $U(p)\to
U_{I(p)}$, where for $|I(p)| \ge 2$ the neighbourhood  $U_{I(p)}$
is \hl{$U^\cap_\nu $ assuming that the component of 
$Y_{I(p)} $ containing $\sigma(p)$ belongs to the intersection
complex $Y^\cap_\nu$}, 
and for 
$I(p) =\{i\}$ it is an open neighbourhood of $Y_i$.   
\item[(b)] In particular, in the case of a 1-simplex $\sigma:
[0,1]\to X$, there exists a finite sequence $0=p_0,p_1,\ldots,p_n=1$, such that
each set 
$\sigma([p_i,p_{i+1}])$ is either contained in some space $Y_j$
or in some neighbourhood $U^\cap_\nu$.
\end{itemize}
\end{lemma}

\begin{pf}
Let us first consider only a 1-simplex $\sigma:\Delta^1\to X$. 

We have a countable, closed covering \hl{$              
\left\{\sigma^{-1}(Y_i)\mid i\in \N\right\}\cup
 \left\{\sigma^{-1}(Y^\cap_\nu) \mid \nu\in\Ups\right\}
$ of $\Delta^1$}. 
The sets in this covering may be disconnected, so let us look at their
path components. We obtain the              covering ${\mathcal U}_1$
by adding to this covering by path-components the following sets: 
If $[P_1,P_2] $ and $[P_3,P_4]$ 
are path-components of  
$\sigma^{-1}(Y^\cap_\nu)$ for the same intersection component $Y^\cap_\nu \;\;
(P_1\le P_2 < P_3 \le P_4)$,    
and if the entire segment $[P_2,P_3] \subset \sigma^{-1}(U^\cap_\nu)$,
then 
$[P_1,P_4]$ will also be a covering set in ${\mathcal U}_1$.
Note that,              by making use of the deformation retraction
according to \hyperref[thm2]{Theorem \ref*{thm2}(iii)},
the path could be homotoped on $[P_2,P_3]$ into $Y^\cap_\nu$.

Define an equivalence relation on $\left[0,1\right]$ by declaring $P\sim Q$ if and only if there exists a finite sequence $P=p_0,p_1,\ldots,p_n=Q$ such that
$p_i$ and $p_{i+1}$ belong to  the same set $A_i\in{\mathcal U}_1$, for $i=0,\ldots,n-1$. Let ${\mathcal U}_2$ be the covering by equivalence classes of $\sim$. 
This is a disjoint covering by path-connected subsets of
$\left[0,1\right]$, which by 
the classification of connected subsets of the real line 
must be either points, closed, open or half-open intervals. 
However, since the covering sets of $\mathcal{U}_1$ are closed,
infinitely many of them must have been united to obtain an
open or half-open interval as a $\mathcal{U}_2$-covering set.
Looking at the above definition of ``$P\sim Q$\/", this means, 
although the equivalence of any two points in the (half-)open 
interval must be based on a finite list of points $
p_0,p_1,\ldots,p_n$, that this list of points could be prolonged 
at least to one side to an infinite list. Each of these points belongs
to an intersection complex.
If such a list is composed economically, then $\sigma|_{[p_i,p_{i+1}]}$ 
must have
\begin{itemize}
\item            either run between two \hl{different intersection
components $Y^\cap_\nu$}, or
\item            left one           \hl{intersection component} 
$Y^\cap_\nu$ and its deformation retractible neighbourhood $U^\cap_\nu$
(but remained within the spaces $Y_i$ \hl{intersecting there}) 
and then returned to $Y^\cap_\nu$, or
\item            gone from some part of $Y^\cap_\nu\cap Y_i$
into some part of $Y^\cap_\nu$ that is not covered by $Y_i$.
\end{itemize}
However, no continuous path, that is parametrized on a compact interval, 
can do infinitely many of such moves within the same finite number of 
spaces $Y_i$ and $Y^\cap_\nu$. 
Therefore, the only chance how one could get this infinite
list of points $p_i$ is, that they come from infinitely
many different intersection components. But then the end-point
not belonging to the (half-)open interval would be an accumulation
point of such a sequence, and this is ruled
out by assumption (v) of \hyperref[thm2]{Theorem \ref*{thm2}}. 

So we have a covering by disjoint points and closed intervals. 
If $\mathcal{U}_2$ is not the trivial covering 
by only one set, it must be infinite.

If $\mathcal{U}_2$ is an infinite covering, 
consider the quotient space $\left[0,1\right]/\sim$ with the countable covering by the closed sets $\sigma^{-1}(Y_i)/\sim$. Since equivalence classes of $\sim$ are points or closed intervals, the quotient space 
$\left[0,1\right]/\sim$ is again 
homeomorphic to $\left[0,1\right]$. On the other hand, the sets $\sigma^{-1}(Y_i)/\sim$ are totally disconnected 
and thus have empty interior. So we obtain a covering of $\left[0,1\right]$ by countably many closed sets of empty interior. But this contradicts Baire's category theorem.

In the remaining case of $[0,1]\in \mathcal{U}_2$ being the only covering 
set, we directly get Statement~(b).

To show (b)$\Rightarrow$(a)
for 1-simplices, it does not suffice to observe that, since
we have $0=P\sim Q=1$, any $p\in[0,1]$ has some neighbourhood which is either
free of all subdivision points $p_0,p_1,\ldots,p_n$ or $p$ is the only one
of them contained there, since the following case is possible: If $p\in (P_2,P_3)$,
where we have $P_1\le P_2 < P_3 \le P_4$, and $[P_1,P_4]$ has
become, as above described, a covering set of $\mathcal{U}_1$
of the second kind, then the claim which needs to be proven 
for such a parameter $p$ (provided that $\sigma(p)$ is outside
of $Y^\cap_\nu$) is not that $\sigma(U(p))$ is contained 
in the neighbourhood $U^\cap_\nu$ (as we know $\sigma$ to fulfil
on the entire interval $[P_1,P_4]$), but we need to show that
the $\sigma$-values are contained in some $Y_i$ for a suitable neighbourhood
of $p$. However, already closedness of $Y^\cap_\nu$ and continuity of
$\sigma$ imply that some neighbourhood $U(p)$ exists such that
$\sigma(U(p))\cap Y^\cap_\nu=\emptyset$. Then we can apply to 
$\sigma|_{U(p)}$ the same argument as in the main part of the proof,
and by assumption and 
\hyperref[thm2]{condition \ref*{thm2}(iv)}
knowing that $\sigma$ can on this             
segment not go into any of the intersection complexes, it results that 
$\sigma|_{U(p)}$ must have remained in the $Y_i$ that contains $\sigma(p)$.

We now show (a) for simplices of arbitrary dimension
by contradiction: Assume that in the neighbourhood $p\in\Delta^k$
there exists a sequence $q_i\to p$ with $\sigma(q_i)$ not being contained
in the set as asserted in the statement of the lemma. Connect
the $q_i$-points by line-segments and observe that they plug together
to a continuous curve that can be extended to the endpoint $p$. 
The restriction to that curve would be a 1-simplex 
contradicting (b).                                           
\end{pf}                     
\begin{remrk}\label{simplnu1po}
The purpose of this remark is to point out, that the conclusions of \hyperref[simplnu1]{Lemma \ref*{simplnu1}}
can also be applied in the course of the proof of \hyperref[negcurvinj]{Corollary \ref*{negcurvinj}(b)}: 

In the proof of \hyperref[simplnu1]{Lemma \ref*{simplnu1}} we did not use \hyperref[thm2]{Condition \ref*{thm2}(iii)}, but just the 
\hyperref[thm2]{Conditions \ref*{thm2}(v)}, \hyperref[thm2]{ \ref*{thm2}(iv)} and indirectly, when referring to $U^\cap_\nu$, also \hyperref[thm2]{Condition \ref*{thm2}(iii)}. 

Now \hyperref[thm2]{Condition \ref*{thm2}(iv)} and \hyperref[thm2]{\ref*{thm2}(iv)} are automatically satisfied in the
case to be considered by \hyperref[negcurvinj]{Corollary \ref*{negcurvinj}(b)}, since we are having only 
one intersection complex. Concerning \hyperref[thm2]{Condition \ref*{thm2}(iii)}, we did in the
above proof not use the existence of the deformation retraction, but just 
the existence of the neighbourhood. And since also the assumptions of 
\hyperref[negcurvinj]{Corollary \ref*{negcurvinj}(b)} contain the existence of some neighbourhood $U$ of $l$,
this $U$ can be used instead of $U^\cap_\nu$, when applying the arguments
of \hyperref[simplnu1]{Lemma \ref*{simplnu1}} in case of \hyperref[negcurvinj]{Corollary \ref*{negcurvinj}(b)}.
\end{remrk}

\subsection{Reductions for the proof of 
\hyperref[thm2]{Theorem \ref*{thm2}}
}\label{reduction}
\neuzl
In the setting of \hyperref[thm2]{Theorem \ref*{thm2}}, for $i\in\N$ 
we let $Y_{i}^\cup$ be the union of \hl{intersection  components of $Y_J
\,,\,(J\subset \N)$} having nonempty intersection with $Y_{J\cup i}$.
We denote \hl{$$Z_i=Y_i \cup Y^\cup_i
                $$} 
i.e., we add to $Y_i$ all those cells that belong to one \hl{intersection
component that}  shares at least one vertex with $ Y_i$. \hl{%
By \hyperref[stattviLemma]{Lemma \ref*{stattviLemma}} the} sets
$Z_i$ are also finite CW-complexes, and the  subdivision         
of $X$ into the $Z_i$ satisfies also all other assumptions
as required for the $Y_i$ in         \hyperref[thm2]{Theorem \ref*{thm2}}.
However, 
every component of every intersection $Z_I=\bigcap_{i\in I}Z_i$, $ \vert I\vert\ge 2$, \hl{ already coincides with the component of
$ Z^\cap := \{z\in Z \mid  \exists{i\neq j\in\N} 
\hbox{ with }z \in Z_i\cap Z_j\}$
(in contrast to the situation as described in            
\hyperref[eq:AnzTeilm]{(\ref{eq:AnzTeilm})})   
and is therefore contractible by assumption (ii) of
\hyperref[thm2]{Theorem \ref*{thm2}}}. 
 So we can reduce the proof of \hyperref[thm2]{Theorem \ref*{thm2}} to the case where all components of all $Y_I\,,\;\vert I\vert\ge 2$, are already 
 \hl{intersection components as defined in   
\hyperref[thm2]{Theorem \ref*{thm2}} satisfying the conditions stated 
there}. We will {\it henceforth make this assumption without losing generality}.

Also, we {\it will henceforth assume}\/ that index sets 
denoted by $I \subset \N$ and $J \subset \N$ always have to contain
at least two indices, if not explicitly stated otherwise.

Despite having reduced to \hl{intersection components $Y^\cap_\nu$,
we still need to distinguish between 
the entire intersection
$ Y_I = \bigcap_{i\in I}Y_i$, that might comprise distinct 
complexes $Y^\cap_\nu$ and considering just one of them. 
Depending on the context, we might denote one of our intersection
components as 
\begin{equation}
Y^\cap_\nu~,~~    Y^0_I ~\hbox{ or as }~ Y^\nu_I,   
\label{eq:int-compsY}      
\end{equation} 
depending on whether we are only giving the index $\nu\in\Ups\not\ni0$
of the enumeration of intersection components, only the subset
$I\subset\N$  of the indices of spaces $Y_I$ intersecting there, 
or both kinds of information, respectively. However by our prior reduction
each symbol $Y^0_I$ and $Y^\nu_I$ will just denote one
of the components $Y^\cap_\nu$ (or $\emptyset$ in case of mismatching
$I$ and $\nu$).}

Recall that in \hyperref[thm2]{Theorem \ref*{thm2}(iii)}
we demanded the existence of \hl{deformation retractible
neighbourhoods, that now with the analogous
philosophy will either be denoted by
$U^\cap_\nu$, $U^0_I$ or $U^\nu_I$ (with the understanding
that an empty $Y^\nu_I$ has only an empty neighbourhood). 
Based on this we now define
\begin{equation}
U_i:=Y_i\cup\bigcup_{J:J\owns i \atop  \nu\in\Ups }U_{J}^\nu\,,~~~ U_I := \bigcup_{\nu\in\Ups} U^\nu_I.
\label{eq:int-compsZ}      
\end{equation} 
}

This means we are adding to $Y_i$  (or taking, resp.) not only the \hl{ 
entire intersection components, }
but also their deformation retractible neighbourhoods. 
Recall that by \hyperref[thm2]{Condition \ref*{thm2}(iii)}
and \hyperref[Hyperbolas]{Example \ref*{Hyperbolas}}
they need 
not to be open. 

Observe that as a $T_4$-space satisfying the first axiom of countability,
$X$ is by 
\cite[1.6.14--1.6.15]{Egkl}
a sequential space. Hence we can and will make use
of the well-known characterizations of 
continuity from metric topology and of open and closed sets via sequences.
\neuzl\hbox{~}\neuzl
{\it In the remainder of this subsection we will construct 
 the maps $r_t$ and $r_{i,t}$}.%
\neuzl\hbox{~}\par
First recall the following construction that can be done
for any $T_4$-space.\neuzl
Let $A$ be a closed subset, $U$ an (open) neighbourhood of $A$, that
deformation retracts to $A$ via a retraction $r_t$ with 
$r_1 = $id$|_U$ and $r_0(U)=A$. We separate the closed sets
$X\setminus U$ and $A$ by an Urysohn function $f$ with $
f|_A = 0$ and $f|_{X\setminus U}=1$. We then define 
$f^{-1}([0,1/2)) $ as a smaller neighbourhood $V$ of $A$
satisfying $\overline{V}\subset U$. Based on this we hereby
construct the following modification of $r$:

\begin{equation}    
\widehat{r_t}(x) := 
\left\{\vcenter{\halign{$#$\hfil&~#~&$#$\hfil\cr
r_t(x)&if&0\le f(x)\le{1\over2},\cr
r_{(2-2t)\cdot f(x) + 2t-1}(x)&if&{1\over2}\le f(x)\le1.\cr
}}\right.
\label{eq:r-zweizl2}      
\end{equation}  
Observe that $\widehat{r}$, \hl{other than} $r$, only deforms  
a smaller neighbourhood $V$ into $A$, and by its nature of definition
it cannot even be ruled out that it may take values which are
only in $U \setminus V$, but \hl{other than} the original $r$
it is the identity for all $t$ on $\partial U$. That way it can be 
continuously extended
to the outside of $U$ by the identity map. 

We now apply an analogous mechanism to our deformation-retractible
neighbourhoods $U^\cap_\nu$, that we required to exist in assumption (iii)
of \hyperref[thm2]{Theorem \ref*{thm2}}.
That way we define on $X$ globally a system of functions $r=(r_t)_{t\in[0,1]}$
via
\begin{equation}    
 r_t(x) := \left\{\vcenter{\halign{#&#\hfil\cr
id &$\forall t\in[0,1]$, if $x$ is outside of all $U^\cap_\nu
~~~(\forall \nu\in\Ups)$,\cr
resu&{\vtop{\noindent\hsize=0.65\hsize lts from the adaptation of the deformation
retraction, that is given in each of the $U^\cap_\nu$, as 
in \hyperref[eq:r-zweizl2]{(\ref{eq:r-zweizl2})}.
}}\cr}}\right.              
\label{eq:r-zweizl3}      
\end{equation}  
Since the $U^\cap_\nu$ 
need not be open (but just relatively open, cf.\     
\hyperref[CA-bis]{Example \ref*{CA-bis}}), 
for each of them 
we find a neighbourhood $U$ that is
open in $X$, so that $U \cap \bigcup_{I:I\owns i} Y_i= U^\cap_\nu$.
With the help of this $U$ we define 
the necessary 
Urysohn function $f$. 
On the other hand, 
the second line of 
\hyperref[eq:r-zweizl2]{(\ref{eq:r-zweizl2})}
will only be applied 
in  $U\cap\bigcup_{I:I\owns i} Y_i$, at the other points of $U$ 
the definition of $r_t$ does not get changed. 

However, where, as in our
\hyperref[Hyperbolas]{Example \ref*{Hyperbolas}}, we have $Y_j$-spaces
accumulating to an intersection complex without intersecting there,
we cannot expect that the globally defined map $r_t$, will be actually
continuous. 
On the other hand for each subspace $U_i$ our map
\begin{equation}  
r_i := r |_{U_i}
\label{eq:r-zweizl4}      
\end{equation}  
will be continuous. 
By construction we have $r_{i,t} : U_i \to U_i\,,\;\;r_{i,1} = {\rm id}|_{U_i}
\,,\;\; r_{i,0}(U_i) = Y_i$, and therefore we will sometimes refer to these
maps as ``deformation retractions", although they are no retractions
in the strict sense: Since the adaptation 
as performed in \hyperref[eq:r-zweizl2]{(\ref{eq:r-zweizl2})}
needs for continuity reasons also to be
applied to $Y_i\cap U^\cap_\nu$, we cannot expect to have
$r_{i,0}|_{Y_i}={\rm id}$, but only $\hl{r_{i,1}}|_{Y_i}\simeq{\rm id}$.
However, this does not mess up any of the proof-constructions to come. 

We {\it hereby will assume},\/ without 
loss of generality, that, similarly as in the first paragraph of this
construction, a smaller neighbourhood has been defined,                   and
$U^\cap_\nu$, as assumed to exist by \hyperref[thm2]{Theorem \ref*{thm2}(iii)},
has been reduced to $U^\cap_\nu\cap f^{-1}([0,1/2))$. 
The restriction of the globally defined map $r_t$ to the accordingly
adapted neighbourhood $U^\cap_\nu$ is still  a  deformation retraction
onto $Y^\cap_\nu$; the fact that it might for $t \in (0,1)$ take some
values outside the correspondingly adapted neighbourhood $U^\cap_\nu$
does not mess up
any of the proof-constructions that will come.

\subsection{Barycentric subdivision}

\begin{lemma}\label{sd}Under the assumptions of \hyperref[thm2]{Theorem \ref*{thm2}} every measure chain $\mu$ has an iterated barycentric subdivision $sd^m(\mu)$ (for some large enough $m$) which can be decomposed as $sd^m(\mu)=\sum_{i=1}^\infty\mu_i$ with the \hl{compact} determination set of $\mu_i$ contained in $U_i$ for every $i\in\N$.\end{lemma} 

\begin{pf} Recall that $X=\bigcup^\infty_{i=1} Y_i$
and that by \hyperref[reduction]{Subsection \ref*{reduction}}
$I\subset\N\,,\;\;|I|\ge2$. Let \hl{% 
$$M_i = \{\sigma \in map(\Delta^k,Y_i)\mid \sigma(\Delta_k) \cap
Y_I = \emptyset ~~~ \forall{I:I\owns i}  \}$$
and 
$$M_I:= \{\sigma \in map(\Delta^k,U_I) \mid \exists{x \in \Delta^k}
\hbox{ with }\sigma(x) \in Y_I\}.$$}
\hl{Observe that} \hyperref[simplnu1]{Lemma \ref*{simplnu1}}
in principle says that every simplex $\sigma \in map(\Delta^k,X)$
locally behaves as if it 
has to belong to one of the sets $M_i$ or $M_I$.   

Let us, however, assume that this does \hl{not} suffice, so that an $m\in\N$
would exist so that $D_m$, i.e., 
the determination set of $sd^m(\mu)$, only would contain
simplices in $\bigcup_i M_i\cup\bigcup_I M_I$. 
Such singular simplices in $map(\Delta^k,X) \setminus
(\bigcup_i M_i\cup\bigcup_I M_I)$ have to contain points in at least
one of the intersection \hl{components} $Y^0_I$, but also outside
the corresponding neighbourhood $U^0_I$. Accordingly, by our assumption, we are
able to choose for each $m\in\N$ some simplex $\sigma_m\in D_m$,
together with $x_m,y_m\in\Delta^k$ and $\sigma_m(x_m)\in Y_{I(m)}^0$,
but $\sigma_m(y_m)\notin  U_{I(m)}^0$. Actually, we can even assume 
$\sigma_m(y_m)$ to belong to $(\bigcup_{i\in I(m)} Y_i)\setminus U_{I(m)}^0$,
because by the conditions of \hyperref[simplnu1]{Lemma \ref*{simplnu1}(b)} 
a simplex $\sigma_m$ with $\sigma_m(x_m)\in Y^0_{I(m)}$ \hl{can only take
a value in $Y_j$ with $j\notin I(m)$ if a path} (that is defined
by the restriction of $\sigma_m$) has taken a value in another
intersection complex than $Y^0_{I(m)}$ before reaching $Y_j$. 
Since intersection complexes and their neighbourhoods are disjoint
(see \hyperref[thm2]{Condition \ref*{thm2}(iv)}), when a path
exits $U^0_{I(m)}$, it will first take values in 
$\bigcup_{i\in I(m)} Y_i$, before it might go into another $Y_j$.

Using standard techniques of sequences we can reduce the situation
to one of two cases: Either 
\hl{the corresponding intersection component}
$Y^0_{I(m)}$ is for all $m$ the same complex, or for each
$m$ a different complex. However, our condition
\hyperref[thm2]{Theorem \ref*{thm2}(v)} rules out that the    
compact image of our map $F\colon\Delta^k\times D_1\to X\,,\;\;
(x,\sigma) \mapsto\sigma(x)$
can intersect with infinitely many 
intersection complexes, and hence we can reduce to the case
of only one intersection complex $Y^0_I$ and simplify the
notation $I(m)$ to $I$ accordingly. 

Now, each simplex $\sigma_m$ is the restriction of some simplex 
$\sigma'_m\in D_1$ to a subsimplex of $\Delta^k$ from its 
$m$-th barycentric subdivision, and via this identification
of the domain of definition of $\sigma_m$ with a subdomain
for $\sigma'_m$, we define points $x'_m$ and $y'_m\in \Delta^k$,
so that these points as arguments for $\sigma'_m$ correspond 
to $x_m$ and $y_m$ as arguments for $\sigma_m$.

Having defined this sequence of singular simplices $\sigma'_m$,
we can, by compactness of $\Delta^k$ and the image of the above-defined
function $F$, without loss of generality, take
a subsequence so that the following sequences converge: 
$ x'_m \to x \in \Delta^k,\;\; y'_m \to y \in \Delta^k,\;\; 
\sigma'_m(x'_m) \to P_0\in Y^0_I,\;\;\sigma'_m(y'_m) \to Q_0\notin U^0_I$.
Compactness of $D_1$ allows then to choose a singular simplex
$\sigma_\infty$ as an accumulation point of this subsequence.
We can then consider $P_\infty := \sigma_\infty(x) = 
\lim(\sigma_\infty(x'_m))$ and $Q_\infty := \sigma_\infty(y) = 
\lim(\sigma_\infty(y'_m))$. If $P_0   \neq                 P_\infty$,
they could be separated by disjoint neighbourhoods $U(P_0)$ 
and $U(P_\infty)$, and in addition a compact neighbourhood $U(x)$
with nonempty interior around $x$ could be chosen so as to
contain almost all of the           $x'_m$ from the selected subsequence. 
That way clearly $map(U(x),U(P_0))$ and $map(U(x),U(P_\infty))$ would be
disjoint open sets in $map(\Delta^k,X)$, the latter of which would contain
(at least, after some adaptation of these neighbourhoods) our simplex $
\sigma_\infty$, but conversely could not contain almost all
of the simplices $\sigma'_m$. But this would then rule out that
$\sigma_\infty$ could be an accumulation point of the
priorly chosen subsequence of $\sigma'_m$. 

From that contradiction, $P_\infty = P_0$ and 
analogously     $Q_\infty = Q_0$, so finally 
$P_\infty= P_0\neq Q_0 = Q_\infty$.

On the other hand, $x'_m$ and $y'_m$ are contained in the same
simplex of the $m$-th iterated subdivision of $\Delta^k$, and hence
for their distance we get dist$(x'_m,y'_m)\to 0$. Consequently,
$x=y$ and also $P_\infty = \sigma_\infty(x) = \sigma_\infty(y) =
Q_\infty$. 

Since the final formulae of the two preceding paragraphs
contradict each other, our assumption on $m$ at the
beginning of this proof must have been wrong. 
Therefore we conclude that there exists some $m \in \N$
so that all simplices $\sigma \in D_m$ are contained in one
of the sets $M_i$ or $M_I$.     

Now $M_i$ and $M_I$ are measurable sets, so if we 
define 
$$\mu'_i := \mu|_{M_i} \hbox{~~~and~~~}\mu_I := \mu|_{M_I},$$
we obtain well-defined measures. Since, as noted before, our 
measure chain can only go into finitely many intersection complexes, 
let us accordingly restrict to a finite list $\{Y_{I_1},\ldots,
Y_{I_n}\}$. The measures $\mu'_i$ have by closedness of $Y_i$
compact determination sets $D_m \cap map(\Delta^k,Y_i)$. This
allows us to consider these measures as elements of the Milnor-Thurston
chain groups of the spaces $Y_i$. In order to obtain 
determination sets for $\mu_I$ we take 
$D_m\cap M_I$; the $U_I$ are not closed, but simplices having at least
one point in $Y_I$ cannot converge \hl{towards} any other simplices
than those having at least one point in $Y_I$. For each of the
relevant index sets $I_1,\ldots,I_n$ we pick an arbitrary element
$i_1\in I_1,\ldots , i_n\in I_n$ and let 
$$\mu_i := \mu'_i + \sum^n_{s=1} \delta_{i,i_s}\cdot \mu_{I_s},$$
where $\delta_{i,i_s}$ denotes the Kronecker symbol. 
This means that we add each of the finitely many measures 
$\mu_{I_s} $ to one (but only one) of the $\mu_i$-measures,
where $Y_i$ passes through the according intersection complexes
$Y_{I_s}$. That way we finally obtain that
$$sd^m(\mu) = \sum^\infty_{i=1} \mu_i$$ with $\mu_i$ being compactly
determined on $U_i$.

\end{pf}

\subsection{Proof of \hyperref[thm2]{Theorem \ref*{thm2}}}\label{proo}

\begin{pf}
We want to prove that
$\iota_*\colon H_k(X;\R)\to {\mathcal{H}}_k(X)$ is injective for $k\ge 2$.

Assume that $\iota_*$ is not injective, which means that there is some $\alpha\in H_k(X,\R)$, represented by a cycle $z\in C_k(X,\R)$ such that $z=\partial\mu$ for some measure chain $\mu\in {\mathcal C}_{k+1}(X)$. 

By the discussion in \hyperref[reduction]{Subsection \ref*{reduction}} we can  without losing generality assume that all components $Y^0_I$ of all $Y_I
$ are already 
some \hl{intersection components $ Y^\cap_\nu$} 
(as defined in the statement of  \hyperref[thm2]{Theorem \ref*{thm2}}).

Recall the notation $U_I^0, U_I  $ and $U_i$ from \hl{% 
\hyperref[eq:int-compsY]{(\ref{eq:int-compsY})}%
--\hyperref[eq:int-compsZ]{(\ref{eq:int-compsZ})}}.
By their definition 
        and adaptation 
at the end of that subsection,
each
$U_i$ is a (not necessarily open) neighbourhood of $Y_i$. Even though the $U_i$ may be not open, \hyperref[sd]{Lemma \ref*{sd}} gives that for sufficiently large $m$ the $m$-th barycentric subdivision $sd^m\mu$ is determined on simplices with image each contained in one of the $U_i$. That is, we have a decomposition $sd^m(\mu)=\sum_{i=1}^\infty \mu_i$ with $\mu_i$ determined on $U_i$. The classical construction of a chain homotopy between $sd^m$ and the identity (on singular chains) implies that $\partial sd^m\mu=sd^m z$ is a cycle in the homology class $\alpha$. Similarly, the iterated subdivision of $z$ decomposes as $sd^m(z)=\sum_{i=1}^\infty z_i$ with $z_i$ a linear combination of simplices with image in $U_i$. (The sum is actually finite because $sd^m(z)$ is a finite linear combination of simplices.)

Note that $\partial\mu=z$ does not imply that individually $\partial \mu_i=z_i$ holds, 
but some simplices in $\partial\mu_i-z_i$ might 
cancel against simplices in other $\partial \mu_j-z_j$ when
$\partial\mu-z=0$ is computed. 
Such a cancellation can however only occur inside the intersection $U_i\cap U_j$. So 
on $ U_i \setminus \bigcup_{I:I\owns i} U_I$ we get that
$\partial\mu_i-z_i$ must vanish.    

Recall that in 
\hyperref[eq:r-zweizl3]{(\ref{eq:r-zweizl3})}--\hyperref[eq:r-zweizl4]{(\ref{eq:r-zweizl4})}                                         
we constructed
a homotopy $r_i$ satisfying up to homotopy all properties
of a deformation retraction. 
Thus $r_{i,*}(\mu_i)$ is a measure chain determined on $Y_i$. Moreover, 
$\bigcup_{I : I\owns i}U_I$
is sent to $\bigcup_{I : I \owns i }Y_I$ and therefore $r_{i,*}(\partial\mu_i)-r_{i,*}(z_i)$ is determined on $map(\Delta^k, \bigcup_{ I:I\owns  i     }Y_I)$.

Similarly, $\partial z_i$ consists of singular cycles in $\bigcup_{I : I\owns i}U_I$, hence $\partial r_{i,*}(z_i)$ consists of singular cycles in the contractible sets 
$\bigcup_{I : I\owns i}Y_I$.  
Now $\partial z_i$ is a cycle, and by disjointness of the intersection
complexes $Y^0_I$ whatever falls from $\partial z_i$ into one of these
complexes $Y^0_I$ must also be a cycle. 
Since each of the finitely many complexes $Y^0_I$, where our singular chain 
$z_i$ can go in, is contractible,
the parts of the boundary of $\partial z_i$ that are in this
complex must be the boundary of some singular chain, that together 
give a singular chain 
$v_i \in {C}_k(\bigcup_{I:I\owns i}Y_I)$. It is best 
to construct such $v_i$ from the coning construction
of \hyperref[1new]{Lemma \ref*{1new}}, in order to make sure for later cancellation
purposes that they are performed by the same method. Observe that this step requires $k\neq1$
(i.e., dim$(\partial z_i)\neq 0$), as  noted
in \hyperref[1new]{Lemma \ref*{1new}}, because otherwise 
this boundary chain $v_i$ is not known to exist by this general 
argument. However, for 
$k\geq2$ we can obtain $v_i$
such that $\partial v_i=\partial r_{i,*}(z_i)$, i.e., $r_{i,*}(z_i)-v_i$ is a 
         singular cycle  in the $Y_I$ with $i\in I$.

From $\partial v_i=\partial r_{i,*}(z_i)$ it follows that $r_{i,*}(\partial\mu_i)-r_{i,*}(z_i)+v_i$ is a measure cycle in $Y_i$. It belongs to the disjoint union of the contractible sets $Y^0_I$ with 
$i \in I$,
which have vanishing measure homology in degrees $k\ge 1$. Therefore it must be a measure boundary, i.e., there 
is a signed measure $\nu_i$ with $$\partial\nu_i=r_{i,*}(\partial\mu_i)-r_{i,*}(z_i)+v_i,$$
where $\nu_i$ could be also constructed with the coning construction
from \hyperref[1new]{Lemma \ref*{1new}}.

Thus, on $Y_i$ we have $$\partial (r_{i,*}(\mu_i)-\nu_i)=r_{i,*}(z_i)-v_i.$$
On the right hand side, we have a singular cycle, which according to the left hand side is the boundary of a measure chain.

But $Y_i$ is a CW-complex, thus the canonical homomorphism
$$\iota_*\colon H_*(Y_i;\R)\to \mathcal{H}_*(Y_i)$$
is an isomorphism according to \cite{zas} and \cite{han}. In particular it is injective, that is, a singular cycle which is the boundary of a measure chain must be 
the boundary of a singular chain $s_i$.
 
Recall that $sd^m(z)=\sum_{i=1}^\infty z_i$ was a finite sum, hence there are only finitely many $v_i\not=0$ and we can set all but finitely many $s_i$ equal zero. That means, $\sum_{i=1}^\infty s_i$ is a finite sum, thus defines a singular chain.

In order to understand $\sum_{i=1}^\infty v_i$, recall that \hl{$z=\sum_{i=1}^\infty z_i$ and that the various
homotopies $r_i$ are just restrictions from one globally defined
map $r$, that was only non-identically (but consistently)
defined on each neighbourhood $U^0_I$ of an intersection 
\hl{component}, hence
$$\sum_{i=1}^\infty \partial v_i=\sum_{i=1}^\infty\partial(r_{i,*}(z_i))
=\partial(r_{*}(z))
=0.$$}
This means, that if the filling chains $v_i$ have sufficiently
consistently been constructed in each intersection complex 
$Y^0_I$ (e.g., by the coning construction from 
\hyperref[1new]{Lemma \ref*{1new}}), $\sum_{i=1}^\infty v_i$ will automatically
cancel to zero. 

We obtain $$\partial(\sum_{i=1}^\infty s_i)=\sum_{i=1}^\infty(r_{i,*}(z_i)-v_i)=\sum_{i=1}^\infty r_{i,*}(z_i),$$ that is, $\sum_{i=1}^\infty r_{i,*}(z_i)$ is the boundary of a singular chain. Each $r_i$ is homotopic to the identity, thus $r_{i,*}(z_i)$ is homologous to $z_i$ and $\sum_{i=1}^\infty r_{i,*}(z_i)$ is homologous to $z$. This yields that $z$ is the boundary of a singular chain, what we wanted to prove.

\end{pf}

\subsection{Example: convergent Y-spaces}\label{hypex}

Recall from \hyperref[negcurv]{Section \ref*{negcurv}} the \hyperref[convy]{Definition \ref*{convy}} of convergent $Y$-spaces for a 
metric space $Y$.

For $Y=\left[0,1\right]$ with exactly the two points $0$ and $1$ along which the identification is done,
this yields the convergent arc space drawn in the introduction. 

Such a convergent $Y$-space space will usually 
not be a CW-complex even if $Y$ is: although $X$ inherits a cell decomposition from that of $Y$, the accumulation property $\lim_{n\to\infty} {\rm dist}(f_n(y),y)=0$ implies that $X$ does not have the weak topology with respect to that cell decomposition.

The assumption on nonvanishing of the Gromov norm, that we required in \hyperref[thm1]{Theorem \ref*{thm1}} is in general not satisfied for convergent $Y$-spaces. For example, it does not hold for the convergent arc space, or when the Gromov norm on $Y$ is not an actual norm, e.g., when $Y$ is simply connected. 

We mention that for $Y$ homeomorphic to a simplicial complex 
or to a CW-complex (e.g.\ a smooth manifold, see \cite[Chapter 2, Theorem 10.6]{munkr})
and $X$ a convergent $Y$-space,  
for all $k$ the homotopy classes (rel.\ vertices) of $k$-simplices with vertices in $x_0$ are Borel sets in $map(\Delta^k,X)$. Thus the assumptions of \hyperref[lem2]{Lemma \ref*{lem2}} are satisfied in this case. We will not include the rather lengthy argument because it is not needed for the proof of \hyperref[ls1]{Corollary \ref*{ls1}}.

\begin{cor}\label{ls1}If the            metric space $Y$ is homeomorphic to a finite CW-complex and if $X$ is a convergent $Y$-space,
then $$\iota_*\colon H_k(X;\R)\to {\mathcal{H}}_k(X)$$
is injective in degrees $k\ge2$.\end{cor}

Concerning the validity of this corollary for $k=0$ and $k=1$, the same
remark applies as given immediately after the statement 
of \hyperref[thm2]{Theorem \ref*{thm2}}. 

\begin{pf}
That a convergent $Y$-space satisfies the conditions of 
\hyperref[thm2]{Theorem \ref*{thm2}} is obvious apart from condition
(iii), i.e., the existence of deformation retractible relatively
open neighbourhoods. But this has been shown in 
\hyperref[dlaCorl3-3]{Lemma \ref*{dlaCorl3-3}}.
\end{pf}

\section{ 
Finishing the proof of \hyperref[negcurvinj]{Corollary \ref*{negcurvinj}(b)}.}   
%--------------------------------------------   
\label{newSectSam}

Throughout this section we treat treat all notations from \hyperref[negcurvinj]{Corollary \ref*{negcurvinj}(b)} 
and \hyperref[newEndRem13a]{Remark \ref*{newEndRem13a}}
as set, and all assumptions from \hyperref[negcurvinj]{Corollary \ref*{negcurvinj}(b)} as given.
   
For reaching the goal of this section we will need 
to make use of our assumptions that each surface has
a hyperbolic structure, and start by fixing the convention, how
to measure angles. By ``differentiable curve" we mean a curve
with a differentiable parametrization that has a nowhere-vanishing
derivative.

\begin{defi}\label{newSectDef0}
Let $c,d$ be two differentiable curves in a surface, where there
is an oriented Riemannian metric structure defined. When both curves
meet at the point $P$, we can measure the angle between them. We denote
          the angle by the following operator  
$$
      \angle(P , c , d)
$$      
and want to understand the angle by the orientation of the underlying
surface as being defined uniquely in 
$\R\over2\pi\ZZ$. Reversing the order 
of the second and third argument of the angle-operator will  
reverse the sign of the angle-operator.    
If at least one of the curves is not differentiable at the 
point $P$, we might have to distinguish between the angle, at which
the curves approach $P$, and at which they depart. In this case
the above formula refers to the angle, at which the curves depart,
while in order to measure the angle at which they approach, we 
write:
$$ 
    \angle(P, c^{-1}  , d^{-1}  ), 
$$    
i.e., put the reversed curves  
into our operator.
\end{defi}

In order to apply a similar straightening mechanism as in the proof
of \hyperref[negcurvinj]{Corollary \ref*{negcurvinj}(a)}, we first have to give a definition, when 
a curve is straight in our space $X$ that was defined in 
\hyperref[negcurvinj]{Corollary \ref*{negcurvinj}(b)} by gluing a countable 
number of hyperbolic surfaces $F_i$ along a simple closed geodesic.

\begin{defi}\label{newSectDef1}
A curve $c: [a,b] \to X$ is straight, if it either coincides with $l$ and is
parametrized with constant speed, or if it falls apart into finitely
many segments   
           $c|_{[a_{i-1},a_i]}$, $(a=a_0 < a_1 <\ldots< a_n = b)$  such that
\begin{itemize}              
\item   each segment 
                $c|_{[a_{i-1},a_i]}$ is a hyperbolic geodesic in $F_{j(i)}$
         parametrized with the same constant speed for all $i \in \{1,....,n\}
         $.
         
\item $c(a_i) \in l$ ($i \in \{1,...,n-1\}$, only), and for some suitably
         small $\eps$ we have that $c((a_i - \eps , a_i)) \subset F_{j(i)}$ and
         $c((a_i,a_i + \eps )) \subset F_{j(i+1)}$, with $j(i) \neq 
         j(i+1)$.        
                  
\item    Assume that we equip $F_{j(i)}$ and $F_{j(i+1)}$ in the neighbourhood
         of $c(a_i)$ with some local orientation, so that this local orientation
         does not change along the segment              
            $ c|_{[a_i - \eps , a_i + \eps]}
       $ and is the basis of angle-measurement according to 
       \hyperref[newSectDef0]{Definition \ref*{newSectDef0}}. Then
       $$          \angle( c(a_i) , (c|_{[a_i - \eps , a_i]})^{-1},l^{-1}) =
         \angle( c(a_i) , c|_{[a_i, a_i + \eps ]},l)                     $$

        In words: a straight line segment is a hyperbolic straight line
        para\-metrized with the same constant speed in each surface $F_i$,
        and when our straight line approaches the common geodesic
        $l$ and there changes 
                     from one of our surfaces into another, it approaches 
        $l$ with that same angle at which it departs from $l$ in the other
        surface.
\end{itemize}        
Observe that it cannot be expected that each surface $F_i$ can be equipped
with a fixed orientation, so that the condition on local orientations
as demanded at the beginning of the third of the preceding items
would always be satisfied. Conversely, if we try to define a local
orientation in the neighbourhood of our straight line maintaining 
this condition whenever the straight line changes between 
surfaces, and our straight line should have self-intersections,
the local orientation might be mismatching at the self-intersections. 

We will have to show that any curve can be homotoped into a straight
position connecting the endpoints, and that this process
even extends to systems of curves forming higher dimensional simplices. 
In order to prepare for such this proof, we
will make use of the generalized universal covering $\widetilde{X}$
that $X$ possesses. Also, because unlike the case of the
proof of  \hyperref[negcurvinj]{Corollary \ref*{negcurvinj}(a)},
we do not have here an analogue to 
\hyperref[asphe]{Lemma \ref*{asphe}}, we will have to construct 
this space explicitly:  
\end{defi}      

\begin{construction}\label{newSectEndConstr4.3}
In order to obtain our generalized covering space $\widetilde{X}$, we first recall
that every surface-subspace $F_i$ has a universal cover $\widetilde{F_i}\approx 
\R^2$. Since, according to the classification of surfaces, the cyclic 
subgroup generated by the closed geodesic $l_i$ has infinite index, 
the plane $\widetilde{F_i}$ will contain countably many lifts of $l$,
that we enumerate as $\widetilde {l_{i,1}}, \widetilde {l_{i,2}}, 
\widetilde{l_{i,3}},\ldots$.
Because we did not set any assumption that our surfaces $F_i$ 
have to be homeomorphic to each other or how to identify them,
in case they should have the same genus and actually be homeomorphic
(and in case they should be homeomorphic: whether they have
the same hyperbolic structure and are isometric), we know that 
all $\widetilde{F_i}$ are hyperbolic planes and that way
isometric, but the distribution of the geodesics 
$\widetilde {l_{i,j}}\subset\widetilde
{F_i}$, can for each $i$ be a different one. Of course, for each 
$\widetilde{F_i}$, there is a covering map $p_i : \widetilde{F_i} \to
F_i$ defined. In addition, we choose  for our basepoint $P \in l$
one lift $\widetilde{P_{i,j}}\in \widetilde{l_{i,j}} \subset 
\widetilde{F_i}$. 
Next, we first define the following set of index-tuples
$$\Ups:= \{(k_1,k_2,\ldots,k_m) \mid m,k_j \in \N
\hbox{ satisfying the following conditions (1.)--(3 .)}\}$$
\begin{itemize}
\item[(1.)] $m$ is odd.
\item[(2.)] $k_j \neq k_{j-2}  \hbox{ for all odd } j \neq1.$  
\item[(3.)] $k_j \neq 1  \hbox{ for all even } j. $ 
\end{itemize}
Let $\Ups_i$ be 
the subset of tuples from $\Ups$, having the final index $k_m=i$. 
In order to construct our generalized universal
covering $\widetilde{X}$ we first define for every $\iota\in \Ups_i$
one isometric copy of $\widetilde{F_i}$
to be denoted by $\widetilde{F_\iota}$, so that ultimately for each $\iota\in \Ups$
some hyperbolic plane $\widetilde{F_i}$ gets associated with $i$
coinciding with the final index in the $\iota$-tuple. This final index
we will denote by $\omg(\iota)$, while the length of a $\iota$-tuple
will be denoted by $\ell(\iota)$. Removing the final index
index from $\iota$ and shortening the tuple accordingly will be denoted
by $\iota\setminus\omg(\iota)$. Our isometry $F_i =F_{\omg(\iota)}
\cong F_\iota$ also carries the definitions of the projections $p_i$, the geodesics $\widetilde{l_{i,j}}$ and
the points $\widetilde{P_{i,j}}$ into $F_\iota$, and accordingly we will
henceforth also use the notation $p_\iota$, $\widetilde{l_{\iota,j}}$ and
$\widetilde{P_{\iota,j}}$. Now 
$\widetilde{X}$ arises by making the appropriate identification in 
the subsequent disjoint sum
$$\widetilde{X}:= (\bigsqcup_{\iota\in\Ups}{\widetilde F_{\iota}})/\sim.$$
The relation ``$\sim$" identifies points of 
$\bigsqcup_{\iota\in\Ups}\widetilde{F_\iota}$ in the following
cases:  
\begin{itemize}
\item[(i)] If $\ell(\iota_1)=\ell(\iota_2)$ and if 
$\iota_1\setminus\omg(\iota_1)= \iota_2\setminus\omg(\iota_2)$, %=k>1$,
we identify the geodesic $\widetilde{l_{\iota_1,1}}\subset \widetilde{F_{\iota_1}}$
with $\widetilde{l_{\iota_2,1}}\subset \widetilde{F_{\iota_2}}$ in a way that
the metric and the lifted orientation from $l$ is preserved and in particular 
$\widetilde{P_{\iota_1,1}}$ and $\widetilde{P_{\iota_2,1}}$ get identified.
This will make sure that the projections
$p_{\iota_1}$ and $p_{\iota_2}$ agree along identified points. 
\item[(ii)] If $\ell(\iota_2)=\ell(\iota_1)+2$, also the penultimate index
of $\iota_2$ will be relevant and will be denoted by 
$\omg'(\iota_2)$. If in addition we have that 
$\iota_1 = (\iota_2\setminus\omg(\iota_2))\setminus(\omg'(\iota_2))$,
we identify the geodesic $\widetilde{l_{\iota_2,1}}\subset \widetilde{F_{\iota_2}}$
with $\widetilde{l_{\iota_1,\omg'(\iota_2)}}\subset \widetilde{F_{\iota_1}}$ in a way that 
the metric and the lifted orientation from $l$ is preserved and in particular 
$\widetilde{P_{\iota_1,\omg'(\iota_2)}}$ and $\widetilde{P_{\iota_2,1}}$ get identified.
This will make sure that the projections
$p_{\iota_1}$ and $p_{\iota_2}$ agree along identified points.
\end{itemize}
If we understand the symbol ``$\bigcup$" as 
``$\bigsqcup$ with the above described identification only",  then
$\widetilde{X} := \bigcup_{\iota\in\Ups}\widetilde{F_\iota}$  describes the point-set
of $\widetilde{X}$, and, since all identifications were non-contradicting
with respect to the $p_\iota$, we obtain a covering map 
$p : \widetilde{X} \to X$. The topology of $\widetilde{X}$ has been
defined according to \cite{fz} as the topology that was later called
``whisker topology" by \cite{DyCo}. Note that it does near
the attaching zones of $\widetilde{X} = \bigcup_{\iota\in\Ups} 
\widetilde{F_\iota}$ not coincide with the identification topology.

The combinatorics of $\widetilde{X} = \bigcup_{\iota\in\Ups}
\widetilde{F_\iota}$ is to some extent analogous to that of an 
infinitely branching tree, with the lifts of $l$ corresponding
to the vertices, and the planes $\widetilde{F_\iota}$ corresponding
to edges, only that our ``edges" connect countably many vertices. 
The geodesic arising from identification of the geodesics 
$\widetilde{l_{\iota,1}}$ for the $\iota$ with $\ell(\iota)=1$
is considered to be as the 
\begin{equation}
\hbox{``root" } \widetilde{l_0} \ni \widetilde{P_0}
\label{eq:root}
\end{equation}
of our tree-like structure, where $\widetilde{P_0}$
is a global notation for the points $\widetilde{P_{\iota,1}}$
that were to be identified when constructing $\widetilde{l_0}$.
We call 
\begin{equation}
(\ell(\iota)-1)/2 := d(\iota)\hbox{ the ``distance from the root"}
\label{eq:distfromroot}
\end{equation}
that the
plane $\widetilde{F_\iota}$ has. In all planes $\widetilde{F_\iota}$
with $d(\iota)\geq1$ the geodesic no.1 (i.e.\ $\widetilde{l_{\iota,1}}$)
is the only intersection with planes with lower distance to the
root (and therefore the ``exit into the direction of the root")
and with planes having  the same distance to the root, while 
at all other geodesics $\widetilde{l_{\iota,k} }$ $(k\neq1)$ planes 
which have by one higher distance to the root are intersecting.\par 
If $\iota = (k_1,\ldots,k_m)$, then a direct path from the  root
to the plane $\widetilde{F_\iota}$ can be described as follows:
\neuzl
When starting at the root, it goes into the plane lying over 
$F_{k_1}$, then at geodesic no.$k_2$ it changes into the plane
lying over $F_{k_3}$. In that plane at geodesic no.$k_4$
it changes into the plane lying over $F_{k_5}$, and that way it continues,
till all entries of our tuple are executed, and finally 
some plane over $F_{k_m}$ is reached. From that point of view the restrictions
(1.)--(3.) can be understood: The odd length of the tuples
arises from the alternate meaning of the symbols in our tuples, 
and (2.)--(3.) precisely ensure that at each geodesic arising
from identification exactly one plane lying over each of the surfaces
$F_i$ meets.\par
Similarly, as the generalized universal coverings
that we constructed in 
\hyperref[ExGenCov]{Lemma \ref*{ExGenCov}}, the space $\widetilde{X}$
may be interpreted as being homotopy equivalent to 
$\widetilde{X^{lpc}}$, i.e., 
if we forget, that there might be some accumulation between the
various $F_i$-layers of $X$ and that way obtain the space
$X^{lpc}$ of $X$ that is homotopy equivalent to a CW-complex, 
to the classical covering that this CW-complex has.
\par
That with $\widetilde{X}$ we actually constructed our generalized
universal covering of $X$ will be clear by the end of the following 
\hyperref[newSectLem2]{Remark \ref*{newSectLem2}}. It will follow, since by 
\hyperref[newSectLem2]{Remark \ref*{newSectLem2}} we can rule out such unwanted
effects as we saw in \hyperref[CAstrich]{Example \ref*{CAstrich}}, from
the fact that our space $\widetilde{X}$ has the same combinatorial structure as the classical
covering space $\widetilde{X^{lpc}}$, and due to the convention to use 
the whisker topology which ultimately guarantees
that around each of the lifts of $l$ we have the same topology 
as around $l\subset X$~--- and this altogether will make sure that
$p : \widetilde{X}\to X$ fulfils the unique path-lifting property,
which by \cite[Proposition 5.1]{fz} for first countable spaces
suffices to replace \hyperref[Def1.3]{Property \ref*{Def1.3}(ii)}.
\end{construction} 

\begin{remrk}\label{newSectLem2}
{\it Already by the construction of our space $\widetilde{X}$
it should be clear, that if this construction should have any chance 
to be the generalized universal covering space of $X$, it must
be true that $X$ does not permit to construct  any other paths than 
those which have finite or infinite segmentations where 
each of the segments exclusively lies in one of the $F_i$, and the endpoints 
of the segments (and maybe accumulation points of the segmentation)
belong to the joint geodesic $l$.} \par
In order to see this, please recall,  
that  by \hyperref[simplnu1po]{Remark \ref*{simplnu1po}}, the 
conclusions of \hyperref[simplnu1]{Lemma \ref*{simplnu1}}  can be used in our context, but that there we only
showed (\hyperref[simplnu1]{Lemma \ref*{simplnu1}(b)}) a weaker statement, namely that any path has a finite 
segmentation, and the segments alternately are belonging to one of the
$F_i$ or remained in the neighbourhood $U$ of the joint geodesic $l$,
where the existence of this neighbourhood $U$ with a special
topological type belonged to the assumptions of \hyperref[negcurvinj]{Corollary \ref*{negcurvinj}(b)}. Thus, what still
needs to be proven, is that if our path $c$ takes with $c(t)$ some point in
$U\setminus l$, i.e., for the appropriate $i$ we have  $c(t)  \in l_i 
\times(-1,+1)$, but distinct 
from $l$, then our path $c$  must spend an entire segment
$ c|_{(t-\eps,t+\eps)}$  only in $l_i \times (-1,+1)$ for some $\eps>0$. \par
Actually this conclusion follows analogously as our \hyperref[simplnu1]{Lemma \ref*{simplnu1}(b)} by 
just simplifying the argument:\neuzl
We take a closed neighbourhood $V$ of $c(t)$ which remains disjoint from 
$l$, and know by continuity that our path $c$ will for an 
appropriate $\eps$-neighbourhood of $t$ of the parameter domain remain
in this neighbourhood. If during this segment  $
c|_{(t-\eps,t+\eps)}$ our path should take values 
on other strips $l_j \times (0,1)$, 
by closeness of the $F_j$ the sets $c^{-1}(l_j \times (0,1))$ give a closed at most
countable covering of $(t-\eps,t+\eps)$. Now, with the paths-components
of this countable covering we can start to argue as in the proof of 
\hyperref[simplnu1]{Lemma \ref*{simplnu1}}, and unlike the proof of \hyperref[simplnu1]{Lemma \ref*{simplnu1}} there is no need
to look at some intersections components or deformation retractible
neighbourhoods of them. There is no need to enlarge this first 
covering somehow with other kinds of covering sets, just the
need, analogously as in the proof of \hyperref[simplnu1]{Lemma \ref*{simplnu1}}, to collapse those 
covering sets that are non-degenerate closed intervals to points,
and then we get a contradiction to Baire's theorem, unless we assume, that 
there was only one covering set, i.e.\ that the entire path $c$
on the entire interval $(t-\eps,t+\eps)$ remained
in its strip $l_i \times (-1,+1)$, i.e.\ in the surface $F_i$. Therefore we know
that, with respect to a finite segmentation of $c$ that comes from
\hyperref[simplnu1]{Lemma \ref*{simplnu1}}, each $F_i$-segments extents to a segment that has its
start- and its endpoint on $l$. 
                
     Now considering an entire $U$-segment of our curve $c$ again:
If a curve at the beginning and at the end of the corresponding $U$-segment
should have been in the same surface $F_i$, it might be, that it does
not take some value on $l$ at all. However, in this case the two
neighbouring $F_i$-segments and the $U$-segment in between could be just
united, to become one segment. If a curve $c$ at the beginning
and end of some $U$-segment belongs to different surfaces $F_i$, by what we have
just shown, it must have taken some value on $l$, and by closedness 
of $l$, there exists a first parameter $t_1$ and a last parameter $t_2$
($t_1 \le t_2$ with ``=" not being excluded) with $t_1 \in l$ and 
$t_2 \in l$ and $c([t_1,t_2]) \subset U$. We now want to homotope our curve
so, that after the homotopy we have $c([t_1,t_2]) \subset l$. 

Having with \hyperref[negcurvinj]{Conditions (i)--(ii)
from Corollary\ref*{negcurvinj}(b)} an assumption on the topological
type of $U$, namely that  
$$
  U \approx \hbox{ (countably many crossing open intervals 
                                        centered around 0) } \times l, 
$$
we know that the restriction of $c$ to $[t_1,t_2]$ can in 
principle be used to define coordinate functions for $c$
according to 
           $$    c(t) = (c_1(t),c_2(t),c(3)), \hbox{ where }$$ 
\begin{itemize}               
\item $c_1$ describes the $l$-coordinate  of $c(t)$,
\item $c_2$ is a real number that describes an element in one of the 
       crossing intervals according to where $c(t)$ takes its value, and
\item $c_3$ is the index $i$ so that $c(t)$ belongs to $F_i$ (or, as we denoted
       it at other places: to the strip ``interval $\times l_i$".). 
       This coordinate $c_3(t)$ may remain undefined for $t$ with 
       $c_2(t) = 0$, when by assumption $c(t)$ is at the crossing point. 
\end{itemize}

Based on these coordinates, we can now define the homotopy
pushing $c$ onto $l$, according to        
$$   h(t,s) :=   (c_1(t)\,,\, s \cdot c_2(t)\,,\, c_3(t)),$$
which will be the identity for $s=1$, and ensure that all values lie on 
$l$ for $s=0$. We do not claim that this construction will give
an overall continuous function on the entire neighbourhood $U$
(and in case of \hyperref[negcurvinj]{Corollary \ref*{negcurvinj}(b)} we do not have such a global assumption, 
i.e.\ we cannot rule out that the topology of the various
crossing line segments might look quite similar to that of our \hyperref[sin-touch]{Example \ref*{sin-touch}}),
but based on the fact that we know that $c$ is a continuous
path and that all points get just pushed towards $l$, we can conclude
that this homotopy will be a continuous function on its
parameter domain $[t_1,t_2] \times [0,1]$. 
\end{remrk}               

\begin{lemma}\label{newSectLem3}
Let $\T_X$ be the topology that we defined on $\widetilde{X}$ according
        to its construction in \hyperref[newSectEndConstr4.3]{Construction \ref*{newSectEndConstr4.3}}. The inner metric induced
        by path-lengths, when measured segment-wise 
        as long as a path was in one of the hyperbolic planes 
        $\widetilde{F_\iota}$  
         of $\widetilde{X}$
        also defines some topology $\T_\IH$ on $\widetilde{X}$. In general, $\T_\IH$ will
        be a richer topology than $\T_X$.
\end{lemma}
\neuzl{~}\neuzl
{\it Remark:} % (preceding the proof)
\neuzl{~}\par
In many cases it can be expected, that 
        $\T_X$ and $\T_\IH$ agree. However, 
         after the proof of this lemma we will show by presenting  
\hyperref[newSectExplEndLem3]{Example \ref*{newSectExplEndLem3}},
        that they need not necessarily agree.\par
Conversely, for paths in $\widetilde{X}$ (and other maps $Y \to \widetilde{X}$),
the inclusion of topologies to be shown in this lemma causes, 
that continuity with respect to $\T_\IH$ implies continuity with respect to 
$\T_X$ and therefore for such proof-constructions we can and will
sometimes work with the $\T_\IH$-topology. 
\begin{pf}        
Since the topology of the hyperbolic plane is known to be induced
       by such a metric, there is only one issue, that needs a closer
       analysis in order to agree to the claim of this lemma, 
       namely that the topologies are in the desired relation 
       along the lifts of $l$: It is clear that according to both definitions
       each point on $l$ has a neighbourhood of type:
\begin{equation}
\bigcup_i (P,Q) \times (-b_i,b_i)       
\label{eq:vierFis}   
\end{equation}
where $(P,Q)$ is a segment of $l$,  
        $b_i \in \R$,  $b_i > 0$, and the $b_i$-scale is adjusted in the 
        surface $F_i$, around $l$ perpendicular to $l$, so 
        that $|b_i|$  gives the hyperbolic distance to $l$. 
        The crucial question is just: Will the $b_i$-sequence be forced
        to have a positive infimum, or might it have 0 as infimum?
        According to the $\T_\IH$-convention ``inf=0" is definitely
        ruled out. In the classical identification topology it is,
        of course, permitted. Therefore, if $\T'_X$ denotes 
        the identification topology on $X$, we obtain 
        $\T'_X\supset \T_\IH$. However, here we have our Condition (ii) from
        \hyperref[negcurvinj]{Corollary \ref*{negcurvinj}(b)}
        in addition, and 
        due to this assumption we are able
        to show $\T_X\subset \T_\IH$ as follows:\neuzl
        Although 
        \hyperref[negcurvinj]{Condition \ref*{negcurvinj}(ii)}        
        is stated for neighbourhoods of the
entire closed geodesic $l$, since any local neighbourhood of some point
$P\in l$ can easily be extended to a global neighbourhood that still violates
this condition of minimal width, it also holds for any local neighbourhood. 
And then the minimal width as required by 
\hyperref[negcurvinj]{Condition \ref*{negcurvinj}(ii)}
gives a radius of a $\T_\IH$-neighbourhood, that squeezes into any $\T_X$-neighbourhood,
and its existence shows the desired inclusion $\T_X \subset \T_\IH$.
\end{pf}

As mentioned in the remark preceding the just finished proof,
only the subset-relations between the two topologies, but
not their equality can be proven, as the following example shows. 
        
\begin{expl}\label{newSectExplEndLem3}
  Let us assume, that at $l$ countably many surfaces are crossing,
        and let us further assume, that they are just enumerated
        by the naturals. Further assume, that each of them has a hyperbolic
        structure, and assume that 
        $\eps>0$ is fixed so that
        in each of the surfaces $F_i$ the $\eps$-neighbourhood
        around $l$ is just an ordinary band (i.e.\ that in particular 
        $2\eps$ is still smaller than the injectivity radius of each $F_i$).
        And then define the countable neighbourhood-basis, that according
        to one of the assumptions of \hyperref[thm1]{Theorem \ref*{thm1}} each point $\in X$ must have, 
        by giving, with respect to formula \hyperref[eq:vierFis]{ (\ref*{eq:vierFis})} the sequence
        of the $b_i$-values for the $k$-th neighbourhood of our countable
        basis as follows: 
$$ \Bigl({\eps\over k},{\eps\over k-1},{\eps\over k-2},\ldots,{\eps\over 3},
    {\eps\over 2},\eps,\eps,\eps,\ldots\Bigr)$$
    
    A neighbourhood from the above basis cannot be squeezed into
    any $\T_\IH$-neighbourhood, that has a diameter smaller
    than $\eps$. 
       On the other hand, demanding for a neighbourhood-basis
       these combinations of widths, with which each surface is 
       represented in each strip does not contradict to any other 
       assumption that was made in   \hyperref[newSectSam]{this Section \ref*{newSectSam}}.
 
Taking into account that by \cite[Proposition 5.1]{fz} the unique
path-lifting property (instead of the more general 
property % Def1.3(ii)
\hyperref[Def1.3]{Definition \ref*{Def1.3}(ii)})
decides, whether or not we
managed to construct a generalized covering space, we wish to 
point out that the subtle differences between these topologies do have
their impact which paths do or do not count as continuous:\neuzl   
In order to prepare for the definition of the corresponding path, 
let us first look at the function 
$f : [0,1] \to \R$
$$ f(x) :=
\left\{\vcenter{\halign{$#$~&#\hfil\cr        
0&for $x=0$\cr
 0.9 \eps \sin({\pi\over x})&for $x > 0$.
\cr}}\right.$$
       We now turn this real-valued function into a path 
       $c_f : [0,{1\over\pi}] \to X$ by fixing a point $P \in l$, and defining
       for each of the surfaces a perpendicular segment to $l$ through $P$
       of at least length $\eps$ in both directions. We demand that the
       path $c_f$ to be defined shall only take its values on these 
       segments. Once having agreed on this, the following
       defines $c_f$ (at least, up to some signs, which, however, are  irrelevant).  
 $$ c_f(t) :=
\left\{\vcenter{\halign{#~&#\hfil\cr
 $P$, &whenever $f(t)=0$,\cr
 $\subset F_m$&if $t \in [{1\over m+1} , {1\over m}]$,\cr
 satisfies&dist($c_f (t),P) = |f(t)|.$
\cr}}\right.$$ 
      Observe, that we have set up the definition of $c_f$ so, that
      it will not be continuous with respect to the $\T_\IH$-metric
      for $t=0$, but it counts as continuous with respect to the 
      $\T_X$-metric. This in particular implies that the $\T_X$-metric
      might indeed be more coarse than the more natural metric
      induced by hyperbolic geometry.     
\end{expl}                                

\begin{lemma}\label{newSectLem4}
Any path connecting two different points in $\widetilde{X}$ can be 
relatively homotoped
        into a subspace of $\widetilde{X}$ that is isometric to $\IH^2$.
\end{lemma}
\neuzl{~}\neuzl{\it
Remark:}        
\neuzl{~}\par
``Subspace isometric to $\IH^2$" does not mean         
     the ``hyperbolic planes $\widetilde{F_\iota}$" 
     out of which we built our space $\widetilde{X}$ in 
     \hyperref[newSectEndConstr4.3]{Construction \ref*{newSectEndConstr4.3}},
     and to which we will henceforth refer as ``$\IH^2$-plane(s)".
     Rather it can mean, e.g.,  a space that comprises a half-plane from some
     lift $\widetilde{F}_i$, and another half-plane from some other lift $\widetilde{F}_j$, 
     and if they should meet at some lift of $l$, it can mean
     the union of both half-spaces together with this lift of $l$. One might wonder,
     whether we have at such ``glueing line" still a differentiable
     structure; but if we metrize such a subspace with the inner metric
     induced by the hyperbolic metric on $\widetilde{F}_i$ and $\widetilde{F}_j$, the metric
     structure will not notice, that such an artificial gluing has taken
     place. Therefore we may call such constructions
     ``subspaces isometric to $\IH^2$". In addition, since all lifts of $l$ are disjoint,
     such a subspace isometric to $\IH^2$ can have any finite number of
     gluing lines inside, and even also a countable number.
     
\begin{pf} 
We start by applying the same mechanisms of simplification,
that we discussed in \hyperref[newSectLem2]{Remark \ref*{newSectLem2}}
for paths in $X$, to our space $\widetilde{X}$. That way,
we get in any case a finite segmentation. If for one of those segments
that our path spends in $\widetilde{F_\iota}$, it  should have entered and left 
$\widetilde{F}_\iota$ via the same lift of $l$, we can homotope
into the corresponding lift of $l$, and that way reduce
the segmentation. If neighbouring segments should belong to the same
$\iota$-index, by uniting them, we also simplify our segmentation. 
Due to  its finiteness, 
after at most finitely many of these simplification
steps of both kinds, this possibility is exhausted, and our path will 
on the one hand leave 
each of the $\widetilde{F}_\iota$-planes that it still enters through a different
lift of $l$ than it has entered it, and on the other hand spend neighbouring
$\widetilde{F}_\iota$-segments on planes with different indices
$\iota$ and $\iota'$. This definitely means, that 
$\omg(\iota)\neq\omg(\iota')$, because by 
\hyperref[newSectEndConstr4.3]{Construction \ref*{newSectEndConstr4.3}}, 
planes lying over the same surface $F_i$ do not meet $\widetilde{X}$. 
\par
Let $\iota_1, \iota_2, \iota_3,\ldots,\iota_n$ be the indices of the
$\IH^2$-planes that our path still sees after all simplifications
processes. For $\iota \in \{\iota_2,\ldots,\iota_{n-1}\}$,
we can then cut from 
$\widetilde{F_\iota}$ those half-spaces away that lie 
behind those lifts of $l$ where our path has entered, and where the path
has left again, treating the remainder as a ``hyperbolic strip".  
Simultaneously we can, even if our path should have on its way 
crossed the lifts of $l$ where it entered and left, homotope
it into a position, where its entire trace that it spends in
$\widetilde{F_\iota}$ 
lies in the hyperbolic strip that it has somehow to 
cross. Since each hyperbolic strip has two boundaries which are just
(infinitely long) hyperbolic straight lines, there is no obstruction
to glue the hyperbolic strips that are 
crossed by our path  
directly to each other, maintaining the same identifications,
that were also used for the construction of $\widetilde{X}$ in 
\hyperref[newSectEndConstr4.3]{Construction \ref*{newSectEndConstr4.3}}. 
In addition, we take 
from $\widetilde{F}_{\iota_1}$ (where the start-point of our path lies and where therefore
there is only one lift of $l$ in this plane  
distinguished by 
our path leaving the plane there), instead of a hyperbolic strip the 
entire half-plane bounded by this lift of $l$,  
and analogously 
for      $\widetilde{F}_{\iota_n}$, where there is the endpoint. There is no obstruction
to fill up an entire $\IH^2$-plane by isometrically embedding there
(disjointly, up to the necessary common boundaries) the half-plane
taken from $ \widetilde{F}_{\iota_1}$, then by the hyperbolic
strips taken consecutively from $ \widetilde{F}_{\iota_2}, \widetilde{F}_{\iota_3},\ldots, \widetilde{F}_{\iota_{n-1}}$, 
and finally the half-plane from $ \widetilde{F}_{\iota_n}$,
and then arrange these isometric embeddings, 
in such a way (just by maintaining the same identifications,
that were also used for the construction of $\widetilde{X}$ in 
\hyperref[newSectEndConstr4.3]{Construction \ref*{newSectEndConstr4.3}})
that the points on the common boundary where the path has crossed 
the corresponding lift of $l$ are embedded to the same point in $\IH^2$.
Since an $\IH^2$-plane
will be filled up that way with embeddings, we see, that if we just 
unite inside $\widetilde{X}$ the corresponding parts of the corresponding 
$\widetilde{F}_\iota$-planes, this will give a subspace isometric to $\IH^2$, and
having that way constructed this subspace, we have reached the
goal of this proof.
\end{pf} 
                                      
\begin{remrk}\label{newSectRem5} 
Recall that when stating \hyperref[negcurvinj]{Corollary \ref*{negcurvinj}(b)}, the existence of a  
deformation retractible neighbourhood in the sense of \hyperref[thm2]{Condition \ref*{thm2}(iii)}
was not demanded. On the other hand, it would have given an effective 
tool to reach the same conclusion as in
\hyperref[newSectLem2]{Remark \ref*{newSectLem2}}, 
namely deduct the finite segmentation of a homotopic version
of an arbitrary given path in $X$ into segments that alternately 
belong to one of the $F_i$ and inside $U$ directly approaches $l$,
may spend a segment on $l$, and then directly go into some
surface $F_j$ (with maybe $j\neq i$). We did not demand \hyperref[thm2]{Condition \ref*{thm2}(iii)}
for our example from \hyperref[negcurvinj]{Corollary \ref*{negcurvinj}(b)}, because we had the
assumption concerning the topological type of the neighbourhood
of the single intersection complex, and were     
able on that basis to conclude the desired result. 

It should be pointed out that the existence  of such a neighbourhood,
in which each surface $F_i$ is just represented by a strip of minimal 
width, is essential: It is well-known (cf.\ e.g.\ \cite[Thm.3.7.1]{Buser})
that one can consruct to any triple of positive real numbers
a pair of pants that has precisely these numbers as lengths 
of their geodesic boundaries, and by plugging such pairs of pants 
with matching boundary lengths together (as, e.g., shown in
\cite[Fig.3.1]{MonkThs}) one easily constructs hyperbolic surfaces
with arbitrarily short injectivity radius. Consequently,  
if our infinite collection of surfaces  
$F_\iota$ should contain such surfaces, where the return-length of a 
non-nullhomotopic path to the same point of $l$ is not bounded 
from below  (as it is by assumption of \hyperref[negcurvinj]{Corollary \ref*{negcurvinj}(b)} by
requiring the existence of the neighbourhood $U$ with its minimal-width 
assumption), non-trivial
passages through $l$ might accumulate along the parameter domain
of a path, and we would not be able to define the equivalent of the
angle-condition, that we required as last item of \hyperref[newSectDef1]{Definition \ref*{newSectDef1}},
at such points. 

When we argued in the 
\hyperref[newSectLem2]{Remark \ref*{newSectLem2}}
with the topological type of our neighbourhood $U$ instead of the
existence of a deformation retraction, this should also be considered
as a very hard geometric assumption, in particular with respect
to the remarks around equation \hyperref[eq:AnzTeilm]{(\ref*{eq:AnzTeilm})}, since this geometric assumption
again was, that a point lies either in the intersection of all $F_i$,
or only in one of these surfaces. Already \hyperref[stattviExpl]{Example \ref*{stattviExpl}} may be understood
as an example which shows that for other intersection scenarios, that
the structure of the intersection complexes can really get complicated;
and in particular \hyperref[thm2]{Condition \ref*{thm2}(v)} is not able to rule out similar unwanted
accumulation as for our \hyperref[CAstrich]{Example \ref*{CAstrich}}, if it takes place inside one of
the intersection complexes. 
\par
Therefore, although we did not need
it in our case, \hyperref[thm2]{Condition \ref*{thm2}(iii)} might also for situations
that are more close to \hyperref[thm1]{Theorem \ref*{thm1}} than to \hyperref[thm2]{Theorem \ref*{thm2}}, turn out
to be a useful condition to rule the most unwanted effects in the topology
transverse to the layers, similar as we pointed this out in the
penultimate paragraph before stating \hyperref[stattviExpl]{Example \ref*{stattviExpl}}.
\end{remrk}

\begin{lemma}\label{newSectLem6}
Any path in $X$ can in its relative homotopy class be uniquely  
represented by a straight 1-simplex. 
\end{lemma}

\begin{pf}
{\it (Existence:) }
The constant path counts as straight. In any other case
lift the path to $\widetilde{X}$. Use then \hyperref[newSectLem4]{Lemma \ref*{newSectLem4}} to find a subspace $H$ in $\widetilde{X}$ that 
contains a relatively homotopic version of our path and is isometric to $\IH^2$,
and then homotope the path into the position of the straight-line segment
connecting the start- and the end-point in this subspace, that is
isometric to $\IH^2$. If glueing-lines need to be crossed, then as well
the glueing-line as      our straight line are both represented
by straight lines with respect to the geometry of the hyperbolic plane, 
and so the angle-condition (i.e., the last item in the statement of
\hyperref[newSectDef1]{Definition \ref*{newSectDef1}}) will be satisfied.
\neuzl{~}\neuzl
{\it (Uniqueness:)}
Any two relatively homotopic straight connections in $X$ of the
same two points lift to straight connections in $\widetilde{X}$ between the
same two points. We look at the hyperbolic subspace $H$ that we
constructed in the existence proof for the first straight connection, and ask
how the second straight connection might appear there. If it should have
stayed in $H$, it would have also to look like
a straight line-segment, and therefore would have to coincide
with the first segment. \par    
As for options to find a straight way
that lies partially outside of $H$: Recall by the tree-like
way $\widetilde{X}$ was constructed (in particular cf.\
the paragraph preceding \hyperref[eq:root]{formula (\ref*{eq:root})}), 
that it is impossible to leave a subspace like
$H$ and to return to $H$ in any other way than through the same
lift of $l$. On the other hand, any path that performs such a turn,
must at least once enter an $\IH^2$-plane and leave it again through
the same lift of $l$ through which it entered. But when our path
is straight, it enters this $\IH^2$-plane crossing the lift of $l$
with a non-zero angle. The lift of $l$ is in the $\IH^2$-plane
a hyperbolic straight line, and our path emanates as a hyperbolic
straight ray from it with a non-zero angle and can by 
hyperbolic geometry not return
to this straight line. Since that way also straight connections partially 
outside of $H$ could be ruled out, uniqueness is proven.
\end{pf} 

\begin{defi}\label{newSectDef7}
A metric space which results by taking 
three hyperbolic half-planes and identifying them along their boundaries,
that are straight lines, we will call a {\it ``hyperbolic tripoid"}. 
Usually we will denote these half-planes by $H_1$, $H_2$ and $H_3$, and the
common straight-line segment by $g$. 
\end{defi}

\begin{remrk}\label{newSectRem8a}
If we try to imagine a hyperbolic tripod as an object 
in three-space, independently of whether and which of the 
hyperbolic models we use, we have to choose dihedral angles. 
Standard choices might be: either three times $120^\circ$, but also
once $180^\circ$ and twice $90^\circ$, so that two hyperbolic half-planes
together look like a hyperbolic plane, and a third sheet as
coming in from the side. But we may also try to imagine
one of the angles very small (i.e.\ ``$= \eps$"), and the others
to be $\pi - \eps/2$ so that two of the three possibilities to combine
two of the hyperbolic half-spaces to a hyperbolic space come
out almost looking like a hyperbolic half-space,  
and the third half-space, then of course, looks quite distorted
by being folded. However, in any case the choice of these dihedral angles
is not relevant for the geometry of such a tripod. The tripod is 
metrized by an inner metric induced by path-length, and the 
``branching to the side" in three-space has no impact on this inner metric.
Cf.\ the picture at the end of the proof of 
\hyperref[newSectLem9]{Lemma \ref*{newSectLem9}}.
\end{remrk}
 
\begin{remrk}\label{newSectRem8b}
 Recall that, as briefly indicated in the 
proof of \hyperref[newEndProp15]{Proposition \ref*{newEndProp15}}, 
having the straightening mechanism for 
lines, we can also define ``straight 2-simplices"
(= ``straight triangles"): \neuzl
If a curved  triangle is given
as a function $\sigma : \Delta^2 \to X$ (the vertices of the
standard 2-simplex $\Delta^2$ are denoted by $v_0$, $v_1$ and $v_2$,
and its images under $\sigma$ we define as $P_i := \sigma(v_{i-1})$),
then we define the straightened version of $\sigma$ 
by first straightening the path      
$ \sigma|_{[v_0,v_1]}$ to the straight 
segment $[P_1,P_2]$ as described in \hyperref[newSectLem4]{Lemma \ref*{newSectLem4}} and \hyperref[newSectLem6]{Lemma \ref*{newSectLem6}}, and then analogously
straightening the paths that are described by the corresponding restrictions
of $\sigma$ to the segments between $v_2$ and each of the points on the 
line-segment $[v_0,v_1]$ to straight segments that connect $P_3$
with any point on the straight segment $[P_1,P_2]$ that we constructed before. 
The straightened version of $\sigma$ we will call ``$str \circ \sigma$", and
we call any singular simplex ``straight", if it coincides from
the very beginning with its straightened version.  

Although straightening mechanisms can be defined
completely abstractly (cf.\ the remarks at the beginning
of \hyperref[4.4-5.4]{Section \ref*{4.4-5.4}}),  
we are only applying here the straightening mechanism as defined
for our space $X$ in \hyperref[newSectDef1]{Definition \ref*{newSectDef1}}.
\end{remrk} 
 
\begin{lemma}\label{newSectLem9} 
Any non-degenerate straight triangle $\sigma$ in our space $\widetilde{X}$
satisfies one of the following three conditions: 
\begin{itemize}
\item[(i)]
    it lies in a subspace of $\widetilde{X}$ that is isometric to $\IH^2$ and
    there the image of $\sigma$ is just the geometric triangle that has
    three different non-collinear points as vertices.
\item[(ii)]
{\it (``the non-folded tripod case")}\/ There exists a subspace of $\widetilde{X}$ that 
     is isometric to a hyperbolic tripod $H_1 \cup H_2 \cup H_3$,
     so that $\sigma(v_i) = P_{i+1} \in H_{i+1}$   $(i \in \{0,1,2\})$. 
     Let $Q$ be the point where the straight line-segment $[P_1,P_2]$
     intersects the common geodesic $g$ of the tripod. In this case
     our straight triangle is the union of two geometric triangles
     namely $[P_1,Q,P_3]$ that is contained in the hyperbolic 
     subspace $H_1 \cup H_3$, and $[P_2,P_3,Q]$ that is contained
     in the hyperbolic subspace $H_2 \cup H_3$. Now $im(\sigma) \cap H_3$
     will be again some hyperbolic triangle bounded by some segment
     of $g$, where the line $[P_3,Q]$ falls into the interior of this
     triangle. This line separates precisely those parts that lie
     in the image of the subtriangles $[P_1,Q,P_3]$ and $[P_2,P_3,Q]$.
\item[(iii)]     
{\it (``the folded tripod case")}. As in the preceding case, in this case again 
     $\sigma(\Delta^2)$ is the union of two geometric triangles,
     $[P_1,Q,P_3] \subset H_1 \cup H_2$ and $[P_2,P_3,Q] \subset H_2 \cup H_3$.
     As in the preceding case $im(\sigma) \cap H_3$ is a geometric
     triangle bounded by some segment of $g$. But in contrast to the 
     preceding case the line-segment $[P_3,Q]$ is one of the edges
     of this triangle, and the straight-line connections from 
     $P_3$ to $P_1$ and to $P_2$ emanate to the same side of this
     straight-line segment. That way a subtriangle of $im(\sigma) \cap H_3$
     has two preimages with opposite orientation in $\Delta^2$, one
     belonging to the subtriangle $[P_1,Q,P_3] \subset H_1 \cup H_2$, and the 
     other in the subtriangle $[P_2,P_3,Q] \subset H_2 \cup H_3$. This phenomenon
     we call ``fold". 
\end{itemize}
\end{lemma}                        

\begin{pf}
Let $P_1, P_2$ and $ P_3$ be given in $\widetilde{X}$, and let us trace the
     construction of the straight triangle between these points
     in order to see that only the above mentioned three cases can occur:
         According to \hyperref[newSectRem8b]{Remark \ref*{newSectRem8b}}%
, we first have to 
     straighten the edge $[P_1,P_2]$. By \hyperref[newSectLem4]{Lemma \ref*{newSectLem4}}, we find
     a subspace of $\widetilde{X}$, that is glued together from subspaces
     of $\IH^2$-subspaces, that are meeting at some lifts if $l$  in $\widetilde{X}$, and
     itself has also the metric of a hyperbolic plane. Let us denote
     this subspace by $H$ henceforth in this proof. In $H$ we can actually 
     realize the straight edge $[P_1,P_2]$ as a hyperbolic straight 
     line-segment. 
     
     Then we have the following possibilities:
\begin{itemize}     
\item[(1.)] 
the third vertex $P_3$ of $im(\sigma)$ 
   might already lie in $H$. If this case should
   be given, then clearly we are in case (i), and our straightening
   mechanism will just generate the system of straight lines that
   connect the point $P_3$ inside $H$  
   with the straight line segment $[P_1,P_2]$, and the area covered by all
   these lines together will just give a triangle as in hyperbolic
   geometry. 
\end{itemize}
However, if $H$ should not contain $P_3$,
of course somewhere in $\widetilde{X}$ exists the system of straight lines
that connects $[P_1,P_2]$ with $P_3$. And since all these lines
start in our subspace $H$, we can ask, where they leave $H$.
In principle they can only leave $H$ via a lift of $l$, and 
by the tree-like way, how $X$ is built (cf.\ 
the paragraph preceding  
\hyperref[eq:root]{formula (\ref*{eq:root})}
and 
the uniqueness-proof of
\hyperref[newSectLem6]{Lemma \ref*{newSectLem6}}
) it has
for all of them to be the same lift of $l$ that will be denoted by 
$\widetilde{l}$.\par
Having agreed on that, there are just two further possibilities:
\begin{itemize}
\item[(2.)]
       $\tilde l$ remains disjoint from $[P_1,P_2]$ --- and
\item[(3.)] there is an intersection point $\tilde l \cap [P_1,P_2]$.
\end{itemize} 

\noi\underbar{Ad (2.):} In this case we will still get case (i), but 
   in contrast to (1.) the subspace isometric to $\IH^2$ in which 
   straightened triangle $[P_1,P_2,P_3]$ lies is not $H$, but some other subspace $H'$ 
   that we will construct in the following. 
      Since $\tilde l$ does not intersect $[P_1,P_2]$, it subdivides $H$ into
   two halfspaces, one of which entirely contains $[P_1,P_2]$. 
   We construct $H'$ by taking this half-space of $H$, and replacing the
   other half-space of $H$ by other parts of $\widetilde{X}$. As such parts we will
   naturally choose the parts of $\widetilde{X}$ through which our straight lines
   pass between their intersection with $\tilde l$ on their way to $P_3$.
   With the same arguments that we used in the paragraph      after having 
   discussed Case~(1.), 
all these lines have after 
   their intersection with $\tilde l$ on their way
   to $P^3$ to enter the same $\IH^2$-planes, and if the corresponding
   $\IH^2$-plane should not contain $P_3$ already, leave it 
   at another (but the same for all our lines) lift of $l$ again,
   entering there another (but for all lines the same $\IH^2$-plane),
   and such changes of $\IH^2$-plane might occur a finite number of times
   and be performed analogously for all or lines till they finally 
   reach $P_3$.  \par
   Having that way by tracing our lines found their 
   ``joint way", we can, analogously as we constructed in \hyperref[newSectLem4]{Lemma \ref*{newSectLem4}} our 
   subspace $H$, construct the partial replacement of $H$
   in order to obtain $H'$. Having cut $H$ along $\tilde l$ down to a hyperbolic
   half-space, we fill the empty half up by a finite number of
   hyperbolic strips, according to the joint way that our straight
   line segments take on their way to $P_3$, and close up the filling
   by the half-space, that is taken from that half-space
   containing $P_3$, bounded by that lift of $l$ through 
   which our straight lines have entered that $\IH^2$-plane.
   In the so-constructed subspace $H'$ our straight triangle
   $[P_1,P_2,P_3] = im(str \circ \sigma)$ will lie.
   
\noi\underbar{Ad(3.):}   
           Assume that the system of straight lines leading to $P_3$
   leave $H$ at a lift of $l$ (again denoted by $\tilde l$) that 
   intersects $[P_1,P_2]$. Denote the intersection point by $Q$. 
      Formally, in discussing this case, we still might distinguish two 
   subcases, namely $\tilde l$ might be one of the glueing lines that was
   used when composing $H$, or it might not be. However, we do not need to
   discuss these subcases separately, although in the second of these subcases,
   the straight segment $[P_1,P_2]$ was, when crossing $\tilde l$ not 
   changing its $\IH^2$-subplane, but in the first subcase
   it did precisely that. On the other hand, both subcases have in common,
   that the system of lines from $[P_1,P_2]$ to $P_3$ does go into 
   another of the countably many subspaces hyperbolic half-planes
   that meet at $\tilde l$, than those containing the straight-line segments
   $[P_1,Q]$ or $[Q,P_2]$.
   
   So, when we now have to define our hyperbolic tripod, we will
   take $g := \tilde l$ as the common straight line, $H_1$ will be that half-space
   of $H$ cut along $\tilde l$ containing $[P_1,Q]$, and $H_2$ be the complementary
   half-space of $H$. The missing $H_3$ will then have
   to be constructed by tracing the straight line-segments from $[P_1,P_2]$ 
   on their way to $P_3$ after their passage through $\tilde l$. In constructing 
   $H_3$ its first segment (at least some hyperbolic strip) will be precisely
   that segment into which our straight segments from 
   $[P_1,P_2]$ to $P_3$ disappear from $H = H_1 \cup H_2$ at $\tilde l$. 
   With similar arguments as we used in the paragraphs after having 
   discussed Case~(1.),  from the way $\widetilde{X}$ is built  
   we conclude, that all these segments after having disappeared from $H$
   into the next hyperbolic strip will have on their way to $P_3$ change
   $\IH^2$-planes at precisely the same geodesics in the same order
   (because they could not be homotopic while keeping $P_3$ fixed
   and sliding the other endpoint only via $[P_1,P_2]$, and simultaneously be 
   straight lines, otherwise). 
   Apart from that, their way  
   through $H_1 \cup H_2$ will for similar reasons be, 
   that they from their start-point
   to their meeting point with $\tilde l$ have to cross the same lifts of $l$ 
   in the same order, as the corresponding subsegment of $[P_1,Q]$ or
   $[P_2,Q]$ between this start-point and $Q$. And since after having crossed
   $\tilde l$, all these paths have to cross the same lifts of $l$ in the
   same order and change between the same $\IH^2$-subspaces of $\widetilde{X}$, 
   it clear that via glueing hyperbolic strips together and 
   a hyperbolic half-plane that contains $P_3$, similarly as we have
   constructed $H$ in the existence-proof for \hyperref[newSectLem6]{Lemma \ref*{newSectLem6}}, 
   we can find a hyperbolic half-space
   in $\widetilde{X}$ that takes up all the corresponding segments  of our 
   straight-line segments and prove that in these cases at least
   a hyperbolic tripod as a subspace of $\widetilde{X}$ containing our straight
   triangle $[P_1,P_2,P_3]$ exists. \par 
      Now we have to analyse how a straightened triangle in our 
   hyperbolic tripod looks like. This becomes best visible, if we, 
   as discussed in \hyperref[newSectRem8a]{Remark \ref*{newSectRem8a}}, use the freedom to imagine $H_1 \cup H_2$
   as almost almost completely folded along the line $\tilde l$ (i.e.\ so that
   the dihedral angle between $H_1$ and $H_2$ is really a small $\eps$, and
   imagine it maybe as the lower half-space of our projection from the
   3-dimensional embedding, while $H_3$ we imagine as the upper 
   half-space of our projection from three dimensions, so that the
   dihedral angles between $H_3$ and $H_1$ or $H_2$ are $\pi-\eps/2$. 
      Let us first discuss the subcase, where our straight-line segment
   $[P_1,P_2]$ intersects $\tilde l$ at a right angle: How the straight triangle
   $[P_1,P_2,P_3]$ in this case will look like, depends on, where in 
   $H_3$ the point $P_3$ lies. Recalling that we decided to call  
   $Q := \tilde l \cap [P_1,P_2]$, if $P_3$ should lie precisely on the perpendicular
   raised at the point $Q$, then the entire triangle will just be degenerate.
   If it should lie to the left or to the right of this perpendicular 
   raised at the point $Q$, then we are precisely in case (iii)
   of the lemma to be proven here, that we called a ``folded tripod
   triangle". Namely, it is just from the geometric assumptions
   that we have made up to this moment clear, that from three emanating 
   straight connections from $P_3$ (to $Q$, $P_1$ and $P_2$) the one
   to $Q$ cannot be in the middle between those  $P_1$ and $P_2$,
   but it must be the one approaching our perpendicular line most quickly,
   and therefore we will just obtain a folded triangle in this case.\neuzl  
      Let us now discuss the general case where $\angle(Q,[P_1,P_2],l) = \beta 
   \neq \pi/2$. Since we did not fix any orientation on either of the
   $H_i$, w.l.o.g.\ $0<\beta<\pi/2$. 
   \par
Naturally $l$ is subdivided by $Q$ into two 
   (semi-infinite) subsegments, and
   one of them will be more close to $[P_1,Q]$ than to $[Q,P_2]$, and
   for the other subsegment it is vice versa. We will call the
   first subsegment $l_1$, the second $l_2$, and henceforth assume that
   $l$ (and also $l_1 \& l_2$) are oriented in the direction from $l_1$
   to $l_2$. We then draw two auxiliary rays $r_1$ and $ r_2$ into $H_3$, both starting
   at the point Q, and being determined by the following angle-conditions:
   $$
         \angle(Q,r_2, l_2) = \beta~~\hbox{and}~~
         \angle(Q,r_1, l_1^{-1}) = \beta.
$$         
                                     %     v---- was: green    
In Euclidean perspective this would be the accordingly labled rays  
in the figure below:

%  Beim n\"achsten Bild:  green ---> 0.2   black, opacity=0.2
%                         cyan ---> 0.6    black, opacity=0.6

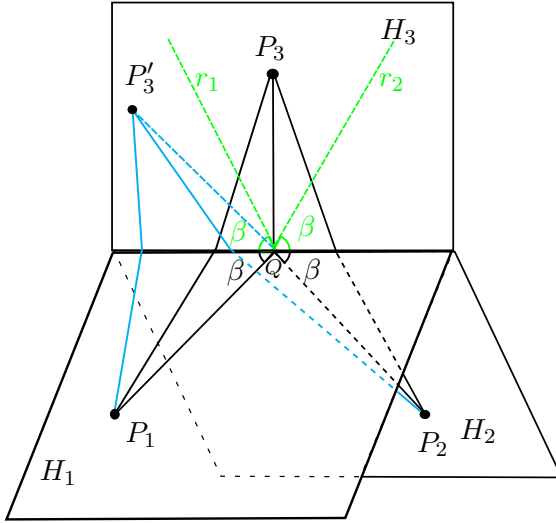
\def \globalscale {2}
\begin{tikzpicture}[y=1cm, x=1cm, yscale=\globalscale,xscale=\globalscale, every node/.append style={scale=1}, inner sep=0pt, outer sep=0pt]

  % H2_2
  \path[draw=black,fill=white,line width=0.0232cm,cm={ 0.0,-0.8552,1.7789,-0.0078,(-0.7317, 2.5722)}] (1.8682, 1.8048) -- (0.1215, 1.409) -- (0.1194, 2.6692) -- (1.8661, 3.065) -- cycle;

  % Right upper angle
  \path[draw=green,fill opacity=0.2667,line width=0.0255cm,rotate around={-180.0:(0.0, 4.5446)}] (-2.8716, 6.5392)arc(116.0001:185.2001:0.105 and -0.1023) -- (-2.8256, 6.6311) -- cycle;

  % H2_1
  \path[draw=black,fill=white,line width=0.0322cm,cm={ 0.0,-1.0102,-1.7686,-0.0092,(6.4882, 2.5876)}] (1.8682, 1.8048) -- (0.1215, 1.409) -- (0.1194, 2.6692) -- (1.8661, 3.065) -- cycle;

  % Dashed plane
  \path[draw=black,line width=0.0163cm,dash pattern=on 0.0653cm off 0.1306cm] (3.676, 2.4505) -- (1.7816, 2.4426) -- (2.4788, 0.9605) -- (3.3968, 0.9571) -- (4.0031, 2.4494) -- cycle;

  % H2_3 ?
  \path[draw=black,line width=0.0234cm] (1.762, 4.1007) rectangle (4.0066, 2.4576);

  % Text nodes
  % H2_1
  \node[anchor=south west,line width=0.0163cm] (text7) at (1.2869, 0.9046){$H_1$};

  % H2_2
  \node[anchor=south west,line width=0.0163cm] (text8) at (4.0511, 1.1788){$H_2$};

  % H23
  \node[anchor=south west,line width=0.0163cm] (text9) at (3.5257, 3.8288){$H_3$};

  % Triangles
  \path[draw=black,line width=0.0255cm] (2.8088, 3.6371) -- (2.4471, 2.4723) -- (1.7782, 1.3705);

    % Edge P2 path no 3
    \path[draw=black,line width=0.0234cm] (3.678, 1.6362) -- (3.8235, 1.3694);

    \path[draw=black,line width=0.0234cm,dash pattern=on 0.0702cm off 0.0702cm] (3.2416, 2.4365) -- (3.387, 2.1698) -- (3.5325, 1.903) -- (3.678, 1.6362);

    \path[draw=black,line width=0.0234cm] (2.8234, 3.6384) -- (3.2416, 2.4365);

  \path[draw=black,line width=0.0238cm] (2.8207, 3.6231) -- (2.8253, 2.4384) -- (1.787, 1.3764);

  % Upper left angle
  \path[draw=green,fill opacity=0.2667,line width=0.0255cm,xscale=1.0,yscale=-1.0,shift={(0,-9.084562708333334)}] (2.7936, 6.536)arc(112.9999:182.2:0.105 and -0.1023) -- (2.8346, 6.6302) -- cycle;

  % Edge - p2 path no 2
  \path[draw=black,line width=0.0234cm] (3.635, 1.5757) -- (3.697, 1.5087) -- (3.7589, 1.4417) -- (3.8208, 1.3748);

  \path[draw=black,line width=0.0234cm,dash pattern=on 0.0702cm off 0.0702cm] (2.8297, 2.4464) -- (2.8916, 2.3794) -- (2.9536, 2.3124) -- (3.0155, 2.2455) -- (3.0775, 2.1785) -- (3.1394, 2.1115) -- (3.2014, 2.0445) -- (3.2633, 1.9775) -- (3.3253, 1.9106) -- (3.3872, 1.8436) -- (3.4492, 1.7766) -- (3.5111, 1.7096) -- (3.5731, 1.6427) -- (3.635, 1.5757);

 % Point
  \path[draw=black,fill=black,line width=0.0078cm] (2.8189, 3.6269) ellipse (0.0365cm and 0.0319cm);

  % Left bottom algle
  \path[draw=black,line width=0.0255cm] (2.7657, 2.3792)arc(131.0:179.2:0.105 and -0.1023) -- (2.8346, 2.4564) -- cycle;

  % Right bottom angle
  \path[draw=black,line width=0.0255cm,xscale=-1.0,yscale=1.0] (-2.8965, 2.3761)arc(131.0002:179.2001:0.105 and -0.1023) -- (-2.8276, 2.4533) -- cycle;

  % Rignt lower beta
  \node[text=black,anchor=south west,line width=0.0255cm,shift={(-0.12,-0.26)}] (text19) at (2.5781, 2.352){$\beta$};

  % Left lower beta
  \node[text=black,anchor=south west,line width=0.0255cm,shift={(0,-0.25)}] (text19-5) at (3.0196, 2.3556){$\beta$};

  \node[text=green,anchor=south west,line width=0.0255cm] (text19-5-0) at (2.9821, 2.5149){$\beta$};

  \node[anchor=south west,line width=0.0255cm] (text22) at (1.8504, 1.1559){$P_1$};

  \node[anchor=south west,line width=0.0255cm] (text23) at (3.7694, 1.095){$P_2$};

  \node[anchor=south west,line width=0.0255cm] (text25) at (2.6984, 3.7199){$P_3$};

  \node[anchor=south west,line width=0.0255cm] (text26) at (1.8431, 3.5037){$P_3'$};

  \path[draw=green,line width=0.0226cm,dash pattern=on 0.0904cm off 0.0226cm] (2.131, 3.8597) -- (2.8292, 2.4855)(2.8292, 2.4855) -- (3.6239, 3.8586);

  \node[text=green,anchor=south west,line width=0.0255cm,shift={(-0.2,-0.1)}] (text19-2) at (2.6315, 2.531){$\beta$};

  \node[text=green,anchor=south west,line width=0.0255cm,shift={(+0.0,+0.0)}] (text27) at (3.51615, 3.51615){$r_2$};
  
   \node[text=green,anchor=south west,line width=0.0255cm,shift={(+0.0,+0.0)}] (text28) at (2.30555, 3.51615){$r_1$};
  
  \node[text=black,below,line width=0.0255cm,shift={(+0.0,-0.18)}] (text29) at (2.8292, 2.4855){$_Q$};

  % Edge path to p2 no 1

  \path[draw=cyan,fill opacity=0.2667,line width=0.0255cm] (3.6449, 1.5161) -- (3.8031, 1.3801);

  \path[draw=cyan,fill opacity=0.2667,line width=0.0255cm,dash pattern=on 0.0764cm off 0.0764cm] (2.5377, 2.4678) -- (2.6959, 2.3319) -- (2.8541, 2.1959) -- (3.0122, 2.0599) -- (3.1704, 1.924) -- (3.3286, 1.788) -- (3.4867, 1.6521) -- (3.6449, 1.5161);

  \path[draw=cyan,fill opacity=0.2667,line width=0.0255cm] (1.7758, 1.3748) -- (1.9589, 2.4651) -- (1.897, 3.3886) -- (2.5377, 2.4678);

  \path[draw=cyan,fill opacity=0.2667,line width=0.0255cm,dash pattern=on 0.1018cm off 0.0255cm] (1.8862, 3.3993) -- (2.8177, 2.4786);

   % Points
  \path[draw=black,fill=black,line width=0.0063cm] (1.8956, 3.3909) ellipse (0.0286cm and 0.0266cm);

  \path[draw=black,fill=black,line width=0.0069cm] (1.7801, 1.3714) ellipse (0.0288cm and 0.0323cm);

  \path[draw=black,fill=black,line width=0.0067cm] (3.8204, 1.3699) ellipse (0.0316cm and 0.027cm);

\end{tikzpicture}

\par\noi
   The two rays $r_1$ and $r_2$ subdivide $H_3$ into three domains that we 
   will call $V_0$, $V_1$ and $V_2$, where $V_0$ will be the domain 
   bounded by $r_1$ and $r_2$, and the other domains $V_i$ are bounded
   by $r_i$ and $l_i$. Please observe that the construction is made
   so that $[P_1,Q] \cup r_2$ and $[P_2,Q] \cup r_1$ are by definition
   straight rays in $\widetilde{X}$. Based on this fact it is clear, depending on  
   where in $H_3$ our point $P_3$ will lie (in particular 
   with respect to the subdivision into the  
   domains $V_0$, $V_1$ and $V_2$), the straight triangle $[P_1,P_2,P_3]$ will
   look differently, because the straight connections from $P_3$ to 
   $P_1$, $P_2$ and $Q$ leave $P_3$ in different orders (with respect
   to angle-measurement):
\begin{itemize}
\item  If $P_3 \in V_0$, we get the case of an (unfolded or ``ordinary")
       hyperbolic tripod-triangle as described in case (ii),
       from a Euclidean perspective, e.g., the black triangle
       $[P_1,P_2,P_3]$ in the above picture.       
       
\item  If $P_3 \in V_1 \cup V_2$, we get folded tripod-triangles,  
       as described in case (iii), from an analogous perspective, e.g.,
       the triangle $[P_1,P_2,P'_3]$, whose fold-line $[P'_3,Q]$
       is drawn dashed.   % was: "blue" 
\end{itemize}
\end{pf}

\begin{lemma}\label{newSectLem10}
Let $X$ be a space satisfying the assumptions of \hyperref[negcurvinj]{Corollary \ref*{negcurvinj}(b)}.
Then the Gromov norm of every non-nullhomologous 
cycle $z \in H_2(X;\R)$ is bounded away from zero. 
\end{lemma}

\begin{pf} 
Recall that we have analysed the algebraic structure of this
homology group in \hyperref[newEndProp14]{Proposition \ref*{newEndProp14}(b)} and in this course already defined 
the index-sets $S$, $T$, $T^+$ and $T^-$. Since $H_2(X;\R)$
was a direct sum, our cycle $z$ can only have non-zero components 
in finitely many summands, and without loss of generality we assume
that they have indices $\in \{1,2,\ldots, m\}$. Let us use
the abbreviations $S_{\le m}$ for $S \cap \{1,2,....,m\}$
and $T_{\le m}$, $T^+_{\le m}$ and $T^-_{\le m}$ analogously.
By what we have shown in \hyperref[newEndProp14]{Proposition \ref*{newEndProp14}},
using square brackets to denote the fundamental cycle,  
we can therefore write
$$ z = \sum_{i\in S_{\le m}} \alpha_i\cdot[F_i] +
       \sum_{i\in T_{\le m}} \beta^+_i\cdot[F^+_i] +
       \sum_{i\in T_{\le m}} \beta^-_i\cdot[F^-_i]  $$
with the additional condition, that $ \sum_{i\in T_{\le m}} \beta^+_i-\beta^-_i
=0$. Assume now, that $c = \sum_{j\in J} \gamma_j\cdot\sigma_j$ 
is a singular chain representing $z$. When attempting to show that 
there exists some lower bound for $\Vert c\Vert_1$, we can, without
loss of generality assume that $c$ is straight, because straightening
does not increase the $l^1$-norm. As opposed to the
situation that we had when proving \hyperref[newEndProp15]{Proposition \ref*{newEndProp15}}, we cannot
assume that straightening alone suffices to make sure, that each
straight singular simplex $\sigma_j$ would go essentially just into
one of our surfaces $F_i$, $F^+_i$ or $F^-_i$; we have analysed in 
the proof of \hyperref[newSectLem9]{Lemma \ref*{newSectLem9}} the situation to that extent, that it can in principle
go into any finite number of $F_i$. 
\par
However, based on this analysis,
we can subdivide any straight simplex further, along those lines, where, 
as discussed in the proofs of  
\hyperref[newSectLem4]{Lemma \ref*{newSectLem4}}
and \hyperref[newSectLem9]{Lemma \ref*{newSectLem9}},
the simplex changed between the
various $F_i$-surfaces. This subdivides our simplices usually not into 
triangles, but into convex $n$-gons. For not having to invent a new
homology theory (where continuous images of $n$-gons could occur
as generators of homology groups), we decide to subdivide these 
convex $n$-gons by inserting $n-3$ diameters between the vertices,
so that each $n$-gon gets subdivided into triangles. These ``subtriangles"
lie in one hyperbolic plane and  are therefore in any case geometric triangles, 
whose volume can be measured. However they might not fulfil the classification
that we proved in \hyperref[newSectLem9]{Lemma \ref*{newSectLem9}}, because by being cut out from a straightened
triangle, it might well be, 
that the straightening mechanism did not work in that way, that we first
pulled one edge straight, and then all line-connections from the opposed
vertex of our subtriangle to the pulled-straight edge. Such a mechanism
was applied to the original triangle $\sigma_j$, but not to the 
triangles $\tau_\nu$ that arise from this additional 
subdivision process. If one of these subtriangles should get subdivided
by the fold-line of a tripod-triangle, 
(cf.\ \hyperref[newSectLem9]{Lemma \ref*{newSectLem9}(iii)})
we will 
additionally subdivide this triangle along its fold-line. 

Having defined the set of according subdivision triangles $\tau_\nu$,
we hereby decide to avoid double-indexing; noting that the subdivision
is in any case finite, we know that we get a finite set $K $ 
that enumerates all the obtained subdivision triangles, $\tau_\nu$,
when $\nu$ runs through $K $. Following these 
principles we obtain a homologous representation $c'$ of our cycle $z$ by
$$ c \sim  c' = \sum_{\nu\in K }\eta_\nu \cdot \tau_\nu$$
Naturally each index $\nu$ is clearly associated with one index 
$j\in J$ by the fact that $\tau_\nu$ is defined by being cut out 
from the straight simplex $c_j$, and we will denote this index as
 ``$j(\nu)$". In addition, each index $\nu$ is associated to some 
index  $i\in \{1,2,\ldots,m\}$ by the fact that this subdivision triangle sits 
in the corresponding surface $F_i$, and we denote this index as ``$i(\nu)$". 
Actually, in case $i(\nu)\in T$, we treat $i(\nu)$ even to that 
extent as defined (just by the fact whether $\tau_{i(\nu)}$ lies in
$F_i^+$ or $F_i^-$, that it is even defined whether $i(\nu)$ is an
element of $T^+$ or of $T^-$. Since they have just come from subdividing a
chain, the new coefficients $\eta_\nu$ are 
in principle just the old coefficients $\gamma_j$ and it holds that
$\eta_\nu = \gamma_{j(\nu)}$. The new chain $c'$ has now the 
advantage, that it can be subdivided into subchains that can be separately
interpreted in the separate $F_i$, where we have a better control
over  their Gromov norm. On the other hand it has to be noted, that $c'$,
as resulting from the subdivision of $c$, will have a much higher $l^1$-norm
than $c$. Conversely, even if we have to sum over the $c'$-summands,
we can still obtain the old $l^1$-norm $\Vert c\Vert_1$, when summing
\begin{equation}
\Vert c\Vert = \sum_{\nu\in K } |\eta_\nu|\cdot {\Vol(\tau_\nu)
                                    \over           \Vol(\sigma_{j(\nu)})}
\label{eq:zweiStar}                                         
\end{equation}                                                    
Here ``$\Vol$" denotes the two-dimensional hyperbolic volume
(or ``area"). The preceding formula holds, 
because, when taking the subsums of those $\nu$ having the same
$j(\nu)$-value, the $\Vol$-ratios will sum up to $1$, and the
identical coefficients $\eta_\nu$ can then be replaced by 
$\gamma_{j(\nu)}$. 

We now consider the components of our cycle in the corresponding 
surfaces $F_i$, $F^+_i$ or $F^-_i$:
$$ z_i \sim c_i = \sum_{\nu\in K\cap\{\nu\mid i(\nu)\in S\}}
                         \eta_\nu \cdot \tau_\nu~~,~~                  
 z^\pm_i \sim c^\pm_i = \sum_{\nu\in K\cap\{\nu\mid i(\nu)\in T^\pm\}}
                         \eta_\nu \cdot \tau_\nu.                   $$                      
Based on the fact that by assumptions and constructions we are led
to the conclusion, that these cycles $c_i$ and $c_i^\pm$ need to realize a
corresponding multiple of the fundamental cycle of $F_i$ or 
$F_i^\pm$, the following inequalities are to be understood that way, 
that on the left hand side the volume of the corresponding surface
with its multiple according to the cycle $z$ is given, and that the 
right hand side, if we first ignore that last factor, would describe 
precisely the same multiple of this volume if we assume that the 
subdivision triangles 
$\tau_\nu$ are optimally placed, and cover the entire area
of $F_i$ or $F_i^\pm$ without some subarea being covered twice
with reverse orientation. Of course, this assumption need not be
fulfilled, but if it should be not fulfilled, since the chain on the
right hand side represents the correct multiple of the 
fundamental cycle $[F_i]$ or $[F_i^\pm]$,
on the left-hand side, the entire area of $F_i$ or $F_i^\pm$ 
must be covered up, and so the only other assumption to the one
from  the preceding sentence would be, that some areas are 
multiply covered with reverse orientation, so that in the algebraic
homology-count some cancellation will take place, but on right 
hand side, where we computed with absolute values, accordingly not,
and therefore the following inequalities will hold as real inequalities
(still ignoring the last factor):
\begin{equation}
\vcenter{$$
|\alpha_i| \cdot \Vol(F_i) \le \sum_{\nu : i(\nu)=i} |\eta_\nu| \cdot
\Vol(\tau_\nu) \cdot { 2\cdot v_2 \over \Vol(\sigma_{j(\nu)})}
\hbox{ for } i\in S   ,      $$$$ |\beta^\pm_i|  \cdot \Vol(F^\pm _i) \le \sum_{\nu 
: i=i(\nu)\in  T^\pm} |\eta_\nu| \cdot \Vol(\tau_\nu) \cdot { 2\cdot v_2 
\over \Vol(\sigma_{j(\nu)})}              \hbox{ for } i \in T^\pm  . $$}
\label{eq:dreiStar}
\end{equation}
In order to justify now the last factor, recall that 
by showing in \hyperref[newSectLem9]{Lemma \ref*{newSectLem9}} that even the tripod-triangles are just
unions of two geometric triangles, 
we showed that even in this context no straightened triangle can have
more volume than $2\cdot v_2$, where $v_2=\pi$ denotes the volume
of the regular ideal triangle. Therefore the quotient on the right hand
of the above inequalities is bigger than $1$, and so does it
not overturn the inequality that we have justified before, 
while ignoring it.
Having that way shown \hyperref[eq:dreiStar]{(\ref*{eq:dreiStar})}, we continue our estimates as follows:\neuzl
We bring ``$2\cdot v_2$" to the left hand side, and then observe
that on the right hand side the quotients between 
$\Vol(\tau_\nu)$ and $\Vol(\sigma_{j(\nu)})$ have already been considered
in equation \hyperref[eq:zweiStar]{(\ref*{eq:zweiStar})}. That way we see: If we now add the inequalities
\hyperref[eq:dreiStar]{(\ref*{eq:dreiStar})} over all $i\in S \sqcup T^+ \sqcup T^-$, the left hand side
consists only of constants that are independent of our concrete 
chain $c$, and apart from that it adds up to a positive value, since
by assumption that $z \not\sim 0$ not all $\alpha_i$ and
$\beta_i^\pm$ can be zero. Conversely, on the right-hand side we get 
in principle a sum over all $\nu\in  K$. For that sum
the subsums of those $\nu$ with the same $j(\nu)$-value $j$, as
discussed in the context of discussing \hyperref[eq:zweiStar]{ (\ref*{eq:zweiStar})}, are precisely adding up 
to $|\gamma_j|$. Consequently, the entire sum adds up to the $l^1$-norm
of our chain $c$. These considerations give the desired
estimate that $\Vert z\Vert\neq0$.
\end{pf}
\neuzl{~}\neuzl{\it      
Completing the proof of \hyperref[negcurvinj]{Corollary \ref*{negcurvinj}(b):}}  
\neuzl{~}\par  
In order to finally show that \hyperref[thm1]{Theorem \ref*{thm1}} can be applied to give the
desired properties of the class of spaces described in \hyperref[negcurvinj]{Corollary \ref*{negcurvinj}(b)}, 
it is necessary to verify that this class of spaces satisfies the
assumptions of \hyperref[thm1]{Theorem \ref*{thm1}}:
   It is clear that the spaces that we considered are second countable
$T_1$-spaces. The question of non-trivial Gromov norm is by \hyperref[newEndProp14]{Proposition \ref*{newEndProp14}(b)} only 
relevant for homology-dimension two, where it was verified in \hyperref[newSectLem10]{Lemma \ref*{newSectLem10}}. 
We constructed the generalized universal covering 
$\widetilde{X}$ of our spaces 
in 
\hyperref[newSectEndConstr4.3]{Construction \ref*{newSectEndConstr4.3}}, 
and thus the only aspect still to be verified is its contractibility.

In order 
to construct the corresponding contraction, i.e.\ 
the corresponding system of maps $r_t$ with 
$r_0 = id$ and $r_1$ projecting everything to one point 
$\widetilde{P_0}$ (cf.\ 
\hyperref[eq:root]{(\ref{eq:root})}--\hyperref[eq:distfromroot]{(\ref{eq:distfromroot})})
we subdivide the interval for the homotopy parameter $t \in [0,1]$
by the unit fractions, and define the maps $r_t$ for 
$t \in [{1\over m+1} , {1\over m}]$ so that,
\begin{itemize}
\item $r_t$, when restricted to those planes $\widetilde{F_\iota}$
          with $d(\iota)< m-1$
         will be the identity. 
\item If $x$ belongs to one of those planes $\widetilde{F_\iota}$
          with $d(\iota)= m-1$, then $r_t(x)$ is to be defined so, 
          that the  points will
        wander via lines orthogonal to $\widetilde{l_{\iota,1}}$
        inside the plane                  so that $r_{1/(m+1)}$ is still
        the identity on those planes, and $r_{1/m}$ projects the entire
        plane to $\widetilde{l_{\iota,1}}$. For the planes $\widetilde{F_\iota}$
          with $d(\iota)= 0$ exceptionally the points wander radially to the
        point $\widetilde{P_0}$ instead of orthogonally to the root-line $
        \widetilde{l_0}$. Since the root-line is the only of our 
        intersection lines, where only planes of the same distance
        to the root intersect, such a definition is possible while 
        non-contradicting on this line. 
\item If $x\in\widetilde{F_\iota}$
          with $d(\iota)\ge m$ we define $r_t$ for $t \in [{1\over m+1} , 
          {1\over m}]$ by the following formula:
          $$r_t(x) := (r_t \circ r_{1\over m+2}\circ r_{1\over m+3}
          \circ\ldots\circ  r_{1\over d(\iota)-2}\circ r_{1\over d(\iota)-1}
          \circ r_{1\over d(\iota)}\circ r_{1\over d(\iota)+1})(x)$$
          Observe that any of the middle-terms $r_{1/k} $
          in the above superposition takes its argument 
          from some geodesic $\widetilde{l_{\iota',k'}}$ with 
          $d(\iota')=k+1$ and projects it to the geodesic 
          $\widetilde {l_{\iota',1}} \subset
          \widetilde{F_{\iota'}}$, which by the identifications made
          in \hyperref[newSectEndConstr4.3]{Construction \ref*{newSectEndConstr4.3}}
          also belongs to some plane with distance $d(\iota')-1$
          to the root. Therefore each of these retractions $
          r_{1/ k}$ takes its argument so, that it is in the range
          as defined in the preceding item. Similar arguments also hold
          for the first and the last member of the superposition
          on the right-hand side of the above equation, so that 
          this term is by the prior item already defined and therefore
          able to define the left-hand side.
\end{itemize}
The geometric idea behind the above definition is that each point
has to wander due the orthogonal or radial projections that have
been used in the above item-list and the direct path that we described
in the  penpenultimate paragraph of %Constr.4.3
\hyperref[newSectEndConstr4.3]{Construction \ref*{newSectEndConstr4.3}}
towards the root, while all points for small
homotopy parameters $t$ first wait, till points on planes 
with bigger distance to the root have on their tour reached at least their
plane. 
It is easily noted that this definition is sufficiently consistent to be
continuous for all parameters $t \in (0,1]$. However, because every point
lies in one plane with  distance $m$ to the root (or in the worst 
case on one intersection line with two adjacent distance
parameters), 
the definition implies that our point and all neighbouring points have
to be mapped identically by $r_t$ for all $t \in [0,{1\over m+2}]$, and
therefore this construction is also continuous for the homotopy 
parameter zero. 

Having shown that the generalized covering space that we constructed 
in \hyperref[newSectEndConstr4.3]{Construction \ref*{newSectEndConstr4.3}} is actually contractible, we verified all conditions of 
\hyperref[thm1]{Theorem \ref*{thm1}}, and can by \hyperref[thm1]{Theorem \ref*{thm1}} now conclude that the desired injectivity
of the canonical homomorphism now follows. \blackbox

\section{Reduction to simplices with all vertices in the basepoint}\label{onepoint}
\subsection{Eilenberg's argument}\label{eilenberg}
For a topological space $X$ with basepoint $x_0$ we denote by $C_*(X)$ the complex of singular simplices, i.e., the chain 
complex whose $k$-th group is the free abelian group generated by $S_k(X)=map(\Delta^k,X)$ with the usual boundary operator 
and by $C_*^{x_0}(X)\subset C_*(X)$ the subcomplex
generated by 
$$S_k^{x_0}(X)=\left\{\sigma\colon\Delta^k\to X\mid \sigma(v_0)=\sigma(v_1)=\ldots=\sigma(v_k)=x_0\right\},$$
where $v_0,v_1,\ldots,v_k$ denote the vertices of the standard simplex $\Delta^k$.

It is a classical result of Eilenberg (Corollary 31.2 in \cite{eil}) that for path-connected $X$ the inclusion 
$$\iota\colon C_*^{x_0}(X)\to C_*(X)$$
is a chain homotopy equivalence. It is well-known that this dualizes to give chain homotopy equivalences 
also in cohomology and bounded cohomology. In this section we are going to show that (under a suitable assumption) the argument also yields chain homotopy equivalences for measure homology and measurable bounded cohomology. 

Let us start with recalling Eilenberg's argument (which in \cite{eil} is given in a more general setting).

\begin{lemma}\label{lem1}
For each path-connected space, there is a chain map $\eta_*\colon C_*(X)\to C_*^{x_0}(X)$ such that $\eta\iota=id$ and a chain homotopy $s_*\colon C_*(X)\to C_{*+1}(X)$ such that $$\partial s+s\partial=\iota\eta-id.$$ 
\end{lemma}
\begin{pf}
For $x\in S_0(X)$ we have to define $\eta_0(x)=x_0$.

Because $X$ is path-connected we have a $1$-simplex $s_0(x)\colon \Delta^1\to X$ with $$\partial_0s_0(x)=x_0\,,\;\;\partial_1s_0(x)=x$$
for each $x\in X$. Let us fix a choice of $s_0(x)$ for each $x$.

Now we define $\eta_*$ and $s_*$ by induction on the dimension of simplices. Suppose they are already defined for all simplices in $S_{k-1}(X)$ and let $\sigma\in S_k(X)$. By induction hypothesis we have $\eta_{k-1}(\partial\sigma)\in C_{k-1}^{x_0}(X)$ and $s_{k-1}(\partial\sigma)\in C_k(X)$ such that 
$$\eta_{k-1}(\partial\sigma)-\partial\sigma=\partial s_{k-1}(\partial\sigma)+s_{k-2}\partial(\partial\sigma)=
\partial s_{k-1}(\partial\sigma).$$ 
We will inductively prove the slightly stronger statement that $s_k$ is of the form $s_k=s_k^0+\ldots +s_k^k$ and that the maps $s_k^0,\ldots,s_k^k$ can be defined through some map $F\colon \Delta^k\times\left[0,1\right]\to X$ via the canonical subdivision $$\Delta^k\times\left[0,1\right]=\Delta_0\cup\ldots\cup\Delta_k$$ as the restrictions of $F$ to $\Delta_0,\ldots,\Delta_k$. Recall that 
the simplex $\Delta_j$ from\break \hl{the}
canonical subdivision is defined as the convex hull of the points\break   
$(v_0,0),(v_1,0),\ldots,(v_j,0),(v_j,1),(v_{j+1},1),\ldots,(v_k,1)$,  
where $v_0,v_1,v_2,\ldots,v_k$ is the list of vertices of the
standard $k$-simplex $\Delta^k$.

So consider $\Delta^k\times\left[0,1\right]$. We can use $\sigma$ to define a continuous map $\Delta^k\times\left\{0\right\}\to X$ and by the above inductive hypothesis we have
$$s_{k-1}(\partial\sigma)=(s_{k-1}^0+\ldots+s_{k-1}^{k-1})(\partial\sigma)$$ defined through a continuous map $\partial \Delta^k\times\left[0,1\right]\to X$. These two maps agree on $\partial \Delta^k\times\left\{0\right\}$, so they define a continuous map $$Q\colon \Delta^k\times\left\{0\right\}\cup\partial \Delta^k\times\left[0,1\right]\to X.$$

It is easy to construct a continuous map $$P\colon \Delta^k\times\left[0,1\right]\to  \Delta^k\times\left\{0\right\}\cup\partial \Delta^k\times\left[0,1\right]$$ which is the identity map on $\Delta^k\times\left\{0\right\}\cup\partial \Delta^k\times\left[0,1\right]$. We can compose $P$ with the before-defined map $Q$ to obtain a continuous map $$F\colon\Delta^k\times\left[0,1\right]\to X$$ that on $\Delta^k\times\left\{0\right\}\cup\partial \Delta^k\times\left[0,1\right]$ agrees with $Q$. 
We use the canonical triangulation of $\Delta^k\times\left[0,1\right]$ 
into $k+1$ simplices to consider $F$ as a formal sum of $k+1$ simplices, which we denote by $s_k^0(\sigma),\ldots,s_k^k(\sigma)$. We obtain  
thus an element $$s_k(\sigma):=s_k^0(\sigma)+\ldots+s_k^k(\sigma)\in C_{k+1}(X).$$
In particular 
$F\vert_{\Delta^k\times\left\{1\right\}}$ defines
$\eta_k(\sigma)\in C_k(X)$ which actually belongs to $C_k^{x_0}(X)$ because all vertices are in $x_0$. It is then clear by 
construction that the equality $\partial s_k(\sigma)+s_{k-1}(\partial\sigma)=\eta_k(\sigma)-\sigma$ holds. 
\end{pf}

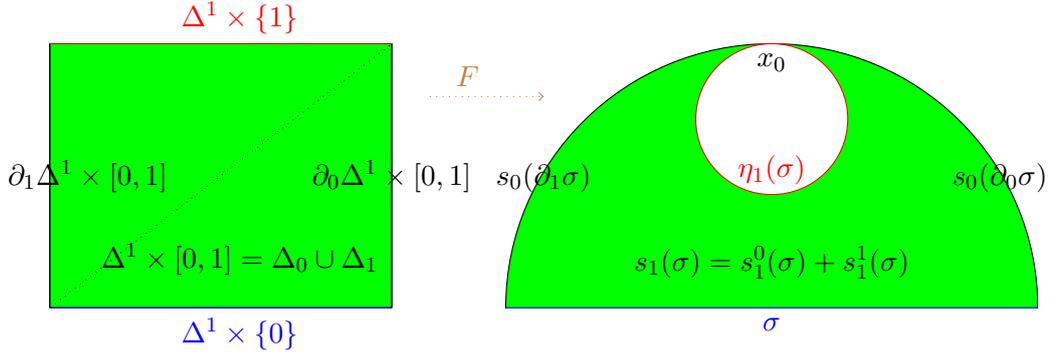
\begin{tikzpicture}
\node[above,red] at (1.5,3.5){$\Delta^1\times\left\{1\right\}$};
\draw[fill=green](-1,0)--(3.5,0)--(3.5,3.5)--(-1,3.5)--(-1,0);
\draw[blue] (-1,0)--(3.5,0);
\node[below,blue] at (1.5,0) {$\Delta^1\times\left\{0\right\}$};
\draw(-1,0)--(-1,3.5);
\draw(3.5,0)--(3.5,3.5);
\node[above] at (-0.5,1.4) {$\partial_1\Delta^1\times\left[0,1\right]$};
\node[above] at (3.5,1.4) {$\partial_0\Delta^1\times\left[0,1\right]$};
\node[above]at(1.5,0.3) {$\Delta^1\times\left[0,1\right]=\Delta_0\cup\Delta_1$};
\draw[red](-1,3.5)--(3.5,3.5);
\draw[dotted](-1,0)--(3.5,3.5);

\draw[brown,dotted,->](4,2.8)--(5.5,2.8);
\node[above,brown] at (4.5,2.8){$F$};

\begin{scope}
    \clip (5,0) rectangle (12,4);
    \draw[fill=green] (8.5,0) circle(3.5);
    \draw[blue] (5,0)--(12,0);
\end{scope}

\node[below,blue] at (8.5,0) {$\sigma$};
\draw[red,fill=white](8.5,2.5)circle(1);
\node[above,red] at (8.5,1.5) {$\eta_1(\sigma)$};
\node[below] at (8.5,3.5) {$x_0$};

\node[above] at (5.5,1.4) {$s_0(\partial_1\sigma)$};

\node[above] at (11.5,1.4) {$s_0(\partial_0\sigma)$};
\node at (8.5,0.6){$s_1(\sigma)=s_1^0(\sigma)+s_1^1(\sigma)$};

\end{tikzpicture}
The figure visualizes
the construction of the map $F$ in case of $k = 1$. \hl{It shows in particular} that
the two vertices (endpoints of the simplex $\Delta^1 \times \left\{1\right\}$) are taken under $F$
to
the same point $x_0\in X$.

\subsection{Pointed measure homology}\label{pmh}

We now want to argue that an analogous result as in \hyperref[lem1]{Lemma \ref*{lem1}} 
holds for measure homology, i.e., that (under suitable assumptions) the inclusion
$$\iota\colon{\mathcal C}^{x_0}_*(X)\to {\mathcal C}_*(X)$$
is a chain homotopy equivalence. Here ${\mathcal C}^{x_0}_*(X)\subset {\mathcal C}_*(X)$ means the subcomplex consisting of those signed measures (of \hl{compact} determination set and bounded variation) which vanish on each measurable subset of the complement of $S_*^{x_0}(X)$.

\begin{lemma}\label{lem2}
If $X$ is a path-connected space that has a finite covering $$X=\bigcup_{i=1}^n U_i$$ such that 

- $U_1,\ldots,U_n$ are Borel-measurable sets

- the closures $\overline{U}_1,\ldots,\overline{U}_n$ are contractible in $X$ and compact, 

\noindent then for any $x_0\in X$ there is a chain map 
$\eta_*\colon \C_*(X)\to \C_*^{x_0}(X)$ such that $\eta\iota=id$, and a chain homotopy 
$s_*\colon \C_*(X)\to \C_{*+1}(X)$ such that $$\partial  s+s\partial =\iota\eta-id.$$ 
\end{lemma}

\begin{pf}
The natural approach to proving this statement would be to define $\eta$ and $s$ as in the proof of 
\hyperref[lem1]{Lemma \ref*{lem1}}. One would have to check then 
that signed measures of compact determination set and bounded variation are mapped to signed measures of compact determination set and bounded variation.

It is clear that a so-constructed $\eta_k$ does not increase the variation and that $s_k$ multiplies the variation by 
at most $k+1$, so the second condition on boundedness of the variation will be satisfied.

To satisfy the first condition on compactness of the determination set it would be sufficient that 
$\eta_k$ and the maps $s_k^0,\ldots,s_k^k$ from the proof of 
\hyperref[lem1]{Lemma \ref*{lem1}} could be defined via some continuous maps on $map(\Delta^k,X)$, because then 
compact determination sets of simplices would be mapped to 
compact sets. In general it will not be possible to define such a continuous map. It would be possible if $X$ were contractible. It is still possible on subsets that are contractible in $X$ and our argument will make use of this fact.

Let $\overline{U}\subset X$ be contractible in $X$. Then there is some continuous map $H:\overline{U}\times\left[0,1\right]\to X$ with $H(x,0)=x$ and $H(x,1)=x_0$ for all $x\in \overline{U}$. Define $$s_0\colon map(\Delta^0,\overline{U})\to map(\Delta^1,X)$$ by $$s_0(x)(t)=H(x,t)$$ 
upon identification $\Delta^1=\left[0,1\right]$. Continuity of $H$ and compactness of $\left[0,1\right]$ imply that $s_0$ is continuous. 

\hl{Now consider the covering $X=\cup_{i\in I} U_i$ (which exists by our assumptions)} by finitely many Borel sets whose closures are compact and contractible in $X$. ($I$ is a finite index set.) W.l.o.g. we can assume that the $U_i$ are disjoint. Indeed, if they were not, we could replace $U_i$ by $V_i=U_i\setminus \cup_{j=1}^{i-1}U_j$ for $i\ge 2$. The closures $\overline{V}_i$ are subsets of $\overline{U}_i$ and hence again compact and contractible in $X$ (although not necessarily in $\overline{U}_i$), and of course the $V_i$ are again Borel sets.

For each ordered $(k+1)$-tuple $(i_0,\ldots,i_k)$ of (not necessarily distinct) 
elements of the index set $I$ we let $S_{i_0,\ldots,i_k}$ be the set of singular simplices with $0$-th vertex in $U_{i_0}$, 1-st vertex in $U_{i_1}$, ..., $k$-th vertex in $U_{i_k}$ and we consider its closure $\overline{S}_{i_0,\ldots,i_k}$ which is contained in the set of singular simplices with $0$-th vertex in $\overline{U}_{i_0}$, 1-st vertex in $\overline{U}_{i_1}$, ..., $k$-th vertex in $\overline{U}_{i_k}$. 

By the above we have defined $\eta_0$ and $s_0$ 
on $S_0=U_0,\ldots,S_k=U_k$ (i.e., on all of $X$), such that the restriction to each $S_i$ extends continuously to $\overline{S}_i$.

Now we assume by induction that for all 
$k$-tuples $(i_0,\ldots,i_{k-1})$ we already have maps 
$$\eta_{k-1}\colon {S}_{i_0,\ldots,i_{k-1}}\to map(\Delta^{k-1},X)$$ and 
$$s^0_{k-1},\ldots,s_{k-1}^{k-1}\colon {S}_{i_0,\ldots,i_{k-1}} \to map(\Delta^k,X)$$
with the desired properties and which all extend continuously to $\overline{S}_{i_0,\ldots,i_{k-1}}$. We claim that $\eta_k$ and $s_k^0,\ldots,s_k^k$ (defined as in the proof of \hyperref[lem1]{Lemma \ref*{lem1}}) are again continuous maps on $\overline{S}_{i_0,\ldots,i_k}$ for each $(k+1)$-tuple $(i_0,\ldots,i_k)$. 

This is seen as follows. Continuity of $s^0_{k-1},\ldots,s_{k-1}^{k-1}$ implies 
that the map $$\overline{S}_{i_0,\ldots,i_{k-1}}\to map(\Delta^k\times\left\{0\right\}\cup\partial \Delta^k\times\left[0,1\right],X)$$ 
which sends $\sigma\colon\Delta^k\to X$ to the ``union" of $\sigma\times\left\{0\right\}$ and 
$s^0_{k-1}(\partial_j\sigma),\ldots,s^{k-1}_{k-1}(\partial_j\sigma), j=0,\ldots,k$, is continuous. 
Moreover, precomposition with the uniformly continuous map $P\colon \Delta^k\times\left[0,1\right]\to  
\Delta^k\times\left\{0\right\}\cup\partial \Delta^k\times\left[0,1\right]$ from the proof of 
\hyperref[lem1]{Lemma \ref*{lem1}}
defines a continuous map  
$$map(\Delta^k\times\left\{0\right\}\cup\partial \Delta^k\times\left[0,1\right],X)\to 
map(\Delta^k\times\left[0,1\right],X),$$
so we obtain a continuous map $$\Phi\colon 
\overline{S}_{i_0,\ldots,i_{k-1}}
\to 
map(\Delta^k\times\left[0,1\right],X).$$ 
Since $\eta_k(\sigma)$ and $s^0_k(\sigma),\ldots,s_k^k(\sigma)$ are all 
defined by restricting $\Phi(\sigma)$ to subsets of $\Delta^k\times\left[0,1\right]$, they also depend 
continuously on $\sigma$.

So we have proved that $\eta_k$ and $s_k^0,\ldots,s_k^k$ (defined on $S_{i_0,\ldots,i_k}$)
can be extended continuously to $\overline{S}_{i_0,\ldots,i_k}$ (although this extension on $\overline{S}_{i_0,\ldots,i_k}\setminus S_{i_0,\ldots,i_k}$ of course does not have to agree with the actual definition of $\eta_k$ and $s_k^0,\ldots,s_k^k$ coming from some other $S_{j_0,\ldots,j_k}$). Since all the $S_{i_0,\ldots,i_k}$ are pairwise disjoint, this allows a (not continuous but measurable) definition of $\eta_k$ and $s_k^0,\ldots,s_k^k$ on $$map(\Delta^k,X)=\bigcup_{(i_0,\ldots,i_k)}S_{i_0,\ldots,i_k}.$$
For any compact subset $K\subset map(\Delta^k,X)$ we obtain 
that the image of  $K\cap S_{i_0,\ldots,i_k}$ under $\eta$ or $s_0,\ldots,s_k$ is contained in the image of $K\cap \overline{S}_{i_0,\ldots,i_k}$ under some continuous extension of $\eta_k$ or $s_k^0,\ldots,s_k^k$
and thus is contained in a compact set. Hence the image of $K\cap S_{i_0,\ldots,i_k}$ has compact closure. 

So the image of $K$ under any of $\eta$ and $s_0,\ldots,s_k$ is a finite union of (subsets of) compact sets, hence has compact closure. 

In particular, because the image of a determination set under any map is a determination set for the push-forward measure, $\eta_k$ and $s_k^0,\ldots,s_k^k$ map measures of compact determination set to measures of compact determination set.
\end{pf}

\begin{cor}\label{corx0}Under the assumptions of \hyperref[lem1]{Lemma \ref*{lem2}} every measure cycle is homologous to a measure cycle with determination set contained in $S_*^{x_0}(X)$.\end{cor}

Recall that we have defined bounded cohomology in \hyperref[relinjbc]{Definition \ref*{relinjbc}} and
measurable bounded cohomology in \hyperref[measbc]{Definition \ref*{measbc}}. Similarly 
one defines 
measurable cohomology.

Let us denote by $H^*_{b,x_0}(X), \mathcal{H}_{x_0}^*(X)$ and $\mathcal{H}_{b,x_0}^*(X)$ the cohomology groups of the complexes of bounded, measurable resp.\ bounded measurable functions from $C_*^{x_0}(X)$ to $\R$. Using \cite[Section 3.4]{loe} there is 
a well-defined pairing between $\mathcal{H}_{b,x_0}^*(X)$ and ${\mathcal{H}}_*^{x_0}(X)$.
\begin{cor}\label{boundedx0}Under the assumptions of \hyperref[lem1]{Lemma \ref*{lem2}}, the canonical restriction induces isomorphisms
$$H^*_b(X)\to H^*_{b,x_0}(X)$$
$${\mathcal H}^*_b(X)\to \mathcal{ H}^*_{b,x_0}(X)$$
$${\mathcal H}^*(X)\to {\mathcal H}^*_{x_0}(X)$$\end{cor}
\begin{pf}The above constructed maps $\eta$ and $s$ are bounded in the sense that $\eta_k$ sends a simplex to a simplex and $s_k$ sends a $k$-dimensional simplex to a formal sum of (at most) $k+1$ $(k+1)$-dimensional simplices. This implies that $\eta^*$ and $s^*$ send bounded cochains to bounded cochains. Moreover $\eta_k$ and $s_k$ are continuous on each of the finitely many disjoint Borel sets $S_{i_0\ldots i_k}$, so they are Borel-measurable on $map(\Delta^k,X)$ and hence $\eta^*$ and $s^*$ send measurable cochains to measurable cochains.\end{pf}

\subsection{Examples}
Let us conclude with some examples fulfilling or not fulfilling the assumptions of \hyperref[lem2]{Lemma \ref*{lem2}}:
\begin{example}{\em CW-complexes}

Any compact manifold or finite CW-complex can be covered by finitely many measurable sets with contractible, compact closures. 
Thus the assumptions of \hyperref[lem2]{Lemma \ref*{lem2}}
are satisfied.
\end{example}

\begin{example}{\em Earring space}

The earring space is the shrinking wedge of circles pictured in the introduction, that is, it can be written
in the form $$HE=\bigcup_{n=1}^\infty C_n\subset\R^2,$$
where $C_n\subset\R^2$ is the circle with centre $(0,\frac{1}{n})$ and radius $\frac{1}{n}$.
Let $A_n^\pm$ be the intersection of $C_n$ with the closure of the upper resp.\ lower half-plane, and let $A^\pm=\bigcup_{n=1}^\infty
A_n^\infty$. Then
$$HE =A^+\cup A^-$$
is a covering by two measurable sets with contractible, compact closures. Thus the assumptions of \hyperref[lem2]{Lemma \ref*{lem2}}
are satisfied for the earring space.
\end{example}

\begin{example}{\em Warsaw circle}

The Warsaw circle is a closed subset $W\subset\R^2$, which is the union of the graph of the function
$$y=\sin(\frac{1}{x})$$
for $0<x\le 1$, the segment 
$$Y=\left\{(x,y)\colon x=0,-1\le y\le 1\right\},$$
and a curve connecting these two parts to get a path-connected space.

This space can be covered by finitely many contractible, measurable, relatively compact sets.
The easiest way to do this is to use the decomposition $$W=Y\cup Y^c$$
into $Y$ and its complement. However the closure of $Y^c$ is all of $W$, which
is known to be not contractible.

On the other hand, $W$ can be covered by countably many contractible, compact sets. For this one has to decompose the graph 
of $y=\sin(\frac{1}{x})$
into its segments for $\frac{1}{n+1}\le x\le\frac{1}{n}$ with $n$ running through
all natural numbers, and then add $Y$ and the connecting curve as two more 
contractible, compact sets to the decomposition.

These two decompositions show that in \hyperref[lem2]{Lemma \ref*{lem2}} the assumption on having contractible closures and the assumption on finiteness of the covering can not
be relaxed by just assuming contractibility of the relatively compact sets themselves or by countability of the covering, respectively.
Indeed for the Warsaw circle $W$, the second author proved in 
\cite[Theorem 4]{pr2} that $\mathcal{H}_0({W})$ is uncountable-dimensional, while of course \hl{$\mathcal{H}_0^{x_0}({W})\cong\R$}. 

The Warsaw circle does however not provide a counterexample to the conclusion of \hyperref[thm2]{Theorem \ref*{thm2}}
in view of $H_0(W;\R)=\R$ and $H_n(W;\R)=0$ for all $n>0$. 

\end{example}

\begin{expl}\label{noninj}{\em A space with non-injective canonical homomorphism}

The following space
can be covered by two contractible sets, but they are not Borel-measurable.
   
Let $Z$ be the space constructed in \cite[Section 5]{pz}. 
There are two points $z_0,z_1\in Z$ such that $$\left[z_1\right]-\left[z_0\right]\in ker(H_0(Z;\R)\to\mathcal{H}_0(Z)),$$ 
see 
\cite[Theorem 5.7]{pz}. 

To get examples in any degree $n$, let $F$ be a closed, orientable manifold of dimension $n\ge 1$ and consider $F\times Z$ 
with the ``boundary" manifolds $F_0=F\times\left\{z_0\right\}$ and 
$F_1=F\times\left\{z_1\right\}$. Then 
$$\left[F_1\right]-\left[F_0\right]\in ker(H_n(F\times Z;\R)\to{\mathcal{H}}_n(F\times Z)).$$ 

To get a path-connected example, let $X$ be the space obtained by gluing one 
arc with the end points to the two different path components of $F\times Z$. This does not change the homology in degrees $\ge 2$ and thus
one has for $n\ge 2$: 
$$\left[F_1\right]-\left[F_0\right]\in ker(H_n(X;\R)\to{\mathcal{H}}_n(X)).$$ 

This space satisfies the other assumptions from \hyperref[lem2]{Lemma \ref*{lem2}}, 
but there is no finite covering by contractible,
measurable sets, though a covering by two contractible, non-measurable sets exists. 
\end{expl}

\subsection{A simplicial construction: straightening}\label{4.4-5.4}
Recall that for a topological space $X$ and a point $x_0\in X$ we denote $$S_k^{x_0}(X)=\left\{\sigma\colon\Delta^k\to X\mid \sigma(v_0)=\sigma(v_1)=\ldots=\sigma(v_k)=x_0\right\},$$
where $v_0,v_1,\ldots,v_k$ denote the vertices of the standard simplex $\Delta^k$. Two simplices $\sigma_0,\sigma_1\in
S_k^{x_0}(X)$ 
are 
said to be homotopic rel.\ boundary if there exists a continuous map $F\colon \Delta^k\times\left[0,1\right]$ with $F(x,0)=\sigma_0(x),F(x,1)=\sigma_1(x)$
for all $x\in\Delta^k$ and $F(x,t)=x_0$ for all $x\in\partial\Delta^k, t\in\left[0,1\right]$.

Let us denote by $C_*^{x_0}(X)\subset C_*(X)$ for some fixed $x_0\in X$ the subcomplex
generated by 
$S_*^{x_0}(X)$. In the following lemma we define a topological analogue of the well-known geometric straightening which we used in the proof
of \hyperref[negcurvinj]{Corollary \ref*{negcurvinj}}. The construction replaces geodesics and straight simplices by a somewhat arbitrary selection of simplices. (Similar constructions in somewhat different settings can be found in \cite[Theorem 9.5]{may} and
\cite[Proposition 3.1]{bfp}.)

\begin{lemma}\label{straight}Let $X$ be a topological space and $x_0\in X$.

i) There is a subset $S_*^{str}(X)\subset S_*^{x_0}(X)$ such that $S_k^{str}(X)$ contains one $k$-simplex in each homotopy class rel.\ boundary of simplices with all boundary faces in $S_{k-1}^{str}(X)$.

ii)  There is a chain map 
$$str\colon C_*^{x_0}(X)\to C_*^{x_0}(X)$$ which is chain homotopic 
to the identity and whose image lies in the chain complex
$C_*^{str}(X)\subset C_*^{x_0}(X)$ spanned by the simplices in $S_*^{str}(X)$. If $Y$ is a subspace of $X$ with $x_0\in Y$ such that $\pi_k(Y,x_0)\to\pi_k(X,x_0)$ is injective for all $k\ge 0$, then $str$ can be chosen to map $C_*^{x_0}(Y)$ to $C_*^{str}(Y)$.

iii) If $X$ is aspherical, then $str$ is constant on homotopy classes rel.\ vertices.

iv) If $X$ is a
space having the topology of a 
countable CW-complex, then $str$ extends to a chain map
$$str\colon {\mathcal C}_*^{x_0}(X)\to  C_*^{str}(X).$$ 

\end{lemma}

\begin{pf}
We recursively construct a subcomplex $C_*^{str}(X)\subset C_*^{x_0}(X)$, whose simplices we call the ``straight simplices", and a map $str\colon C_*^{x_0}(X)\to C_*^{str}(X)$, which we call the straightening map. 
This is similar to well-known constructions, which in slightly different settings can be found in \cite[Theorem 9.5]{may} and
\cite[Proposition 3.1]{bfp}. We recall the construction for convenience of the reader and for later reference.

The $0$-skeleton of $C_*^{str}(X)$ consists of the one 
vertex $x_0$. For the $1$-skeleton of $C_*^{str}(X)$ we choose one $1$-simplex in each homotopy class (rel.\ vertices) 
of $1$-simplices in $C_*^{x_0}(X)$, 
and we define $str$ on $1$-simplices by sending each of them to the unique straight simplices in its homotopy class. For each $1$-simplex $\sigma$
we fix a homotopy (rel.\ vertices) between $\sigma$ and $str(\sigma)$.

Assume now that for some $k>1$ we have already defined $S_{*\le k-1}^{str}(X)$ 
and $str\colon C_{*\le k-1}^{x_0}(X)\to C_{*\le k-1}^{str}(X)$. For $S_k^{str}(X)$
we choose one $k$-simplex with all boundary faces in $S_{k-1}^{str}(X)$ inside 
each homotopy class (rel.\ boundary) of $k$-simplices in $S_k^{x_0}(X)$ with all boundary faces in $S_{k-1}^{str}(X)$.

For a simplex $\sigma\in C_k^{x_0}(X)$ we can assume by induction that we have defined $str(\partial\sigma)$
and that we have a homotopy between $\partial\sigma$ and $str(\partial\sigma)$. By the cofibration
property of the inclusion $\partial\Delta^k\to\Delta^k$ this homotopy extends to a homotopy of $\sigma$ keeping vertices fixed. Let $\sigma^\prime$ be the result of this homotopy.
Among simplices with boundary $str(\partial\sigma)$ we have in the homotopy class (rel.\ boundary) of $\sigma^\prime$ exactly one simplex in $S_k^{str}(X)$. Define this simplex to be $str(\sigma)$. 
By construction we have a homotopy from $\sigma$ to $str(\sigma)$ whose restriction to $\partial\sigma$ is a reparametrisation of the homotopy from $\partial\sigma$ to $str(\partial\sigma)$, which was given by the inductive hypothesis.
This family of compatible homotopies yields the wanted chain homotopy between $id$ and $str$.

For a pair $(X,Y)$ one chooses straight simplices to be in $Y$ whenever this is possible. The assumption on injectivity of $\pi_k(Y,x_0)\to\pi_k(X,x_0)$ implies that homotopies between simplices in $Y$ can be chosen to remain in $Y$. Inductively, for a simplex $\sigma$ in $Y$ this applies (in the above recursive construction) to the homotopy between $str(\partial\sigma)$ and $\partial\sigma$, and then the extending homotopy yields a simplex in $Y$.

For  iii), if $\sigma$ and $\sigma^\prime$ are homotopic rel.\ vertices, then the 1-skeleta of $str(\sigma)$ and $str(\sigma^\prime)$ agree. Assuming inductively that for some $k\ge 2$ the $(k-1)$-skeleta of $str(\sigma)$ 
and $str(\sigma^\prime)$ agree, we get from asphericity of $X$ that the $k$-skeleta of both simplices are homotopic rel.\ boundary, and thus these $k$-skeleta must be equal by property (i). Proof by induction on $k$ 
yields that $str(\sigma)=str(\sigma^\prime)$.

For iv), given $\mu\in{\mathcal C}_*^{x_0}(X)$, we note that for a Borel set $A\subset map(\Delta^k,X)$, its preimage $str^{-1}(A)$ is a countable union of homotopy classes, which are Borel sets,  
as discussed in the introductory paragraph of 
\hyperref[hypex]{Section \ref*{hypex}}.

Thus we can define $str(\mu)$ by $$str(\mu)(A)=\mu(str^{-1}(A)).$$
Its variation is bounded because $$\Vert str(\mu)\Vert =max_A\mu(str^{-1}(A))-min_B\mu(str^{-1}(B))\le \Vert\mu\Vert<\infty.$$
Local contractibility of $X$ implies by \cite{wada} that $map(\Delta^k,X)$ is locally path-connected. In fact, \cite{wada} has a more precise result about subsets of $map(\Delta^k,X)$ mapping subpolyhedra to given closed subspaces, from which we can in particular infer that the subset of $map(\Delta^k,X)$, consisting of simplices mapping all vertices to $x_0$, is locally path-connected. By a theorem of Fox, path components correspond to homotopy classes. This means that every compact set $K\subset map(\Delta^k,X)$ meets only finitely many homotopy classes. So, if $\mu$ has compact determination set $K\subset map(\Delta^k,X)$, then $str(K)$ is a finite set by (iii).
Since $str(\mu)$ is determined on this finite set, it is a finite sum, i.e., an element of $C_*^{str}(X)$. Hence the image of $str$ does indeed belong to the subspace $C_*^{str}(X)$ consisting of singular (rather than measure) chains.
\end{pf}

Of course $str$ need in general not be continuous or measurable.

\subsection{Measurable sections}

We remark that in the setting of CW-complexes \hyperref[pmh]{Section \ref*{pmh}} can be used to give a simpler proof than the one given in \cite{loe} by reducing the problem to the simpler situation of simplices having all their vertices in the basepoint.

In fact, in the setting of countable fundamental groups, in particular for CW-complexes, we can also prove the existence of a measurable section, which was the main technical lemma in \cite{loe}.
Namely, let \hl{$G\cong\pi_1(X,x_0)$} be the (by assumption countable) deck transformation group of the generalized universal covering $p\colon\widetilde{X}\to X$. The lift of $x_0$ to $\widetilde{X}$ is a $G$-orbit $G\tilde{x}_0$ for some $\tilde{x}_0$,
see \hyperref[guc]{Section \ref*{guc}}. 
The lift of a homotopy class (rel.\ $x_0$) of $n$-simplices is a homotopy class inside $S_*^{\gamma_0\tilde{x}_0,\ldots,\gamma_n\tilde{x}_0}(\widetilde{X})$, by which we mean the set of simplices mapping their $i$-th vertex to $\gamma_i\tilde{x}_0$ for $i=0,\ldots,n$. Clearly the projection 
$$p\colon S_*^{\gamma_0\tilde{x}_0,\ldots,\gamma_n\tilde{x}_0}(\widetilde{X})\to S_*^{x_0}(X)$$
maps a homotopy class homeomorphically onto its image, which is also a homotopy class. In particular, the restriction of $p$ to any homotopy class has a continuous right-inverse $s$ defined on the image of that homotopy class. Thus we get a right-inverse $s$ defined on each of the homotopy classes downstairs. 
Since the homotopy classes of simplices (rel.\ vertices) are Borel sets 
(cf.\ \hyperref[hypex]{Section \ref*{hypex}}), and   
there are only countably many homotopy classes, the so-defined $s$ yields a measurable map $$s\colon S_*^{x_0}(X)\to S_*^{G\tilde{x}_0}(\widetilde{X})$$
right-inverse to $p$.

The following example shows that no such section can exist when the fundamental group is uncountable,
like for the earring space.  
\begin{thm}\label{uncount}Let $X$ be a complete, separable metric space admitting a generalized universal covering $p\colon\widetilde{X}\to X$. Assume that $X$ is semi-locally simply connected at $x_0$ and that the deck group \hl{$G\cong\pi_1(X,x_0)$} is uncountable. Then $p_1\colon S_1^{G\tilde{x}_0}\to S_1^{x_0}$ admits no Borel section.\end{thm}

\begin{pf}$X$ being semi-locally simply connected at $x_0$ implies that $G\tilde{x}_0$ is discrete in $\widetilde{X}$. 

Then one can show that the sets $S_1^{g_1\tilde{x}_0,g_2\tilde{x}_0}$ are open and closed, so that we have decomposed $S_1^{G\tilde{x}_0}$ into a disjoint union of open and closed sets.

$G$ is uncountable, hence the cardinality of its power set is bigger than continuum. For a measurable section $s$ and any subset $\Gamma\subset G$ let $U_\Gamma$ be the set of 1-simplices with the same endpoints as $s(\gamma_g)$
for $g\in\Gamma$, where $\gamma_g$ means the loop representing $g$. 

As a disjoint union of sets of the form $S_1^{g_1\tilde{x}_0,g_2\tilde{x}_0}$, any $U_\Gamma$ must be open and hence have Borel-measurable preimage under $s$. Thus the $s^{-1}(U_\Gamma)$ are more than continuum many distinct Borel-measurable subsets of $map(\Delta^1,X)$. 

But the latter is a Polish space by \cite[Theorem 4.19]{kec} and thus has only continuum many Borel sets by \cite[Theorem 3.3.18]{sri}. This is a contradiction.
\end{pf}

\end{document}